\documentclass{gtpart}     
\usepackage[all]{xy}
\usepackage{amscd}
\usepackage[usenames]{color}
\usepackage{enumitem}
\usepackage{hyperref}
\usepackage{amsthm}
\usepackage{graphicx}
\title{Two-generator free Kleinian groups and\\ hyperbolic displacements}
\author{\'{I}lker S. Y\"{u}ce}
\givenname{\.{I}lker Sava\c{s}} \surname{Y\"{u}ce}
\address{TED University\\
Ziya G\"okalp St., No. 48, Kolej 06420, \c{C}ankaya, ANKARA, TURKEY}
\email{} \urladdr{}
\newtheorem{theorem}{Theorem}[section]
\newtheorem{lemma}[theorem]{Lemma}
\newtheorem{proposition}[theorem]{Proposition}
\theoremstyle{definition}
\newtheorem{definition}[theorem]{Definition}

\theoremstyle{remark}

\newtheorem*{thm}{Theorem}     
\begin{document}
\newcommand{\noi}{\noindent}
\newcommand{\tb}{\textbf}
\newcommand{\tr}{\textrm}
\newcommand{\pr}{\partial}
\newcommand{\fr}{\frac}
\newcommand{\ra}{\rightarrow}
\newcommand{\Ra}{\Rightarrow}
\newcommand{\lag}{\langle}
\newcommand{\rag}{\rangle}
\newcommand{\vs}{\vspace{.1in}}
\newcommand{\vvs}{\vspace{.05in}}
\newcommand{\ep}{\epsilon}
\newcommand{\de}{\delta}
\newcommand{\dis}{\displaystyle}
\newcommand{\dgammaa}{\textrm{dist}(z_1,\ \gamma\cdot z_1)}                
\newcommand{\dgammaz}{\textrm{dist}(z,\ \gamma\cdot z)}                 
\newcommand{\dgamma}{\textrm{dist}(z_0,\ \gamma\cdot z_0)}                 
\newcommand{\dgammainv}{\textrm{dist}(z_0,\ \gamma^{-1}\cdot z_0)}         
\newcommand{\dgammai}{\textrm{dist}(z_0,\ \gamma_i\cdot z_0)}              
\newcommand{\dgammaj}{\textrm{dist}(z_0,\ \gamma_j\cdot z_0)}              
\newcommand{\dzetak}{\textrm{dist}(z_0,\ \zeta_k\cdot z_0)}                
\newcommand{\ddeltar}{\textrm{dist}(z_0,\ \delta_r\cdot z_0)}              
\newcommand{\dgammaiprm}{\textrm{dist}(z_0,\ \gamma_i'\cdot z_0)}          
\newcommand{\dxi}{\textrm{dist}(z_0,\ \xi\cdot z_0)}                       
\newcommand{\dxiprm}{\textrm{dist}(z_0,\ \xi'\cdot z_0)}                   
\newcommand{\detaprm}{\textrm{dist}(z_0,\ \eta'\cdot z_0)}                 
\newcommand{\dxizone}{\textrm{dist}(z_1,\ \xi\cdot z_1)}                   %
\newcommand{\dpsi}{\textrm{dist}(z,\ \psi\cdot z)}                         
\newcommand{\dpsiz}{\textrm{dist}(z_0,\ \psi\cdot z_0)}                      
\newcommand{\dpsizinv}{\textrm{dist}(z_0,\ \psi_0^{-1}\cdot z_0)}              
\newcommand{\dphi}{\textrm{dist}(z,\ \phi\cdot z)}                         
\newcommand{\dxiinv}{\textrm{dist}(z_0,\ \xi^{-1}\cdot z_0)}               
\newcommand{\deta}{\textrm{dist}(z_0,\ \eta\cdot z_0)}                     
\newcommand{\detazone}{\textrm{dist}(z_1,\ \eta\cdot z_1)}                 %
\newcommand{\detainv}{\textrm{dist}(z_0,\ \eta^{-1}\cdot z_0)}             
\newcommand{\dxieta}{\textrm{dist}(z_0,\ \xi\eta\cdot z_0)}                
\newcommand{\dxietaprm}{\textrm{dist}(z_0,\ \xi'\eta'\cdot z_0)}           
\newcommand{\dpsiphi}{\textrm{dist}(z_0,\ \psi\phi\cdot z)}                
\newcommand{\detaxi}{\textrm{dist}(z_0,\ \eta^{-1}\xi^{-1}\cdot z_0)}      
\newcommand{\detantwo}{\textrm{dist}(z_0,\ \eta^{-2}\cdot z_0)}            
\newcommand{\detanonexi}{\textrm{dist}(z_0,\ \eta^{-1}\xi\cdot z)}         
\newcommand{\dxitwo}{\textrm{dist}(z_0,\ \xi^{2}\cdot z_0)}                
\newcommand{\dxintwo}{\textrm{dist}(z_0,\ \xi^{-2}\cdot z_0)}              
\newcommand{\dxietanone}{\textrm{dist}(z_0,\ \xi\eta^{-1}\cdot z_0)}       
\newcommand{\detaxinone}{\textrm{dist}(z_0,\ \eta\xi^{-1}\cdot z_0)}       
\newcommand{\disetaxi}{\textrm{dist}(z_0,\ \eta\xi\cdot z_0)}              
\newcommand{\disetatwo}{\textrm{dist}(z_0,\ \eta^2\cdot z_0)}              
\newcommand{\disxiinvetainv}{\textrm{dist}(z_0,\ \xi^{-1}\eta^{-1}\cdot z_0)}  
\newcommand{\disxiinveta}{\textrm{dist}(z_0,\ \xi^{-1}\eta\cdot z_0)}          
\newcommand{\disbetaij}{\textrm{dist}(z_0,\ \beta_{i,j}\cdot z_0)}             
\newcommand{\disgammaij}{\textrm{dist}(z_0,\ \gamma_{i,j}\cdot z_0)}           
\newcommand{\diszetaij}{\textrm{dist}(z_0,\ \zeta_{i,k}\cdot z_0)}             
\newcommand{\disupsilonij}{\textrm{dist}(z_0,\ \upsilon_{i,j}\cdot z_0)}       
\newcommand{\hyp}{\mathbb{H}^3}
\newcommand{\hype}{\overline{\mathbb{H}}^3}
\newcommand{\alphas}{\alpha_*}
\newcommand{\betas}{\beta_*}
\newcommand{\dseven}{\Delta^7}
\newcommand{\deleven}{\Delta^{11}}
\newcommand{\tnr}{\textnormal}
\newcommand{\xb}{\textbf{x}}
\newcommand{\xs}{\textbf{x}^*}
\newcommand{\xc}{\textbf{x}^{\circ}}
\newcommand{\ys}{\textbf{y}^*}
\newcommand{\yc}{\textbf{y}^{\circ}}
\newcommand{\triri}{\triangleright}
\def\co{\colon\thinspace}

\begin{abstract}    
The $\log 3$ theorem, proved by Culler and Shalen, states that every point in the hyperbolic $3$--space $\hyp$ is moved a distance at least $\log 3$ by
one of the non--commuting isometries $\xi$ or $\eta$ of  $\hyp$ provided that $\xi$ and $\eta$ generate a torsion--free, discrete group which is not
co-compact and contains no parabolic.
This theorem lies in the foundation of many techniques that provide lower estimates for the volumes of orientable,
closed hyperbolic 3--manifolds whose fundamental groups have no $2$--generator subgroup of finite index and, as a consequence, gives insights into the
topological properties of these manifolds.

Under the hypotheses of the $\log 3$ Theorem, the main result of this paper shows that every point in $\hyp$ is moved a distance at least $\log\sqrt{5+3\sqrt{2}}$ by one of the isometries $\xi$, $\eta$ or $\xi\eta$.
\end{abstract}

\maketitle



\section{Introduction}

Let $M$ be a closed orientable hyperbolic $3$--manifold. Anderson, Canary, Culler and Shalen prove in \cite{ACCS} that $3.08$ is a lower bound for the volume of $M$ under the assumptions that the first Betti number of $M$ is at least $4$ and $\pi_1(M)$ has no subgroup isomorphic to the fundamental group of a genus two surface. In \cite{CSParadox}, Culler and Shalen show that the volume of $M$ is at least $0.92$ provided that the first Betti number of $M$ is at least $3$ and $\pi_1(M)$ has no two--generator subgroup of finite index. Later Culler, Hersonsky and Shalen improve the previous volume estimate to $0.94$ in \cite{CHS}. These deep results are among a number of theorems stated in  \cite{ACS}, \cite{ACCS}, \cite{CHS}, \cite{CSParadox}, and \cite{CSMargulis} that relate the topology of hyperbolic $3$--manifolds to their geometry.

The common denominator in all of the volume estimates listed above is that they are consequences of one of the fundamental results in the study of Kleinian groups, the so called $\log 3$ theorem proved by Culler and Shalen \cite{CSParadox} and its generalization due to Anderson, Canary, Culler and Shalen \cite{ACCS}. This seminal result can be stated as follows:

\textit{Let $\xi$ and $\eta$ be non--commuting isometries of $\hyp$. Suppose that $\xi$ and $\eta$ generate a torsion--free, topologically tame,  discrete group which is not
co--compact and contains no parabolic. Let $\Gamma_1$ and $\alpha_1$ denote the set of isometries $\{\xi,\eta\}$ and the real number $9$, respectively. Then, for any $z_0\in\hyp$, we have}
\begin{displaymath}
e^{\left(\dis{2\max\nolimits_{\gamma\in\Gamma_1}\left\{\dgamma\right\}}\right)}\geq\alpha_1.
\end{displaymath}
The $\log 3$ theorem and its generalization imply that $(1/2)\log 5$ and $(1/2)\log 3$ are Margulis numbers for the hyperbolic $3$--manifolds which satisfy the conditions in the cases for which the first Betti numbers are at least $4$ or $3$, respectively. Consequently, the lower bounds computed in \cite{ACCS}, \cite{CSParadox}, and \cite{CHS} for the volumes of such manifolds follow. Although the bounds given in \cite{CSParadox} and \cite{ACCS} are superseded by the recent works of Gabai--Meyerhoff--Milley \cite{DRP}, \cite{DRP2} and Milley \cite{Mil} using a newer approach, Mom technology, it is conceivable that an improvement in the lower bound for the displacements under the isometries described in the $\log 3$ theorem will lead to improved Margulis numbers and lower bounds for the volumes of the classes of hyperbolic $3$--manifolds mentioned above through the ideas introduced in \cite{ACCS}, \cite{CSParadox}, and \cite{CHS}. With this motivation, in this paper, we prove the following:

\textbf{Main Result.}\textit{ Let $\xi$ and $\eta$ be non--commuting isometries of $\hyp$. Suppose that $\xi$ and $\eta$ generate a torsion--free discrete group which is not
co--compact and contains no parabolic. Let $\Gamma_{\dagger}$ and $\alpha_{\dagger}$ denote the set of isometries $\{\xi,\eta,\xi\eta\}$ and the real number $5+3\sqrt{2}$, respectively. Then, for any $z_0\in\hyp$, we have}
$$
e^{\left(\dis{2\max\nolimits_{\gamma\in\Gamma_{\dagger}}\left\{\dgamma\right\}}\right)}\geq\alpha_{\dagger},
$$
which is given as Theorem \ref{thm5.1} in Section \ref{sec5}. 

An orientable hyperbolic $3$--manifold may be regarded as the quotient of the hyperbolic $3$--space $\hyp$ by a discrete group $\Gamma$ of orientation--preserving isometries of $\hyp$. If $\Gamma$ is a torsion free Kleinian group and $M=\hyp/\Gamma$, then $\Gamma$ is called
\textit{topologically tame} if $M$ is homeomorphic to the interior of a compact $3$--manifold. In \cite{Agol} and \cite{CG}, Agol and Calegari--Gabai prove that every finitely generated Kleinian group is topologically tame. Therefore, we drop the tameness hypothesis from Theorem \ref{thm5.1}.

The proof of Theorem \ref{thm5.1} requires the use of the same ingredients introduced in \cite{CSParadox} to prove the $\log 3$ Theorem. In the following subsections
of Introduction, we review these ingredients  briefly. In particular, we summarize the proof of the $\log 3$ Theorem in \S\ref{sec1.1} with an emphasis on the calculations required to obtain the number $\log 3$. In \S\ref{sec1.2}, we propose an alternative technique to perform these calculations which makes it possible to extend Culler and Shalen's arguments in \cite{CSParadox} to determine a lower bound for the displacements under any given set of isometries in $\Gamma=\langle\xi,\eta\rangle$ as long as the hypotheses of the $\log 3$ theorem are satisfied. We describe this extension and summarize its application to the set $\Gamma_{\dagger}=\{\xi,\eta,\xi\eta\}\subset\Gamma$ to achieve the lower bound stated in Theorem \ref{thm5.1} in \S\ref{sec1.3}.

In the rest of this manuscript the boundary of the canonical compactification $\hype$ of $\hyp$ will be denoted by $S_{\infty}$, which is homeomorphic to $S^2$. The notation $\Lambda_{\Gamma\cdot z}$ will denote the limit set of $\Gamma$--orbit of $z\in\hyp$ on $S_{\infty}$. By dist($z,\gamma\cdot z$) we will mean the hyperbolic displacement of $z\in\hyp$ under the action of the isometry $\gamma\co\hyp\to\hyp$.  Any isometry $\gamma$ of $\hyp$ extends to a conformal
automorphism $\overline{\gamma}\co\hype\to\hype$. The conformal automorphism of $S_{\infty}$ obtained by restricting
$\overline{\gamma}$ will be denoted by $\gamma_{\infty}$.

The author would like to extend his sincerest thanks to the anonymous referee whose recommendations lead to a much better exposition of the ideas in this paper, shortened the proofs substantially and made this text much more readable as a result. He is deeply grateful to Peter B. Shalen for setting the course of this research and very helpful discussions. He is also grateful to Marc E. Culler for his corrections in an earlier version of this work.


\subsection{A decomposition of $\Gamma=\langle\xi,\eta\rangle$ and Proof of the $\log 3$ theorem}\label{sec1.1}

Let $\xi$ and $\eta$ be two non-commuting isometries of $\hyp$. Suppose that $\xi$ and $\eta$ generate a torsion--free discrete group which is not co--compact and contains no parabolic. Then $\Gamma=\langle\xi,\eta\rangle$ is a free group of rank $2$ (\cite{CSParadox}, Proposition 9.2). This fact allows one to decompose $\Gamma$ as disjoint union of subsets of reduced words. In particular, 
the decomposition
\begin{equation}\label{prdx.d.}
\Gamma=\{1\}\cup\bigcup_{\psi\in\Psi^1}J_{\psi}
\end{equation}
is used in the proof of the $\log 3$ theorem, which is carried out in two cases:
\begin{enumerate}[label=\roman*]
    \renewcommand{\labelenumi}{\roman{enumi}}
       \item\hspace{-.3cm}.\label{I}  when $\Gamma$ is geometrically infinite; that is, $\Lambda_{\Gamma\cdot z}=S_{\infty}$ for every $z\in\hyp$ and,
       \item\hspace{-.3cm}.\label{II} when $\Gamma$ is geometrically finite.
\end{enumerate}
In (\ref{prdx.d.}) each $J_{\psi}$ is defined as the set of all non--trivial reduced words in $\Gamma$ that have the initial letter $\psi\in\Psi^{1}=\{\xi,\eta,\eta^{-1},\xi^{-1}\}$.

In the case (\ref{I}), Culler and Shalen first prove that the Patterson density, a $\Gamma$--invariant conformal density $(\mu_z)_{z\in\hyp}$, constructed by Patterson (\cite{SJP}) and extensively studied by Sullivan (\cite{Su1}, \cite{Su2}, \cite{Su3}), is the area density $(A_z)_{z\in\hyp}$, whose support is $S_{\infty}$ (\cite{CSParadox}, Propositions 6.9 and 3.9). Then, using the decomposition  (\ref{prdx.d.}) together with its group--theoretical properties
\begin{eqnarray}\label{prdxgp}
\psi J_{\psi^{-1}} & = & \Gamma - J_{\psi}
\end{eqnarray}
for $\psi\in\Psi^1$, they construct a decomposition of the area density $(A_z)_{z\in\hyp}$, which in turn gives a decomposition of the area measure $A_{z_0}$ based at $z_0\in\hyp$ into a finite sum of four measures $\nu_{\xi}$, $\nu_{\eta}$, $\nu_{\eta^{-1}}$, $\nu_{\xi^{-1}}$  so that each measure $\nu_{\psi^{-1}}$ is transformed to the complement of $\nu_{\psi}$ for $\psi\in\Psi^{1}$ (\cite{CSParadox}, Proposition 4.2 (ii) and Lemma 5.3 (ii) and (iii)). In other words, they obtain the following:
\begin{theorem}\label{thm1.1}
Let $\Gamma=\langle\xi,\eta\rangle$ be a free, geometrically infinite Kleinian group without parabolics. For any $z_0\in\hyp$, let $A_{z_0}$ be the area measure based at $z_0$. There is a family of Borel measures $\{\nu_{\psi}\}_{\psi\in\Psi^1}$ on $S_{\infty}$ for $\Psi^1=\{\xi,\eta,\eta^{-1},\xi^{-1}\}$ such that
\begin{itemize}
 \item[(1)] $A_{z_0}(S_{\infty})=\sum_{\psi\in\Psi^1}\nu_{\psi}(S_{\infty})$, where $A_{z_0}$ is normalized so that $A_{z_0}(S_{\infty})=1$, and,

\item[(2)]  $\dis{\int_{S_{\infty}}\left(\lambda_{\psi,z_0}\right)^2d\nu_{\psi^{-1}}=1-\int_{S_{\infty}} d\nu_{\psi}}$ for each $\psi\in\Psi^1$.
\end{itemize}
Furthermore, if $z_0$ is on the common perpendicular $\ell(\xi,\eta)$ of the isometries
$\xi$ and $\eta$, then
\begin{itemize}
\item[(3)]  $\dis{\int_{S_{\infty}} d\nu_{\xi^{-1}}=\int_{S_{\infty}} d\nu_{\xi}\ \ \tr{and} \ \ \int_{S_{\infty}} d\nu_{\eta^{-1}}=\int_{S_{\infty}} d\nu_{\eta}}.$
\end{itemize}
\end{theorem}
Theorem \ref{thm1.1} is not explicitly stated in \cite{CSParadox}. But, as summarized above, it follows from Lemma 5.3 using the conclusions of Propositions 4.2, 6.9 and 3.9 in \cite{CSParadox}. The function $\lambda_{\psi,z_0}$ in part (2) is the conformal expansion factor of $\psi_{\infty}$ measured in the round metric centered at $z_0$ (see \cite{CSParadox}, \S2.4 for details). The common perpendicular $\ell(\xi,\eta)$ mentioned in Theorem \ref{thm1.1} is the fixed locus of the involution $\tau\in\tnr{Isom}^+(\hyp)$ that conjugates $\xi$ to $\xi^{-1}$ and $\eta$ to $\eta^{-1}$ when $\xi_{\infty}$ and $\eta_{\infty}$ have no common fixed point on $S_{\infty}$ (\cite{CSParadox}, Proposition 1.8).

Next, they relate the masses of the measures $\nu_{\xi}$, $\nu_{\eta}$, $\nu_{\eta^{-1}}$, and $\nu_{\xi^{-1}}$ to the hyperbolic displacements $\dxi$, $\deta$, $\detainv$, and $\dxiinv$. In particular, they prove the statement below:
\begin{lemma}\label{lem1.2}(\cite{CSParadox}, Lemma 5.5; \cite{CSMargulis}, Lemma 2.1)\label{CSlem5.5} Let $a$ and $b$ be numbers in $[0,1]$ which are not both equal to $0$ and are not both equal to $1$. Let $\gamma$ be a loxodromic isometry of $\hyp$ and let $z_0$
be a point in $\hyp$. Suppose that $\nu$ is a measure on $S_{\infty}$ such that
(i) $\nu\leq A_{z_0}$,  (ii) $\nu\left(S_{\infty}\right)\leq a$, (iii) $\int_{S_{\infty}}(\lambda_{\gamma,z_0})^2d\nu\geq b$.
Then we have $a>0$, $b<1$, and \begin{equation*} \textnormal{dist}(z_0,\gamma\cdot z_0)\geq\frac{1}{2}\log\frac{b(1-a)}{a(1-b)}.\end{equation*}
\end{lemma}
Then, using Theorem \ref{thm1.1} and Lemma \ref{CSlem5.5}, they calculate the lower bound $\log 3$ when $\Gamma$ is geometrically infinite as follows (\cite{CSParadox}, Proposition 5.2): By the geometric fact
$\max_{\gamma\in\{\xi^{-1},\eta^{-1}\}}\{\dgamma\}\geq\max_{\gamma\in\{\xi^{-1},\eta^{-1}\}}\{\dgammaa\}$
for any $z_0\in\hyp$ and the orthogonal projection $z_1$ of $z_0$ onto $\ell(\xi,\eta)$, we may assume, without loss of generality, that $z_0\in\ell(\xi,\eta)$.  Theorem \ref{thm1.1} (1) and (3) imply $1=2\nu_{\xi}(S_{\infty})+2\nu_{\eta}(S_{\infty})$, which in turn shows either $\nu_{\xi}(S_{\infty})\leq 1/4$
or $\nu_{\eta}(S_{\infty})\leq 1/4$. If $\nu_{\xi}(S_{\infty})\leq 1/4$, parts (2) and
(3) of Theorem \ref{thm1.1} give that
$\int_{S_{\infty}}(\lambda_{\xi^{-1},z_0})^2d\nu_{\xi}=1-\nu_{\xi}\geq 3/4$.
Since $\nu_{\xi}\leq A_{z_0}$ by Theorem \ref{thm1.1} (1), it is possible to apply Lemma \ref{CSlem5.5} with the choices $a=1/4$, $b=3/4$, $\nu=\nu_{\xi}$ and
$\gamma=\xi^{-1}$ to conclude that
$\dxiinv \geq (\log 9)/2$.
If $\nu_{\eta}(S_{\infty})\leq 1/4$, an analogous calculation with the choices $a=1/4$, $b=3/4$, $\nu=\nu_{\eta}$ and $\gamma=\eta^{-1}$ shows that
$\detainv\geq(\log 9)/2$. Because $\dxi=\dxiinv$ and $\deta=\detainv$, the $\log 3$ theorem follows in the case (\ref{I}).

In the case (\ref{II}), Culler and Shalen define the function
\[
f^1_{z_0}(\xi,\eta)=\max\{\dxi,\deta\}
\]
for a fixed $z_0\in\hyp$ on the character variety $\mathfrak{X}=\tnr{Isom}^+(\hyp)\times\tnr{Isom}^+(\hyp)$ of the free group on two generators. It is easy to show that $f_{z_0}$ is proper and continuous on the closure of
the set $\mathfrak{GF}$ consisting of $(\xi,\eta)$ such that $\langle\xi,\eta\rangle$ is free on the generators $\xi$ and $\eta$, geometrically finite and without any parabolic. They prove that $f^1_{z_0}$ has no minimum in $\mathfrak{GF}$ (\cite{CSParadox}, Lemma 9.5). Since $\mathfrak{GF}$ is open in $\mathfrak{X}$ (\cite{M}, Theorem 8.1), the function $f^1_{z_0}$ achieves it minimum on the boundary $\overline{\mathfrak{GF}}-\mathfrak{GF}$. Then they show that the set of $(\xi,\eta)$ such that $\langle\xi,\eta\rangle$ is free, geometrically infinite and without any parabolic is dense in $\overline{\mathfrak{GF}}-\mathfrak{GF}$ (\cite{CSParadox}, Proposition 8.2). Therefore, the fact that every $(\xi,\eta)\in\mathfrak{X}$ so that $\langle\xi,\eta\rangle$ is free and contains no parabolic is in $\overline{\mathfrak{GF}}$ (\cite{CSParadox}, Proposition 9.3) reduces geometrically finite case to geometrically infinite case, proving the $\log 3$ theorem.



\subsection{An alternative technique to calculate the number $\log 3$}\label{sec1.2}

When $\Gamma=\langle\xi,\eta\rangle$ mentioned in the $\log 3$ theorem is geometrically infinite, the calculation of the lower bound
$\log3$ in \S\ref{sec1.1} follows from the application of Lemma \ref{CSlem5.5} with the choices $\gamma=\psi^{-1}$, $\nu=\nu_{\psi}$, $a=1/4$ and $b=3/4$ for $\psi\in\Gamma_1=\{\xi,\eta\}$. These choices of $a$ and $b$ are derived from the equalities $\nu_{\xi}(S_{\infty})=\nu_{\xi^{-1}}(S_{\infty})$ and $\nu_{\eta}(S_{\infty})=\nu_{\eta^{-1}}(S_{\infty})$ in Theorem \ref{thm1.1} part (3), which requires the use of a point $z_0$ on the common perpendicular $\ell(\xi,\eta)$ of $\xi$ and $\eta$. Such a geodesic perpendicular to the axes of a collection of more than two isometries in $\Gamma=\langle\xi,\eta\rangle$ does not exist in general.

We calculate $\log 3$ without referring to $\ell(\xi,\eta)$ as follows: If $\nu_{\xi^{-1}}(S_{\infty})=0$, we get $\nu_{\xi}(S_{\infty})=1$ by Theorem \ref{thm1.1} (2). Then we obtain $\nu_{\eta}(S_{\infty})=0$ and $\nu_{\eta^{-1}}(S_{\infty})=0$ by Theorem \ref{thm1.1} (1). Part (2) of Theorem \ref{thm1.1} applied to $\psi=\eta^{-1}$ gives a contradiction.
Similar arguments for $\xi^{-1}$, $\eta$, $\eta^{-1}$ show that $0<\nu_{\psi}(S_{\infty})<1$ for $\psi\in\Psi^1$. In particular, we derive that
$0<\int_{S_{\infty}}\lambda_{\psi,z_0}^2d\nu_{\psi^{-1}}=1-\int_{S_{\infty}} d\nu_{\psi} < 1$
for every $\psi\in\Psi^1$ by Theorem \ref{thm1.1} (2). Hence we can apply Lemma \ref{CSlem5.5} with choices
\[
\begin{array}{llll}
\gamma=\xi, & \nu=\nu_{\xi^{-1}}, & a=\nu_{\xi^{-1}}(S_{\infty}),  & b=\int_{S_{\infty}}(\lambda_{\xi,z_0})^2d\nu_{\xi^{-1}},\\
\gamma=\eta, & \nu=\nu_{\eta^{-1}}, & a=\nu_{\eta^{-1}}(S_{\infty}), & b=\int_{S_{\infty}}(\lambda_{\eta,z_0})^2d\nu_{\eta^{-1}},\\
\gamma=\eta^{-1}, & \nu=\nu_{\eta}, & a=\nu_{\eta}(S_{\infty}), & b=\int_{S_{\infty}}(\lambda_{\eta^{-1},z_0})^2d\nu_{\eta},\\
\gamma=\xi^{-1}, & \nu=\nu_{\xi}, & a=\nu_{\xi}(S_{\infty}), & b=\int_{S_{\infty}}(\lambda_{\xi^{-1},z_0})^2d\nu_{\xi}.
\end{array}
\]
Then Lemma \ref{CSlem5.5} produces the inequalities
\begin{equation}
e^{2\dgamma}\geq\frac{\left(1-\int_{S_{\infty}} d\nu_{\gamma^{-1}}\right)\left(1-\int_{S_{\infty}} d\nu_{\gamma}\right)}{\int_{S_{\infty}} d\nu_{\gamma^{-1}}\cdot\int_{S_{\infty}} d\nu_{\gamma}}\label{8}
\end{equation}
for every $\gamma\in\Psi^1$. We consider the constants on the right-hand side of the inequalities in (\ref{8}) as the values of the functions
\begin{eqnarray*}
f_1^1(x_1,x_2,x_3,x_4)=\frac{1-x_4}{x_4}\cdot\frac{1-x_1}{x_1},&  & f_2^1(x_1,x_2,x_3,x_4)=\frac{1-x_3}{x_3}\cdot\frac{1-x_2}{x_2},\\
f_3^1(x_1,x_2,x_3,x_4)=\frac{1-x_2}{x_2}\cdot\frac{1-x_3}{x_3},&  & f_4^1(x_1,x_2,x_3,x_4)=\frac{1-x_1}{x_1}\cdot\frac{1-x_4}{x_4}
\end{eqnarray*}
at $\left(\int_{S_{\infty}} d\nu_{\xi},\int_{S_{\infty}} d\nu_{\eta},\int_{S_{\infty}} d\nu_{\eta^{-1}},\int_{S_{\infty}} d\nu_{\xi^{-1}}\right)\in\mathbb{R}^4$ with $\sum_{\psi\in\Psi^{1}}\int_{S_{\infty}}d\nu_{\psi}=1.$

Although the total masses of the measures $\nu_{\psi}$ for $\psi\in\Psi^1$ may change for a different $z_0\in\hyp$, Theorem \ref{thm1.1} part (1) still applies for the same decomposition in (\ref{prdx.d.}). Therefore, the calculation of the number $\log 3$ in the proof of the $\log 3$ theorem becomes a
consequence of the statement
\begin{equation*}
\inf\nolimits_{\tb{x}\in\Delta^3}\left\{\max\left(f_1^1(\tb{x}),f_2^1(\tb{x}),f_3^1(\tb{x}),f_4^1(\tb{x})\right)\right\}=9,
\end{equation*}
where $\Delta^3=\left\{(x_1,x_2,x_3,x_4)\in\mathbb{R}^4:x_1+x_2+x_3+x_4=1,\ x_i>0,\ 1\leq i\leq 4\right\}$ (see Lemma \ref{lem2.1}).
The functions $f_1^1$, $f_2^1$, $f_3^1$ and $f_4^1$ will be referred to as \textit{displacement functions} for the decomposition of $\Gamma$ in (\ref{prdx.d.}).

When combined with the arguments developed by Culler and Shalen in \cite{CSParadox} to prove the $\log 3$ theorem, this calculation technique extends to a process to find a lower bound for the maximum of the displacements under the isometries determined by a decomposition of $\Gamma=\langle\xi,\eta\rangle$. This extension is summarized in the next section.



\subsection{Decompositions of $\Gamma=\langle\xi,\eta\rangle$ and Hyperbolic displacements}\label{sec1.3}

For any non-elementary discrete subgroup $\Gamma$ of the isometries of the hyperbolic $n$--space $\mathbb{H}^n$, there exists a $\Gamma$--invariant conformal density $(\mu_z)_{z\in\mathbb{H}^n}$ for $\mathbb{H}^n$, constructed by Patterson \cite{SJP}, whose support is the limit set of $\Gamma$. In the case (\ref{I}) of the proof of the $\log 3$ theorem, it is required to decompose the Patterson density for $\hyp$, which is the area density $(A_z)_{z\in\hyp}$, whose support is $S_{\infty}$ (\cite{CSParadox}, Propositions 3.9 and 6.9), corresponding to the decomposition of $\Gamma=\langle\xi,\eta\rangle$ in (\ref{prdx.d.}). To this purpose, Culler and Shalen prove a more general statement, Proposition 4.2 in \cite{CSParadox}, which establishes the existence of a family of $D$--conformal densities $\left(\mathcal{M}_V\right)_{V\in\mathcal{V}}$, $D\in[0,n-1]$, for $\mathbb{H}^n$ indexed by a countable collection $\mathcal{V}$ of subsets of an infinite, uniformly discrete subset $W$ of $\mathbb{H}^n$ with certain conditions. One of which is that $\mathcal{M}_W$ can be decomposed as a sum of finitely many $D$--conformal densities $\mathcal{M}_{V_i}=(\mu_{V_i,z})_{z\in\hyp}$ provided that $W=\cup_{i=1}^mV_i\in\mathcal{V}$ for disjoint sets $V_i\in\mathcal{V}$ for $1\leq i\leq m$ (\cite{CSParadox}, Proposition 4.2 (ii)). In particular, Proposition 4.2 (ii) \cite{CSParadox} is applied to the disjoint union
\begin{equation*}
W^1=\{z_0\}\cup\bigcup_{\psi\in\Psi^1}\{\gamma\cdot z_0\co\gamma\in J_{\psi}\}\subset\hyp,
\end{equation*}
which leads to the decomposition of the area density $(A_z)_{z\in\hyp}$, and consequently, the decomposition of the area measure $A_{z_0}$ based at $z_0$ into a finite sum of Borel measures as stated in Theorem \ref{thm1.1} (1).

The notion of $D$--conformal densities, $D\in[0,n-1]$, for $\mathbb{H}^n$ is introduced by Sullivan (\cite{Su1}, \cite{Su2}, \cite{Su3}) as a generalization of $\Gamma$--invariant conformal densities for $\mathbb{H}^n$ (\cite{SJP}). Interested readers may refer to \cite{CSParadox}, \cite{Su1}, \cite{Su2}, \cite{Su3}, \cite{N}, and \cite{SJP} for the basics of this subject. In this paper, their use will be limited to the application of Proposition 4.2 \cite{CSParadox} to a carefully chosen infinite, uniformly discrete subset of $\hyp$ in \S\ref{Sec3} (see Lemma \ref{lem3.2}). Therefore, constructions and properties of conformal densities will be assumed without any exploration in the rest of this text. Unless otherwise is stated, we shall assume that $\Gamma=\langle\xi,\eta\rangle$ satisfies the hypotheses given in the $\log 3$ theorem.

The organization of the rest of this paper is as follows: Let $z_0$ be a fixed point in $\hyp$. In Section \ref{S2}, we give the necessary calculations in detail to obtain the number $\log 3$ for the $\log 3$
 theorem by using the approach outlined in \S\ref{sec1.2}. In particular, we show that the infimum of the maximum of $f_1^1(\tb{x})$ $f_2^1(\tb{x})$, $f_3^1(\tb{x})$, and $f_4^1(\tb{x})$ over the simplex $\Delta^3$ is $\alpha_1=9$ which follows from the facts
\begin{enumerate}
\item[a.]
$\inf_{\tb{x}\in\Delta^3}\left\{\max\left(f_1^1(\tb{x}),f_2^1(\tb{x})\right)\right\}=\min_{\tb{x}\in\Delta^3}\left\{\max\left(f_1^1(\tb{x}),f_2^1(\tb{x})\right)\right\}$,
\item[b.]\label{B} $\min_{\tb{x}\in\Delta^3}\left\{\max\left(f_1^1(\tb{x}),f_2^1(\tb{x})\right)\right\}=f_1^1(\tb{x}^*)$ for a  point
$\tb{x}^*\in\Delta_3\subset\Delta^3$,
\end{enumerate}
proved in Lemma \ref{lem2.1} in \S\ref{S2}, where $\Delta_3=\left\{\tb{x}\in\Delta^3:f_1^1(\tb{x})=f_2^1(\tb{x})\right\}$.

When $\Delta^3$ is considered as a submanifold of $\mathbb{R}^3$, the statement
$\tb{x}^*\in\Delta_3$ is deduced from the observation that $f_1^1$ and $f_2^1$ are smooth functions in an open neighborhood of $\Delta^3$ with no local minimums. The coordinates of $\tb{x}^*$ are calculated by using the conclusions that $\tb{x}^*$ is a solution of a certain Lagrange multipliers problem and satisfies the equality $f_1^1(\tb{x})=f_2^1(\tb{x})$. The lower bound $\log 3$ is obtained
by evaluating $f_1^1$ at the point $\tb{x}^*$.

Let $\Psi^{\dagger}=\{\xi\eta, \xi^2,
\xi\eta^{-1}, \eta, \xi^{-1}, \eta^{-1}\xi^{-1}, \eta^{-2}, \eta^{-1}\xi\}$ and $\Psi_r^{\dagger}=\{\xi,\eta^{-1}\}$. In \S\ref{Sec3}, we introduce the decomposition
\begin{equation}\label{decomp2}\Gamma_{\mathcal{D}^{\dagger}}=\{1\}\cup \Psi^{\dagger}_r\cup\bigcup\nolimits_{\psi\in\Psi^{\dagger}}J_{\psi}\end{equation}
of $\Gamma$ which contains the isometries in $\Gamma_{\dagger}=\{\xi,\eta,\xi\eta\}\subset\Psi^{\dagger}\cup\Psi^{\dagger}_r$. We apply Proposition 4.2 \cite{CSParadox} to the infinite, uniformly discrete subset
\begin{equation}\label{decomp2LS}
W^{\dagger}=\{z_0\}\cup\{\gamma\cdot z_0\co\gamma\in\Psi^{\dagger}_r\}\cup\bigcup_{\psi\in\Psi^{\dagger}}\{\gamma\cdot z_0\co\gamma\in J_{\psi}\}
\end{equation}
of $\hyp$ to construct the decomposition of the area measure $A_{z_0}$ based at $z_0$ corresponding the decomposition
$\Gamma_{\mathcal{D}^{\dagger}}$ in Theorem \ref{thm3.4}, an analog of Theorem \ref{thm1.1} for $\Psi^{\dagger}$, in \S\ref{Sec3}. Using Theorem \ref{thm3.4} and Lemma \ref{CSlem5.5}, we determine the displacement functions for the
displacements under the isometries in $\Psi^{\dagger}\cup\Psi^{\dagger}_r$ in \S\ref{Sec3}. There are $18$ displacement functions  $\{f_i\}_{i=1}^8$, $\{g_j\}_{j=1}^6$ and $\{h_k\}_{k=1}^4$ for $\Gamma_{\mathcal{D}^{\dagger}}$ with formulas
\begin{displaymath}
\begin{array}{l}
\dis{f_1(\tb{x})=\frac{1-x_4-x_5-x_6}{x_4+x_5+x_6}\cdot\frac{1-x_1}{x_1}},\  \dis{f_2(\tb{x})=\frac{1-x_4-x_5-x_6-x_7-x_8}{x_4+x_5+x_6+x_7+x_8}\cdot\frac{1-x_2}{x_2}},\\
\dis{f_4(\tb{x}) =\frac{1-x_1-x_2-x_3}{x_1+x_2+x_3}\cdot\frac{1-x_4}{x_4}},\  \dis{f_3(\tb{x})=\frac{1-x_1-x_2-x_3-x_7-x_8}{x_1+x_2+x_3+x_7+x_8}\cdot\frac{1-x_3}{x_3}},\\
\dis{f_7(\tb{x})=\frac{1-x_4-x_5-x_6}{x_4+x_5+x_6}\cdot\frac{1-x_7}{x_7}},\ \dis{f_5(\tb{x})=\frac{1-x_1-x_2-x_3-x_7-x_8}{x_1+x_2+x_3+x_7+x_8}\cdot\frac{1-x_5}{x_5}},\\
\dis{f_8(\tb{x})=\frac{1-x_1-x_2-x_3}{x_1+x_2+x_3}\cdot\frac{1-x_8}{x_8}},\ \dis{f_6(\tb{x})=\frac{1-x_4-x_5-x_6-x_7-x_8}{x_4+x_5+x_6+x_7+x_8}\cdot\frac{1-x_6}{x_6}},
\end{array}
\end{displaymath}
\begin{displaymath}
\begin{array}{l}
\dis{g_1(\tb{x})=\frac{x_7}{1-x_7}\cdot\frac{1-x_1}{x_1}},\  \dis{g_2(\tb{x})=\frac{x_6}{1-x_6}\cdot\frac{1-x_2}{x_2}},\\
\dis{g_3(\tb{x}) =\frac{x_5}{1-x_5}\cdot\frac{1-x_3}{x_3}},\  \dis{g_4(\tb{x})=\frac{x_8}{1-x_8}\cdot\frac{1-x_4}{x_4}},\\
\dis{g_5(\tb{x})=\frac{x_3}{1-x_3}\cdot\frac{1-x_5}{x_5}},\ \dis{g_6(\tb{x})=\frac{x_2}{1-x_2}\cdot\frac{1-x_6}{x_6}},
\end{array}
\end{displaymath}
\begin{displaymath}
\begin{array}{l}
\dis{h_1(\tb{x})=\frac{1-x_6}{x_6}\cdot\frac{1-x_8}{x_8}},\  \dis{h_2(\tb{x})=\frac{1-x_2}{x_2}\cdot\frac{1-x_8}{x_8}},\\
\dis{h_3(\tb{x}) =\frac{1-x_5}{x_5}\cdot\frac{1-x_7}{x_7}},\  \dis{h_4(\tb{x})=\frac{1-x_3}{x_3}\cdot\frac{1-x_7}{x_7}},
\end{array}
\end{displaymath}
for $\tb{x}\in\Delta^7=\left\{(x_1,x_2,\dots,x_8)\in\mathbb{R}^8:\sum_{i=1}^8x_i=1,\ x_i>0,\ 1\leq i\leq 8\right\}.$

We will not be concerned with the functions $\{h_k\}_{k=1}^4$ in this paper. Because they provide information about the displacements under $\eta^{-2}$, $\eta^{-1}\xi$, $\xi\eta^{-1}$, and $\xi^{2}$. Only the functions $\{f_i\}_{i=1}^8$ and $\{g_j\}_{j=1}^6$ are related to the displacements under $\xi$, $\eta$ and $\xi\eta$.  Furthermore, it is possible to show that just the first eight $\{f_i\}_{i=1}^8$ are significant to find a lower bound for the maximum of the displacements $\dxi$, $\deta$ and $\dxieta$.

We consider $\dseven$ as a submanifold of $\mathbb{R}^8$. Let $I^{\dagger}=\{1,2,3,4,5,6,7,8\}$. Each function $f_i$ for $i\in I^{\dagger}$ is smooth in an open neighborhood of $\dseven$. Then the calculation of the infimum of the maximum of the functions $\{f_j\}_{j\in I^{\dagger}}$ follows from the statements
\begin{itemize}
\item[c.]
$\inf_{\tb{x}\in\Delta^7}\left\{\max\left(f_1(\tb{x}),\dots,f_8(\tb{x})\right)\right\}=\min_{\tb{x}\in\Delta^7}\left\{\max\left(f_1(\tb{x}),\dots,f_8(\tb{x})\right)\right\}$,

\item[d.] $\min_{\tb{x}\in\Delta^7}\left\{\max\left(f_1(\tb{x}),f_2(\tb{x}),\dots,f_8(\tb{x})\right)\right\}=f_1(\tb{x}^*)$ for
$\tb{x}^*\in\Delta_7\subset\Delta^7$,
\end{itemize}
proved in Proposition \ref{lem4.2} in \S\ref{sec416}, where $\Delta_7=\left\{\tb{x}\in\Delta^7:f_j(\tb{x})=f_k(\tb{x}),j,k\in I^{\dagger}\right\}$. Although the proof of the
observation that $\tb{x}^*\in\Delta_7$ also uses the fact that none of the functions $\{f_j\}_{j\in I^{\dagger}}$  has a local minimum in $\Delta^7$, it is more involved. It requires a case by case study of the values of the functions $\{f_j\}_{j\in I^{\dagger}}$ at the point $\tb{x}^*$.

In \S\ref{sec411}, we use the first order partial derivatives of the functions $\{f_j\}_{j\in I^{\dagger}}$ to show that there are certain direction vectors in the tangent space $T_{\xs}\dseven$ so that moving along these vectors reduces the number of possible
cases for the values of the functions $\{f_j\}_{j\in I^{\dagger}}$ at the point $\tb{x}^*$ to only five main cases:
\begin{itemize}
\item[I.] $f_j(\tb{x}^*)<\alphas$ for $j\in\{4,5,6,7,8\}$ and $f_j(\tb{x}^*)=\alphas$ for $j\in\{1,2,3\}$,

\item[II.] $f_j(\tb{x}^*)<\alphas$ for $j\in\{1,2,3,4,5,6\}$ and $f_j(\tb{x}^*)=\alphas$ for $j\in\{7,8\}$,

\item[III.] $f_j(\tb{x}^*)=\alphas$ for $j\in\{1,2,3,4,5,6\}$ and $f_j(\tb{x}^*)<\alphas$ for $j\in\{7,8\}$,

\item[IV.] $f_j(\tb{x}^*)=\alphas$ for $j\in\{1,2,3,7,8\}$ and $f_j(\tb{x}^*)<\alphas$ for $j\in\{4,5,6\}$,

\item[V.] $f_j(\xs)=\alphas$ for $j\in I^{\dagger}$,
\end{itemize}
 for $\alphas=\min_{\tb{x}\in\Delta^7}\left\{\max\left(f_1(\tb{x}),f_2(\tb{x}),\dots,f_8(\tb{x})\right)\right\}$. In each of the cases I, II, III and IV, we prove in
\S\ref{sec416} that there exists a piece of curve in $\Delta^7$ passing through the point $\tb{x}^*$ such that moving along this curve produces a point in
$\Delta^7$ at which a smaller minimum for the maximum of the functions $\{f_j\}_{j\in I^{\dagger}}$ is attained. This leaves only the case
$\tb{x}^*\in\Delta_7$ which suggests a method to find the coordinates of the point $\tb{x}^*$. Then we evaluate one of the displacement functions in $\{f_j\}_{j\in I^{\dagger}}$ at $\xs$ to calculate $\alphas$. In other words, we prove the following statement
\begin{thm}
Let $F^{\dagger}\co\dseven\to\mathbb{R}$ be the function defined by $\tb{x}\mapsto\max\{f_i(\tb{x}): i\in I^{\dagger}\}$. Then we have $\inf_{\tb{x}\in\Delta^7}F^{\dagger}(\tb{x})=5+3\sqrt{2}$,
\end{thm}
given as Theorem \ref{thm4.14} in \S\ref{sec416}, which provides the main estimate of Theorem \ref{thm5.1} in the geometrically infinite case.

In Section \ref{sec5}, we show that the proper and continuous function defined for a fixed point $z_0\in\hyp$ on the character variety $\mathfrak{X}$ with the formula
\[f_{z_0}^{\dagger}(\xi,\eta)=\max\{\dxi,\deta,\dxieta\}\]
 has no local minimum in $\mathfrak{GF}$ to reduce the geometrically finite case to the geometrically infinite case, completing the proof of Theorem \ref{thm5.1}. Note that an analogous process applies to a broader class of decompositions.

As summarized above, when $\Gamma=\langle\xi,\eta\rangle$ is geometrically infinite, the infimum of the maximum of the displacement functions $f_1,\dots,f_8,g_1,\dots,g_6$, determined by the decomposition $\Gamma_{\mathcal{D}^{\dagger}}$ of $\Gamma$, over $\dseven$ provides a lower bound for the displacements under the isometries $\xi$, $\eta$ and $\xi\eta$ in Theorem \ref{thm5.1}. Notice that, similar to the displacement functions $f_1^1$, $f_2^1$, $f_3^1$ and $f_4^1$ given by the decomposition $\Gamma_{\mathcal{D}^1}$, the infimum of the maximum of the displacement functions $g_1$, $f_2$, $f_3$, $g_4$, $f_5$, $f_6$, $f_7$ and $f_8$ over $\dseven$ provides a lower bound for the displacements under the isometries $\xi$ and $\eta$. Although we have
\[
\inf_{\tb{x}\in\Delta^3}\left\{\max\left(f_1^1(\tb{x}),f_2^1(\tb{x}),f_3^1(\tb{x}),f_4^1(\tb{x})\right)\right\}=\inf_{\tb{x}\in\Delta^7}\left\{\max\left(g_1(\tb{x}),f_2(\tb{x})\dots,f_8(\tb{x})\right)\right\},
\]
for $\Gamma_{\mathcal{D}^1}$ and $\Gamma_{\mathcal{D}^{\dagger}}$ by Lemma \ref{lem2.1} (see \S\ref{S2}) and the fact that $g_1(\tb{x})=11/4$, $f_2(\tb{x})=9$, $f_3(\tb{x})=9$, $g_4(\tb{x})=11/4$, $f_5(\tb{x})=9$, $f_6(\tb{x})=9$, $f_7(\tb{x})=9$, and $f_8(\tb{x})=9$ for the point $\tb{x}=(1/12,1/12,1/12,1/12,1/12,1/12,1/4,1/4)\in\dseven$, it may be possible to obtain a larger lower bound than $\log 3$ introduced in the $\log 3$ Theorem by examining a carefully chosen decomposition or a sequence of decompositions of $\Gamma$. If a larger lower bound for the displacements under the isometries $\xi$ and $\eta$ can be achieved, all the results involving the $\log 3$ Theorem in \cite{ACS}, \cite{ACCS}, \cite{CSBetti}, and \cite{CSMargulis} can be improved accordingly.



\section{The $\log 3$ theorem revisited}\label{S2}

In this section, we calculate the number $\log 3$ using the view point proposed in \S\ref{sec1.2} in the case (\ref{I}) of the proof of the $\log 3$ theorem. 

Let $\mathcal{F}^1$ be the set of functions $\{f_1^1,f_2^1,f_3^1,f_4^1\}$ introduced in \S\ref{sec1.2}. Let us define the continuous function $G^1 \co  \Delta^3 \to \mathbb{R}$ so that $G^1(\tb{x})=\max\{f(\tb{x})\co f\in\mathcal{F}^1\}$.

We aim to find $\inf_{\tb{x}\in\Delta^3}G^1(\tb{x})$. To this purpose we shall consider $\Delta^3$ as a submanifold of $\mathbb{R}^{4}$. The tangent space $T_{\tb{x}}\Delta^3$ of $\Delta^3$ consists of the vectors in $\mathbb{R}^4$ whose entries sum to $0$ at any $\tb{x}\in\Delta^3$. Note that each displacement function $f_i^1$ for $i=1,2,3,4$ is smooth in an open neighborhood of $\Delta^3$. Therefore the directional derivative of $f_i^1$ in the direction of any $\vec v\in T_{\tb{x}}\Delta^3$ is given by $\nabla f_i^1(\tb{x})\cdot\vec v$ for any $i=1,2,3,4$.

Let us introduce the function $\sigma\co(0,1)\to (0,\infty)$ defined by $\sigma(x)=1/x-1$, where $\sigma'(x)=-1/x^2<0$ for every $x\in (0,1)$. Observe that we have the equality $\inf_{\tb{x}\in\Delta^3}G^1(\tb{x})=\inf_{\tb{x}\in\Delta^3}\{\max(f_1^1(\tb{x}),f_2^1(\tb{x}))\}$ because  $f_1^1(\tb{x})=f_4^1(\tb{x})$
and $f_2^1(\tb{x})=f_3^1(\tb{x})$ for every $\tb{x}\in\Delta^3$. In other words, it is enough to prove the following:
\begin{lemma}\label{lem2.1} If $f_1^1(\tb{x})=\sigma(x_4)\sigma(x_1)$ and $f_2^1(\tb{x})=\sigma(x_3)\sigma(x_2)$ are the displacement functions defined in \S\ref{sec1.2} for $\tb{x}\in\Delta^3$, then
$\inf_{\tb{x}\in\Delta^3}\left\{\max\left(f_1^1(\tb{x}),f_2^1(\tb{x})\right)\right\}=9.$
\end{lemma}
\begin{proof}  Let $F^1 \co\Delta^3\to \mathbb{R}$ be the function defined by $\tb{x} \mapsto\max\left(f_1^1(\tb{x}),f_2^1(\tb{x})\right)$. Since $f_1^1$ and $f_2^1$ are continuous on $\Delta^3$, the function $F^1$ is also continuous. The number $\inf_{\tb{x}\in\Delta^3}F^1(\tb{x})$ exists and is greater than or equal to $1$ because the inequalities $f_1^1(\tb{x})>1$ and $f_2^1(\tb{x})>1$ hold for every $\tb{x}\in\Delta^3$. We aim to show that $\inf_{\tb{x}\in\Delta^3}F^1(\tb{x})=9$. First we prove that
$\inf\nolimits_{\tb{x}\in\Delta^3}F^1(\tb{x})=\min\nolimits_{\tb{x}\in\Delta^3}F^1(\tb{x}).$

Let $\Delta_n=\left\{\left(x_1,x_2,x_3,x_4\right)\in\Delta^3\co 1/n\leq x_i\leq 1-1/(2n)\  \tr{for}\ i=1,2,3,4\right\}$ for $n\geq 2$,
then $F^1$ has an absolute minimum $\alpha_n$ at a point $\tb{x}_n$ in  $\Delta_n$ as $\Delta_n$ is compact. The sequence $(\alpha_n)_{n=1}^{\infty}$ is  decreasing because we have $\Delta_{n}\subset\Delta_{n+1}$. In other words, the sequence $(F^1(\tb{x}_n))_{n=1}^{\infty}$ tends to an infimum of $F^1$. Assume that $\tb{x}_n$ approaches to a point $\tb{x}_0=(b_1,b_2,b_3,b_4)\in\partial\Delta^3$ as $n$ approaches to infinity. Then we get $b_i=0$ for some $i=1,2,3,4$. Suppose that $b_1=0$. By the definition of $f_1^1$, we must have
$b_4=1$. Otherwise $(F^1(\tb{x}_n))_{n=1}^{\infty}$ would approach to infinity instead of tending to an infimum. Then we conclude that $b_3=0$ and $b_2=0$. But, the function $(F^1(\tb{x}_n))_{n=1}^{\infty}$ tends to infinity by the definition of $f_2^1$ in this case, a contradiction. Thus, we get $b_1\neq 0$. Similar contradictions arise under the assumption $b_i=0$ for any $i\in\{2,3,4\}$ implying that $b_i\neq 0$ for every $i\in\{2,3,4\}$. Therefore $(\tb{x}_n)_{n=1}^{\infty}$ cannot have a limit point on the boundary of $\Delta^3$.

We claim that there exists a positive integer $n_0$ so that $\tb{x}_n=\tb{x}_{n+1}$ for every $n\geq n_0$. Let us assume otherwise that we
have a subsequence $(\tb{x}_{n_j})_{j=1}^{\infty}$ of the sequence $(\tb{x}_n)_{n=1}^{\infty}$ so that $\tb{x}_{n_j}\in\Delta_{n_{j+1}}-\Delta_{n_j}$ for every integer $j>0$. Since we
have $\cup_{n=2}^{\infty}\Delta_n=\Delta^3$, there exists a subsequence of  $(\tb{x}_{n_j})_{j=1}^{\infty}$ which has a limit point on $\partial\Delta^3$, a contradiction. In other words the absolute minimum of $F^1$ is attained at a point $\tb{x}^*=\left(x_1^*,x_2^*,x_3^*,x_4^*\right)\in\Delta^3$ so that
$F^1(\xs)=\inf_{\tb{x}\in\Delta^3}F^1(\tb{x})$.

Let $\Delta_3=\left\{\tb{x}\in\Delta^3 \co f_1^1(\tb{x})=f_2^1(\tb{x})\right\}$. We claim that $\tb{x}^*\in\Delta_3$. Assume otherwise that $f_1^1(\tb{x}^*)> f_2^1(\tb{x}^*)$. Then there exists a neighborhood $U$ of $\xs$ such that $f_1^1(\tb{x})> f_2^1(\tb{x})$ for every $\tb{x}\in U$. In particular we get $F^1(\tb{x})=f_1^1(\tb{x})$ on $U$.  Since $F^1(\xs)$ is the minimum value, the function $f_1^1$ must have a critical point at $\xs$. This is a contradiction. Because the derivative of $f_1^1$ in the direction of $\vec v=\langle 1,-1,0,0\rangle\in T_{\tb{x}}\Delta^3$ calculated as
\begin{equation*}
\nabla f_1^1(\tb{x})\cdot\vec v=-\frac{1-x_4}{x_1^2x_4}<0
\end{equation*}
implies that $f_1^1$ decreases at any $\tb{x}\in\Delta^3$ in the direction of $\vec v$. This argument also applies, mutatis mutandis, to show that the assumption $f_2^1(\tb{x}^*)>f_1^1(\tb{x}^*)$ leads to a contradiction as the directional derivative $\nabla f_2^1(\tb{x})\cdot\vec v>0$ at any $\tb{x}\in\Delta^3$. Hence we obtain that $\xs\in\Delta_3$.

Let $f_1(x)=\sigma(x_4)\sigma(x_1)$ and $f_2(x)=\sigma(x_3)\sigma(x_2)$ be the extensions of $f_1^1$ and $f_2^1$ to the open set $C=(0,1)\times (0,1)\times (0,1)\times (0,1)$. Let us consider the Lagrange multipliers problem of finding the minimum of $f_1$ subject to the constraints $G_1(x)=0$ and $G_2(x)=0$ for $x=(x_1,x_2,x_3,x_4)\in C$ where  $G_1$ and  $G_2$ are the functions defined by
$$G_1(x)=f_1(x)-f_2(x)\ \ \tnr{and}\ \ G_2(x)=x_1+x_2+x_3+x_4-1.$$
Since we have $f_1(\xs)=\min_{\tb{x}\in\Delta^3}F^1(\tb{x})$ and $\xs\in\Delta_3$, the point $\xs$ is a solution to this problem. Therefore the gradient vector $\nabla f_1(\xs)$ is in the subspace of $\mathbb{R}^4$ spanned by the vectors $\nabla G_1(\xs)$ and $\nabla G_2(\xs)$. Equivalently the matrix
\begin{equation*}
\left[\begin{array}{c}
        \nabla f_1(\xs)\\
        \nabla G_1(\xs)  \\
        \nabla G_2(\xs)
       \end{array}
\right]=\left[\begin{array}{cccc}
                 \sigma'(x_1^*)\sigma(x_4^*) & 0 & 0 & \sigma(x_1^*)\sigma'(x_4^*)\\
                 \sigma'(x_1^*)\sigma(x_4^*) & -\sigma'(x_2^*)\sigma(x_3^*) & -\sigma(x_2^*)\sigma'(x_3^*) & \sigma(x_1^*)\sigma'(x_4^*) \\
                 1 & 1 & 1 & 1
                \end{array}
         \right]\label{matrix1}
\end{equation*}
has rank less than $3$. By applying the row reduction operations $R_2\mapsto (-1)R_1+R_2$, and $R_1\mapsto (-\sigma(x_1^*)\sigma'(x_4^*))R_3+R_1$, and $R_1\mapsto (\sigma(x_1^*)\sigma'(x_4^*))/(\sigma(x_2^*)
\sigma'(x_3^*))R_2+R_1$ simultaneously, we see that the matrix above is row equivalent to
\begin{equation*}
\left[\begin{array}{cccc}
                 \dis{\frac{(x_1^*-x_4^*)(1-x_1^*-x_4^*)}{(x_1^*x_4^*)^2}} & \dis{\frac{(x_2^*-x_3^*)(1-x_2^*-x_3^*)\sigma(x_1^*)\sigma'(x_4^*)}{(x_2^*x_3^*)^2\sigma(x_2^*)\sigma'(x_3^*)}} & 0 & 0\\
                 0 & -\sigma'(x_2^*)\sigma(x_3^*) & -\sigma(x_2^*)\sigma'(x_3^*) & 0 \\
                 1 & 1 & 1 & 1
                \end{array}
         \right].\label{matrix1}
\end{equation*}
The reduced matrix above has rank less than $3$ if and only if $x_1^*=x_4^*$ and $x_2^*=x_3^*$. Then it follows from the use of the fact $f_1(\xs)=f_2(\xs)$ that $x_1^*=x_2^*$. As a result we find that $x_i^*=1/4$ for every $i\in\{1,2,3,4\}$. Finally a simple evaluation shows that $\inf_{\tb{x}\in\Delta^3}F^1(\tb{x})=9$.
\end{proof}
It is worth emphasizing a few key points used in the proof of Lemma \ref{lem2.1}. The functions $f_1^1$ and $f_2^1$ have no local minimum on $\Delta^3$. This fact implies that $\tb{x}^*\in\Delta_3$. The proof of Lemma \ref{lem2.1} shows that the main ingredients used to calculate the
number $\log 3$ are
\begin{enumerate}[label=\arabic*]
\renewcommand{\labelenumi}{\alph{enumi}}
\item\hspace{-.3cm}.\label{a} $\inf_{\tb{x}\in\Delta^3}F^1(\tb{x})=\min_{\tb{x}\in\Delta^3}F^1(\tb{x})$,
\item\hspace{-.3cm}.\label{b} there exists a  point $\tb{x}^*\in\Delta_3\subset\Delta^3$ such that $\min_{\tb{x}\in\Delta^3}F^1(\tb{x})=f_1^1(\tb{x}^*)$.
\end{enumerate}
Lemma \ref{lem2.1} also establishes that the point $\tb{x}^*\in\Delta^3$ is unique. If we assume the uniqueness of the point $\tb{x}^*$ a priori
together with (\ref{a}), it is possible to suggest an alternative way of finding the coordinates of the point $\tb{x}^*$: Let $T_1, T_2\co\mathbb{R}^4\to \mathbb{R}^4$ be the maps defined by $(x_1,x_2,x_3,x_4)\mapsto
(x_2,x_1,x_4,x_3)$ and $(x_1,x_2,x_3,x_4)\mapsto (x_4,x_3,x_2,x_1)$, respectively. We have  $T_i(\tb{x})\in\Delta^3$
and  $f_i^1(T_1(\tb{x}))=f_i^1(\tb{x})$ for every $\tb{x}\in\Delta^3$. Let  $H_1\co\Delta^3\to \mathbb{R}$ be the function so that $H_1(\tb{x})=\max\{(f_i^1\circ T_1)(\tb{x})\co i=1,2\}$. We see that
$F^1(\tb{x})=H_1(\tb{x})$ for every $\tb{x}\in\Delta^3$. Since $F^1$ takes
its minimum value at the point $\tb{x}^*$, the function $H_1$ takes its minimum value at the point
$T_1^{-1}(\tb{x}^*)$. Therefore we obtain $T^{-1}_1(\tb{x}^*)=\tb{x}^*$ which implies that $x_1^*=x_2^*$,
and $x_3^*=x_4^*$.

Let $H_2\co\Delta^3\to \mathbb{R}$ be the function defined by $H_2(\tb{x})=\max\{(f_i^1\circ T_2)(\tb{x})\co i=1,2\}$. We repeat an analog of the argument above for $H_2$ to derive that $x_1^*=x_4^*$, and $x_2^*=x_3^*$. As a result we find $x_i^*=x_j^*$ for every $i,j=1,2,3,4$. In
this calculation method, we don't refer to the statement $\tb{x}^*\in\Delta_3$ given in (\ref{b}).



\section{Decomposition of $\Gamma=\lag\xi,\eta\rag$ for the isometries $\xi$, $\eta$ and $\xi\eta$}\label{Sec3}

Let $\Gamma$ be a group which is free on a finite generating set $\Xi$. Let $\Xi^{-1}=\{\psi^{-1}:\psi\in\Xi\}$. For $m\geq 2$, every element $\gamma$ of $\Gamma$ can be written uniquely as a reduced word
$\psi_{1}\cdots\psi_{m}$, where each $\psi_i$ is an element of $\Xi\cup\Xi^{-1}$ and $\psi_{i+1}\neq\psi_{i}^{-1}$ for $i=1,\dots ,m-1$. If
$k\leq m$ is a positive integer and $\gamma\neq 1$, we shall call $\psi_1\dots\psi_k$ the \textit{initial word of length $k$} of $\gamma$.

Let $\Psi^*$ be a finite set of words in $\Gamma$.  For each word $\psi\in\Psi^*$, let $J_{\psi}$ denote the set of non-trivial elements of $\Gamma$ that have
initial word $\psi$. Depending on the number of elements in $\Xi$ and lengths of words in $\Psi^*$ there may be a set of words which are not contained in
any of $J_{\psi}$. Let us call this set the \textit{residue set} of $\Psi^*$ and denote it by $\Psi_r^*$.
\begin{definition}\label{Def3.1}
For a given pair $\mathcal{D}^*=(\Psi^*, \Psi_r^*)$ of finite, disjoint sets of words $\Psi^*$ and $\Psi_r^*$ in $\Gamma$, if $\Gamma=\{1\}\cup
\Psi_r^*\cup\bigcup_{\psi\in\Psi^*}J_{\psi}$, then $\Gamma_{\mathcal{D}^*}$  will be called a \textit{decomposition} of $\Gamma$.
\end{definition}

\begin{definition}
A decomposition $\Gamma_{\mathcal{D}^*}$ with $\mathcal{D}^*=(\Psi^*,\Psi^*_r)$ is \textit{symmetric} if $\Psi^*$ and $\Psi^*_r$ are preserved by every bijection of $\Xi\cup\Xi^{-1}$.
\end{definition}

We know that $\Gamma=\langle\xi,\eta\rangle$ described in the $\log 3$ theorem is a free group on $\Xi=\{\xi,\eta\}$ by Proposition $9.2$ in \cite{CSParadox}. For $\mathcal{D}^1=(\Psi^1,\Psi^1_r=\emptyset)$, the decomposition
$\Gamma_{\mathcal{D}^1}$ corresponds to the symmetric decomposition of $\Gamma$ in (\ref{prdx.d.}).

We introduce another decomposition of $\Gamma$ that contains the set $\Gamma_{\dagger}=\{\xi,\eta,\xi\eta\}$. Let $\Psi^{\dagger}=\{\xi\eta, \xi^2, \xi\eta^{-1}, \eta, \xi^{-1}, \eta^{-1}\xi^{-1}, \eta^{-2}, \eta^{-1}\xi\}$. Then it is straightforward to see that 
\begin{eqnarray}\label{decomp}
\Gamma =  \{1\}\cup \Psi_r^{\dagger}\cup\bigcup\nolimits_{\psi\in\Psi^{\dagger}}J_{\psi},
\end{eqnarray}
where the residue set is $\Psi^{\dagger}_r=\{\xi,\eta^{-1}\}$. Note that $\Gamma_{\mathcal{D}^{\dagger}}$ for $\mathcal{D}^{\dagger}=(\Psi^{\dagger},\Psi^{\dagger}_r)$ is not symmetric. We shall use $\Gamma_{\mathcal{D}^{\dagger}}$ in the proof of Theorem \ref{thm5.1}. In particular, we aim to prove Theorem \ref{thm3.4}, an analog of Theorem \ref{thm1.1}, for $\Psi^{\dagger}$.

We first state a more general result Lemma \ref{lem3.2}, an analog of Lemma 5.3 in \cite{CSParadox}, for $\Psi^{\dagger}$. Although Lemma \ref{lem3.2}
follows directly from the same arguments used in the  proof of Lemma 5.3 \cite{CSParadox}, its proof is included for the sake of completeness. The main tool is Proposition 4.2 in \cite{CSParadox}. 
In the following, $J_{\Psi}$ denotes the disjoint union $\bigcup_{\psi\in\Psi}J_{\psi}$ for $\Psi\subset\Psi^{\dagger}$:

\begin{lemma} \label{lem3.2}
Let $\Gamma$ be a Kleinian group which is free on a generating set $\{\xi,\eta\}$.  Let $z_0$ be any point of $\hyp$. Then there exists a number $D\in[0,2]$, a $\Gamma$--invariant $D$--conformal density $\mathcal{M}=(\mu_z)$
for $\hyp$ and a family $\{\nu_{\psi}\}_{\psi\in\Psi^{\dagger}}$ of Borel measures  on $S_{\infty}$ such that
\begin{itemize}
\item[(1)]  $\mu_{z_0}(S_{\infty})=1$,\quad (2) \   $\mu_{z_0}=\dis{\sum_{\psi\in\Psi^{\dagger}}\nu_{\psi}}$, 
\item[(3)] \begin{itemize}
      \item[(a)] $\dis{\int_{S_{\infty}}(\lambda_{\xi,z_0})^Dd\nu_{\xi^{-1}}=1-\int_{S_{\infty}} d\nu_{\xi\eta}-\int_{S_{\infty}} d\nu_{\xi^2}-\int_{S_{\infty}} d\nu_{\xi\eta^{-1}}}$,\\
      \item[(b)] $\dis{\int_{S_{\infty}}(\lambda_{\xi^{-1},z_0})^Dd\nu_{\xi\eta^{-1}}  =  \int_{S_{\infty}} d\nu_{\eta^{-1}\xi^{-1}}+\int_{S_{\infty}} d\nu_{\eta^{-1}\xi}+\int_{S_{\infty}} d\nu_{\eta^{-2}}}$,\\
      \item[(c)] $\dis{\int_{S_{\infty}}(\lambda_{\xi^{-1},z_0})^Dd\nu_{\xi^{2}}  =  \int_{S_{\infty}} d\nu_{\xi\eta}+\int_{S_{\infty}} d\nu_{\xi^2}+\int_{S_{\infty}} d\nu_{\xi\eta^{-1}}}$,\\
      \item[(d)] $\dis{\int_{S_{\infty}}(\lambda_{\xi^{-1},z_0})^Dd\nu_{\xi\eta}  =  \int_{S_{\infty}} d\nu_{\eta}}$,
      \end{itemize}
\item[(4)] \begin{itemize}
      \item[(a)] $\dis{\int_{S_{\infty}}(\lambda_{\eta^{-1},z_0})^Dd\nu_{\eta}  =  1-\int_{S_{\infty}} d\nu_{\eta^{-1}\xi}-\int_{S_{\infty}} d\nu_{\eta^{-2}}-\int_{S_{\infty}} d\nu_{\eta^{-1}\xi^{-1}}}$,\\
      \item[(b)] $\dis{\int_{S_{\infty}}(\lambda_{\eta,z_0})^Dd\nu_{\eta^{-2}}  =  \int_{S_{\infty}} d\nu_{\eta^{-1}\xi^{-1}}+\int_{S_{\infty}} d\nu_{\eta^{-1}\xi}+\int_{S_{\infty}} d\nu_{\eta^{-2}}}$,\\
      \item[(c)] $\dis{\int_{S_{\infty}}(\lambda_{\eta,z_0})^Dd\nu_{\eta^{-1}\xi}  =  \int_{S_{\infty}} d\nu_{\xi\eta}+\int_{S_{\infty}} d\nu_{\xi^2}+\int_{S_{\infty}} d\nu_{\xi\eta^{-1}}}$,\\
      \item[(d)] $\dis{\int_{S_{\infty}}(\lambda_{\eta,z_0})^Dd\nu_{\eta^{-1}\xi^{-1}}  =  \int_{S_{\infty}} d\nu_{\xi^{-1}}}$,
      \end{itemize}
\item[(5)] \begin{itemize}
      \item[(a)] $\dis{\int_{S_{\infty}}(\lambda_{\eta^{-1}\xi^{-1},z_0})^Dd\nu_{\xi\eta}  =  1-\int_{S_{\infty}} d\nu_{\eta^{-1}\xi}-\int_{S_{\infty}} d\nu_{\eta^{-2}}-\int_{S_{\infty}} d\nu_{\eta^{-1}\xi^{-1}}}$,\\
      \item[(b)] $\dis{\int_{S_{\infty}}(\lambda_{\eta^{-1}\xi^{-1},z_0})^Dd\nu_{\xi\eta^{-1}}  =  \int_{S_{\infty}} d\nu_{\eta^{-2}}}$,\\
      \item[(c)] $\dis{\int_{S_{\infty}}(\lambda_{\eta^{-1}\xi^{-1},z_0})^Dd\nu_{\xi^{2}}  =  \int_{S_{\infty}} d\nu_{\eta^{-1}\xi}}$,\\
      \item[(d)] $\dis{\int_{S_{\infty}}(\lambda_{\xi\eta,z_0})^Dd\nu_{\eta^{-1}\xi^{-1}} = 1-\int_{S_{\infty}} d\nu_{\xi\eta}-\int_{S_{\infty}} d\nu_{\xi^2}-\int_{S_{\infty}} d\nu_{\xi\eta^{-1}}}$,\\
      \item[(e)] $\dis{\int_{S_{\infty}}(\lambda_{\xi\eta,z_0})^Dd\nu_{\eta^{-1}\xi}  =  \int_{S_{\infty}} d\nu_{\xi^2}}$, \ (f) $\dis{\int_{S_{\infty}}(\lambda_{\xi\eta,z_0})^Dd\nu_{\eta^{-2}}  =  \int_{S_{\infty}} d\nu_{\xi\eta^{-1}}}$, 
      \end{itemize}
\item[(6)] \begin{itemize}
      \item[(a)] $\dis{\int_{S_{\infty}}(\lambda_{\eta^{-1}\xi,z_0})^Dd\nu_{\xi^{-1}} =  1-\int_{S_{\infty}} d\nu_{\eta^{-1}\xi}}$,\\
      \item[(b)] $\dis{\int_{S_{\infty}}(\lambda_{\xi\eta^{-1},z_0})^Dd\nu_{\eta}  =  1-\int_{S_{\infty}} d\nu_{\xi\eta^{-1}}}$,\\
      \item[(c)] $\dis{\int_{S_{\infty}}(\lambda_{\xi^{2},z_0})^Dd\nu_{\xi^{-1}}  =  1-\int_{S_{\infty}} d\nu_{\xi^{2}}}$, \ (d) $\dis{\int_{S_{\infty}}(\lambda_{\eta^{-2},z_0})^Dd\nu_{\eta}  =  \int_{S_{\infty}} d\nu_{\eta^{-2}}}$.
      \end{itemize}
\end{itemize}
\end{lemma}
\begin{proof} 
Since $\Gamma$ acts
freely on $\hyp$ and it can be decomposed as in (\ref{decomp}), the orbit $W=\Gamma\cdot z_0$ is a disjoint union
\begin{equation}\label{decomp2LS_2}
W^{\dagger}=\{z_0\}\cup V_0\cup\bigcup\nolimits_{\psi\in\Psi^{\dagger}}V_{\psi},
\end{equation}
\noi where $V_0=\{\gamma\cdot z_0\co\gamma\in\Psi^{\dagger}_r\}$ and $V_{\psi}=\{\gamma\cdot z_0\co\gamma\in J_{\psi}\}$. Note that $V_0$ is the finite set $V_{\xi}\cup V_{\eta^{-1}}=\{\xi\cdot z_0\}\cup\{\eta^{-1}\cdot z_0\}$. Let $\mathcal{V}$ denote the finite collection 
of all sets of the form $\bigcup_{\psi\in\Psi}V_{\psi}$ or $V_0\cup\bigcup_{\psi\in\Psi}V_{\psi}$ or $\{z_0\}\cup\bigcup_{\psi\in\Psi}V_{\psi}$ or $\{z_0\}\cup
V_0\cup\bigcup_{\psi\in\Psi}V_{\psi}$ for $\Psi\subset\Psi^{\dagger}$. We apply Proposition 4.2 \cite{CSParadox} to $W$ and $\mathcal{V}$.

Let $D$ be a number in $[0,2]$, and $(\mathcal{M}_V)_{V\in\mathcal{V}}$ be a family of conformal densities, for which conditions ($i$)--($iv$) of Proposition
4.2 \cite{CSParadox} are satisfied. Let $\mathcal{M}_V=(\mu_{V,z})_{z\in\hyp}$. We set $\mathcal{M}=\mathcal{M}_W$, and
$\nu_{\psi}=\mu_{V_{\psi},z_0}$ for each $\psi\in\Psi^{\dagger}$. By Proposition 4.2 ($iii$), $\mathcal{M}$ is $\Gamma$--invariant. By Proposition 4.2 ($i$) and the definition of a conformal density, we have $\mu_{z_0}(S_{\infty})=\mu_{W,z_0}(S_{\infty})\neq 0$. Therefore, we may assume that $\mu_{z_0}$ has total mass $1$ after normalization, which gives (1) of Lemma \ref{lem3.2}.

By Proposition 4.2 ($iv$), we have $\mu_{\{z_0\},z_0}=0$ and $\mu_{V_0,z_0}=0$. Applying Proposition 4.2 ($ii$) to the disjoint union in (\ref{decomp2LS_2}), we obtain
\begin{eqnarray*}
\mu_{z_0}  & = & \mu_{\{z_0\}, z_0}+\mu_{V_0,z_0}+\sum\nolimits_{\psi\in\Psi^{\dagger}}\mu_{V_{\psi},z_0}.
\end{eqnarray*}
Hence, we get conclusion (2) of Lemma \ref{lem3.2}.
In order to complete parts ($3$)--($6$) of the lemma, we need to determine all of the group theoretical relations
between the sets of words $J_{\psi}$ for $\psi\in\Psi^{\dagger}$:
%
We know that $\xi^{-1}\eta\in J_{\xi^{-1}}$. Therefore, we have $1\in\eta^{-1}\xi J_{\xi^{-1}}$. Let $w$ be a word in $J_{\xi^{-1}}$. Then we have
$w=\xi^{-1} w_1$ for some $w_1\in\Gamma$.  We compute that $\eta^{-1}\xi w=\eta^{-1} w_1$. The first letter of $w_1$ cannot be $\xi$. But it can
be either $\eta$, $\eta^{-1}$ or $\xi^{-1}$. Assume that it is $\eta$ and $w_1\neq\eta$. Then we have $w_1=\eta w_2$ for
some word $w_2\in\Gamma$. The first letter of $w_2$ cannot be $\eta^{-1}$, but it can be either $\eta$, $\xi^{-1}$ or $\xi$. Since we get $\eta^{-1}\xi
w=w_2$, we derive that
$
\{1\}\cup J_{\eta}\cup J_{\xi\eta^{-1}}\cup J_{\xi\eta}\cup J_{\xi^{2}}\cup J_{\xi^{-1}}\subset \eta^{-1}\xi
J_{\xi^{-1}}.$
If the first letter of $w_1$ is $\eta^{-1}$, then we get $w_1=\eta^{-1} w_2$ for some $w_2\in\Gamma$. We see
that $\eta^{-1}\xi w=\eta^{-2} w_2$. This means that $J_{\eta^{-2}}\subset\eta^{-1}\xi J_{\xi^{-1}}$. If the first letter of $w_1$ is
$\xi^{-1}$. Then we get $w_1=\xi^{-1} w_2$ for some $w_2\in\Gamma$ which implies that $\eta^{-1}\xi w=\eta^{-1}\xi^{-1} w_2$. Therefore, we find
that $J_{\eta^{-1}\xi^{-1}}\subset\eta^{-1}\xi J_{\xi^{-1}}$. In other words, $\eta^{-1}\xi J_{\eta}$ contains every word in $\Gamma$ except the ones
start with $\eta^{-1}\xi$. Hence, we conclude that $\eta^{-1}\xi J_{\xi^{-1}}  =  \Gamma-J_{\eta^{-1}\xi}$.

Similar computations show that $\eta^{-1}\xi J_{\xi^{-1}}  =  \Gamma-J_{\xi\eta^{-1}}$, $\xi^{2} J_{\xi^{-1}}  =  \Gamma-J_{\xi^{2}}$, and $\eta^{-2}J_{\eta}  =  \Gamma-J_{\eta^{-2}}$.
\begin{table}[h]
\begin{center}
\begin{tabular}{|r|c|c|c|}
  \hline
 &$\gamma$ & $s(\gamma)$ & $S(\gamma)$\\
  \hline
  (3) (a) &$\xi$ & $\xi^{-1}$ & $\{\xi\eta,\xi^2,\xi\eta^{-1}\}$ \\ 
  \hline
  (b)&$\xi^{-1}$ & $\xi\eta^{-1}$  &  $\{\xi\eta, \xi^2, \xi\eta^{-1}, \eta, \xi^{-1}\}$ \\ 
  \hline
  (c)&$\xi^{-1}$  & $\xi^2$  & $\{\eta, \xi^{-1}, \eta^{-1}\xi^{-1}, \eta^{-2}, \eta^{-1}\xi\}$ \\  
  \hline
  (d)&$\xi^{-1}$ & $\xi\eta$ & $\{\xi\eta, \xi^2, \xi\eta^{-1},  \xi^{-1}, \eta^{-1}\xi^{-1}, \eta^{-2}, \eta^{-1}\xi\}$ \\ 
  \hline
  (4) (a)&$\eta^{-1}$ & $\eta$ & $\{\eta^{-1}\xi^{-1}, \eta^{-2}, \eta^{-1}\xi\}$ \\ 
  \hline
  (b)&$\eta$ & $\eta^{-2}$ & $\{\xi\eta, \xi^2, \xi\eta^{-1}, \eta, \xi^{-1}\}$ \\
  \hline
  (c)&$\eta$ & $\eta^{-1}\xi$ & $\{\eta, \xi^{-1}, \eta^{-1}\xi^{-1}, \eta^{-2}, \eta^{-1}\xi\}$ \\ 
  \hline
  (d)&$\eta$ & $\eta^{-1}\xi^{-1}$ & $\{\xi\eta, \xi^2, \xi\eta^{-1}, \eta, \eta^{-1}\xi^{-1}, \eta^{-2}, \eta^{-1}\xi\}$ \\ 
  \hline
  (5) (a) &$\eta^{-1}\xi^{-1}$ & $\xi\eta$ & $\{\eta^{-1}\xi^{-1}, \eta^{-2}, \eta^{-1}\xi\}$ \\ 
  \hline
  (b) &$\eta^{-1}\xi^{-1}$ & $\xi^{2}$ & $\{\xi\eta, \xi^2, \xi\eta^{-1}, \eta, \xi^{-1}, \eta^{-1}\xi^{-1}, \eta^{-2}\}$ \\ 
  \hline
  (c) &$\eta^{-1}\xi^{-1}$ & $\xi\eta^{-1}$ & $\{\xi\eta, \xi^2, \xi\eta^{-1}, \eta, \xi^{-1}, \eta^{-1}\xi^{-1}, \eta^{-1}\xi\}$ \\ 
  \hline
  (d) &$\xi\eta$ & $\eta^{-1}\xi^{-1}$ & $\{\xi\eta, \xi^2, \xi\eta^{-1}\}$ \\ 
  \hline
  (e) &$\xi\eta$ & $\eta^{-1}\xi$ & $\{\xi\eta, \xi\eta^{-1}, \eta, \xi^{-1}, \eta^{-1}\xi^{-1}, \eta^{-2}, \eta^{-1}\xi\}$ \\ 
  \hline
  (f) &$\xi\eta$ & $\eta^{-2}$ & $\{\xi\eta, \xi^2, \eta, \xi^{-1}, \eta^{-1}\xi^{-1}, \eta^{-2}, \eta^{-1}\xi\}$ \\ 
  \hline
  (6) (a)&$\eta^{-1}\xi$ & $\xi^{-1}$ & $\{\eta^{-1}\xi\}$ \\
  \hline
  (b)&$\xi\eta^{-1}$ & $\eta$ & $\{\xi\eta^{-1}\}$ \\
  \hline
  (c) &$\xi^{2}$ & $\xi^{-1}$ & $\{\xi^2\}$ \\
  \hline
  (d) &$\eta^{-2}$ & $\eta$ & $\{\eta^{-2}\}$ \\
  \hline
\end{tabular}
\caption{Group--theoretical properties of the decomposition $\Gamma_{\mathcal{D}^{\dagger}}$.}
\label{table1}
\end{center}
\end{table}
It follows from the discussion above and definitions of the sets $J_{\psi}$ for each $\psi\in\Psi^{\dagger}$ that, for each row $\gamma$, $s(\gamma)$ and $S(\gamma)$ of Table \ref{table1}, the decomposition $\Gamma_{\mathcal{D}^{\dagger}}$ of $\Gamma=\langle\xi,\eta\rangle$ has the group--theoretical properties
\begin{equation}
\gamma J_{s(\gamma)}=\Gamma-J_{S(\gamma)}.\label{113}
\end{equation}

Let $V_{\Psi}$ denote the union $\bigcup_{\gamma\in\Psi}V_{\gamma}$ where $\Psi$ is a subset of $\Psi^{\dagger}\cup\Psi^{\dagger}_r$. Using the group--theoretical
relations in (\ref{113}), we derive the following relations
\begin{equation}
\gamma V_{s(\gamma)}  =  W - V_{S(\gamma)}\label{eqn14}
\end{equation}
between the orbits $V_{s(\gamma)}$ and $V_{S(\gamma)}$. 
Since we have $W -V_{S(\gamma)}=V_{s(\gamma)}\in\mathcal{V}$, condition $(iii)$ of Proposition 4.2 gives $\mathcal{M}_{V_{s(\gamma)}}=\gamma^*_{\infty}\left(\mathcal{M}_{W -V_{S(\gamma)}}\right).$
On the other hand, by Proposition 4.2 ($ii$), we get $\mathcal{M}=\mathcal{M}_{W -V_{S(\gamma)}}+\mathcal{M}_{V_{S(\gamma)}}$.
We combine the last two equalities to obtain
$ \mathcal{M}_{V_{s(\gamma)}}=\gamma^*_{\infty}\left(\mathcal{M}-\mathcal{M}_{V_{S(\gamma)}}\right), $
which implies that
\begin{equation}
d\mu_{V_{s(\gamma)},\gamma\cdot z_0}  =d\left(\gamma^*_{\infty} \left(\mu_{z_0}-\sum_{\psi\in S(\gamma)}\nu_{\psi}\right)\right). \label{total}
\end{equation}
Since $\mathcal{M}_{V_{s(\gamma)}}$ is a $D$--conformal density and
$d\mu_{V_{s(\gamma)},\gamma\cdot z_0}
                              =   \lambda^{D}_{\gamma,z_0}d\mu_{V_{s(\gamma)}}$
(\cite{CSParadox}, Proposition 2.4), we obtain the equality
\begin{equation*}
\int_{S_{\infty}}(\lambda_{\gamma,z_0})^{D}d\mu_{V_{s(\gamma)}}
                                                  =  1-\sum\nolimits_{\psi\in S(\gamma)}\int_{S_{\infty}} d\nu_{\psi}
\end{equation*}
for every row of Table \ref{table1} by equating the total masses of both sides of (\ref{total}), which provides parts ($3$)--($6$) of the lemma.
\end{proof}

The following is an analog of Theorem \ref{thm1.1} for the set $\Psi^{\dagger}\subset\Gamma=\langle\xi,\eta\rangle$. Notice that Theorem \ref{thm3.4} has no analog for part (3) of Theorem \ref{thm1.1}.

\begin{theorem}\label{thm3.4}
Let $\Gamma=\langle\xi,\eta\rangle$ be a free, geometrically infinite Kleinian group without parabolics. For any $z_0\in\hyp$, let $A_{z_0}$ be the area measure based at $z_0$. There is a family of Borel measures $\{\nu_{\psi}\}_{\psi\in\Psi^{\dagger}}$  for $\Psi^{\dagger}=\{\xi\eta, \xi^2, \xi\eta^{-1}, \eta, \xi^{-1}, \eta^{-1}\xi^{-1}, \eta^{-2}, \eta^{-1}\xi\}$ on $S_{\infty}$ such that
\begin{itemize}
 \item[(1)] $A_{z_0}=\sum_{\psi\in\Psi^{\dagger}}\nu_{\psi}$, where $A_{z_0}$ is normalized so that $A_{z_0}(S_{\infty})=1$, and,
\item[(2)]  $\dis{\int_{S_{\infty}}\left(\lambda_{\gamma,z_0}\right)^2d\nu_{s(\gamma)}=1-\sum_{\psi\in S(\gamma)}\int_{S_{\infty}} d\nu_{\psi}}$ for each row 
    of Table \ref{table1}.
\end{itemize}
\end{theorem}
\begin{proof}
By the conclusions of Propositions 6.9 and 3.9 in \cite{CSParadox} and tameness (\cite{Agol}, \cite{CG}) every $\Gamma$--invariant $D$--conformal density $\mathcal{M}$ is a constant multiple of the area density $\mathcal{A}$, i.e., $D=2$. By Lemma \ref{lem3.2} ($1$), we get $\mathcal{M}=\mathcal{A}$. Then ($2$) follows from Lemma \ref{lem3.2} ($3$)-($6$).
\end{proof}

We shall use Theorem \ref{thm3.4} together with Lemma \ref{lem1.2} to produce the displacement functions for the decomposition $\Gamma_{\mathcal{D}^{\dagger}}$. In the rest of this paper, we will use the bijection $p\co\Psi^{\dagger}\to I^{\dagger}$ defined by 
\begin{equation}\label{sigma}
\begin{array}{rrrr}
\xi\eta\mapsto 1,&  \xi^2\mapsto 2, &  \xi\eta^{-1}\mapsto 3, & \eta\mapsto 7,\\ \eta^{-1}\xi^{-1}\mapsto 4, & \eta^{-2}\mapsto 5, & \eta^{-1}\xi\mapsto 6, & \xi^{-1}\mapsto 8,
\end{array}
\end{equation}
to enumarate the displacement functions and their variables. We have:
\begin{proposition}\label{dispfunc}
Let $\Gamma=\langle\xi,\eta\rangle$ be a free Kleinian group. For any $z_0\in\hyp$ and for each $\gamma\in\{\xi,\eta,\xi^{-1},\eta^{-1},\xi\eta,\eta^{-1}\xi^{-1}\}$, the expression $e^{2\dgamma}$ is bounded below by $f_i(\tb{x})$ or $g_j(\tb{x})$ for at least one of $f_i$ or $g_j$ for $i\in I^{\dagger}=\{1,2,3,4,5,6,7,8\}$ and $j\in\{1,2,3,4,5,6\}$ on the list
\begin{displaymath}
\begin{array}{l}
\dis{f_1(\tb{x})=\frac{1-x_4-x_5-x_6}{x_4+x_5+x_6}\cdot\frac{1-x_1}{x_1}},\  \dis{f_2(\tb{x})=\frac{1-x_4-x_5-x_6-x_7-x_8}{x_4+x_5+x_6+x_7+x_8}\cdot\frac{1-x_2}{x_2}},\\
\dis{f_4(\tb{x}) =\frac{1-x_1-x_2-x_3}{x_1+x_2+x_3}\cdot\frac{1-x_4}{x_4}},\  \dis{f_3(\tb{x})=\frac{1-x_1-x_2-x_3-x_7-x_8}{x_1+x_2+x_3+x_7+x_8}\cdot\frac{1-x_3}{x_3}},\\
\dis{f_7(\tb{x})=\frac{1-x_4-x_5-x_6}{x_4+x_5+x_6}\cdot\frac{1-x_7}{x_7}},\ \dis{f_5(\tb{x})=\frac{1-x_1-x_2-x_3-x_7-x_8}{x_1+x_2+x_3+x_7+x_8}\cdot\frac{1-x_5}{x_5}},\\
\dis{f_8(\tb{x})=\frac{1-x_1-x_2-x_3}{x_1+x_2+x_3}\cdot\frac{1-x_8}{x_8}},\ \dis{f_6(\tb{x})=\frac{1-x_4-x_5-x_6-x_7-x_8}{x_4+x_5+x_6+x_7+x_8}\cdot\frac{1-x_6}{x_6}},
\end{array}
\end{displaymath}
(The functions above are produced from rows (3)(a)-(c), (4)(a)-(c), (5)(a) and (5)(d) of Table \ref{table1}),
\begin{displaymath}
\begin{array}{l}
\dis{g_1(\tb{x})=\frac{x_7}{1-x_7}\cdot\frac{1-x_1}{x_1}},\  \dis{g_2(\tb{x})=\frac{x_6}{1-x_6}\cdot\frac{1-x_2}{x_2}},\\
\dis{g_3(\tb{x}) =\frac{x_5}{1-x_5}\cdot\frac{1-x_3}{x_3}},\  \dis{g_4(\tb{x})=\frac{x_8}{1-x_8}\cdot\frac{1-x_4}{x_4}},\\
\dis{g_5(\tb{x})=\frac{x_3}{1-x_3}\cdot\frac{1-x_5}{x_5}},\ \dis{g_6(\tb{x})=\frac{x_2}{1-x_2}\cdot\frac{1-x_6}{x_6}},
\end{array}
\end{displaymath}
(These come from rows (3)(d), (4)(d), (5)(b), (5)(c), (5)(e) and (5)(f) of Table \ref{table1}), for some $\tb{x}=(x_1,\dots,x_8)\in\dseven=\{\tb{x}\in\R^8_+|\sum_{i=1}^8x_i=1\}$. Under the same hypothesis on $\Gamma$, for any $z_0\in\hyp$ and for each $\gamma\in\{\xi^2,\eta^{-2},\xi\eta^{-1},\eta\xi\}$, the expression $e^{2\dgamma}$ is bounded below by $h_i(\tb{x})$ for at least one of $h_i$ from the list
\begin{displaymath}
\begin{array}{l}
\dis{h_1(\tb{x})=\frac{1-x_6}{x_6}\cdot\frac{1-x_8}{x_8}},\  \dis{h_2(\tb{x})=\frac{1-x_2}{x_2}\cdot\frac{1-x_8}{x_8}},\\
\dis{h_3(\tb{x}) =\frac{1-x_5}{x_5}\cdot\frac{1-x_7}{x_7}},\  \dis{h_4(\tb{x})=\frac{1-x_3}{x_3}\cdot\frac{1-x_7}{x_7}},
\end{array}
\end{displaymath}
for some $i\in \{1,2,3,4\}$ and $\tb{x}\in\dseven$ (The functions $h_i$ are produced from rows (6)(a)-(d) of Table \ref{table1}).
\end{proposition}
\begin{proof}
By Lemma \ref{lem3.2} ($1$), we have $0\leq\nu_{\psi}(S_{\infty})\leq 1$ for every $\psi\in\Psi^{\dagger}$. We aim to show that $0<\nu_{\psi}(S_{\infty})<1$ for any $\psi\in\Psi^{\dagger}$. 
First assume on the contrary that $\nu_{\xi^{-1}}(S_{\infty})=0$. Applying Theorem \ref{thm3.4} (2) to row (6)(a) of Table \ref{table1} implies that $\nu_{\eta^{-1}\xi}(S_{\infty})=1$. By Theorem \ref{thm3.4} (1), we see that $\nu_{\psi}(S_{\infty})=0$ for every $\psi\in\Psi^{\dagger}- \{\eta^{-1}\xi\}$. Using the fact that $\nu_{\eta}(S_{\infty})=0$ and applying Theorem \ref{thm3.4} (2) to row (6)(b) of Table \ref{table1} shows that $\nu_{\xi\eta^{-1}}(S_{\infty})=1$, a contradiction. A similar argument can be repeated for $\nu_{\eta}(S_{\infty})$ by exchanging the roles of $\xi^{-1}$ and $\eta$ above. Therefore, we have $\nu_{\xi^{-1}}(S_{\infty})\neq 0$ and $\nu_{\eta}(S_{\infty})\neq 0$.

Assume that $\nu_{\psi_0}(S_{\infty})=0$ for a given $\psi_0\in \{\xi\eta,\xi\eta^{-1},\xi^2,\eta^{-1}\xi,\eta^{-1}\xi^{-1},\eta^{-2}\}$. Consider the following lists
$$
\begin{array}{c}
( \xi^{-1}, \xi\eta,  \Psi^{\dagger}-\{\eta\}, \xi^{-1}),\quad ( \xi^{-1}, \xi\eta^{-1}, \Psi^{\dagger}-\{\eta^{-1}\xi,\eta^{-1}\xi^{-1},\eta^{-2}\},\eta),  \\
(  \xi^{-1}, \xi^2, \Psi^{\dagger}-\{\xi\}, \xi^{-1}),\quad ( \eta, \eta^{-2}, \Psi^{\dagger}-\{\eta^{-1}\xi,\eta^{-1}\xi^{-1},\eta^{-2}\},\xi^{-1}),  \\
(  \eta,  \eta^{-1}\xi^{-1}, \Psi^{\dagger}-\{\xi^{-1}\}, \eta), \quad ( \eta, \eta^{-1}\xi, \Psi^{\dagger}-\{\xi\eta^{-1},\xi\eta,\xi^{2}\},\xi^{-1}),
\end{array}
$$
where each entry in a list is assigned for $\gamma_0$, $\psi_0$,  $S(\gamma_0)$, $\psi_1$, respectively. By applying Theorem \ref{thm3.4} (2) to Table \ref{table1} with $\psi_0=s(\gamma_0)$, we get $\sum_{\psi\in S(\gamma_0)}\nu_{\psi}=1$.
We have $\psi_1\notin S(\gamma_0)$. Therefore, we obtain $\nu_{\psi_1}(S_{\infty})=0$ for some $\psi_1\in\{\xi^{-1},\eta\}$, a contradiction.
As a result, we conclude that $0<\nu_{\psi}(S_{\infty})<1$ for every
$\psi\in\Psi^{\dagger}$. Since we have
$\psi=s(\gamma)$ for some $\gamma$ in Table \ref{table1} and $S(\gamma)\subset\Psi^{\dagger}$, we also conclude
$$0<\int_{S_{\infty}}(\lambda_{\gamma,z_0})^{2}d\mu_{V_{s(\gamma)}}  =  1-\sum\nolimits_{\psi\in S(\gamma)}\int_{S_{\infty}} d\nu_{\psi}<1$$
by Theorem \ref{thm3.4} (2). In other words, $\nu_{s(\gamma)}$ and $\int_{S_{\infty}}\lambda^{2}_{\gamma,z_0}d\mu_{V_{s(\gamma)}}$ satisfy the hypothesis of Lemma \ref{CSlem5.5} for every $\gamma$ in Table \ref{table1}.

We apply Lemma \ref{CSlem5.5} to every row of Table \ref{table1} with $\nu=\nu_{s(\gamma)}$, $a=\nu_{s(\gamma)}(S_{\infty})$ and
$b=\int_{S_{\infty}}\lambda^{2}_{\gamma,z_0}d\mu_{V_{s(\gamma)}}$. 
Using Theorem \ref{thm3.4} (2), we calculate the lower bounds as
\begin{eqnarray}
e^{2\dgamma}                & \geq    & \frac{\left(1-\sum_{\psi\in
S(\gamma)}m_{p(\psi)}\right)\cdot\left(1-m_{p(s(\gamma))}\right)}{\left(\sum_{\psi\in S(\gamma)}m_{p(\psi)}\right)\cdot
              m_{p(s(\gamma))}},\label{25}
\end{eqnarray}
where $\int_{S_{\infty}} d\nu_{\psi}=m_{p(s(\gamma))}$ for the bijection $p$ in (\ref{sigma}).
Upon replacing each constant $m_{p(s(\gamma))}$ appearing in (\ref{25}) with the variable $x_{p(s(\gamma))}$ we obtain the functions listed in the proposition.
\end{proof}

Note that we have $18$ lower bounds given in
the expression (\ref{25}) for the displacements under the isometries in $\Psi^{\dagger}_r\cup\Psi^{\dagger}$ 
because, there are $18$ group theoretical relations listed in (\ref{113}). Since we are interested in the displacements under the isometries in $\Gamma_{\dagger}=\{\xi,\eta,\xi\eta\}$, we
will concentrate on the first $14$ displacement functions $f_1,f_2\dots,f_8,g_1, g_2,\dots,g_6$ in Proposition \ref{dispfunc} for the proofs of Lemmas \ref{lem4.1}, \ref{lem4-2}, \ref{lem4.1-1}, \ref{lem4.1-2}, \ref{lem4.1-3}, \ref{lem4.1-4}, \ref{lem4.3}, \ref{lem4.4}, \ref{lem4.5}, \ref{lem4.9} and  Theorems \ref{thm4.14}, \ref{thm4.1} and \ref{thm5.1}.


\section{Lower bound for $\max\{\dgammaz\co\gamma\in\Gamma_{\dagger}\}$ when $\Lambda_{\Gamma\cdot z}=S^{2}$}\label{sec4}

Let $\mathcal{F}^{\dagger}=\{f_1,f_2,\dots,f_8,g_1,g_2,\dots,g_6\}$.
The constants on the right hand side of the inequalities in (\ref{25}) can be considered as the values of the functions in $\mathcal{F}^{\dagger}$ at the point $\tb{m}_{\dagger}=(m_1,m_2,\dots,m_8)\in\dseven=\left\{\left(x_1,x_2,x_3,x_4,x_5,x_6,x_7,x_8\right)\in\mathbb{R}^8_+\co\sum_{i=1}^{8}x_i=1\right\}$.

When $\Gamma=\langle\xi,\eta\rangle$ is geometrically infinite, the lower bound given in Theorem \ref{thm5.1} for the displacements under the isometries in $\Gamma_{\dagger}=\{\xi,\eta,\xi\eta\}$ follows from the calculation of the infimum of the maximum of the functions in $\mathcal{F}^{\dagger}$ over the simplex $\dseven$. Therefore, in this section, we aim to prove the statement below:
\begin{theorem}\label{thm4.1}
Let $G^{\dagger}\co  \Delta^7 \to \mathbb{R}$ be the function defined by $\tb{x}\mapsto\max\{f(\tb{x})\co f\in\mathcal{F}^{\dagger}\}$. Then $\inf_{\tb{x}\in\Delta^7}G^{\dagger}(\tb{x})=5+3\sqrt{2}$.
\end{theorem}

To this purpose, we shall show that it is enough to calculate the infimum of the maximum of the first eight $f_1,f_2,\dots,f_8$ of the displacement functions in $\mathcal{F}^{\dagger}$. Let $I^{\dagger}=\{1,2,3,4,5,6,7,8\}$. Then we first state the following: 
\begin{lemma}\label{lem4.1}
Let $F^{\dagger}\co\dseven\to\mathbb{R}$ be the function defined by $\tb{x}\mapsto\max\{f_i(\tb{x}): i\in I^{\dagger}\}$.
Then $\alphas=\inf_{\tb{x}\in\Delta^7}F^{\dagger}(\tb{x})$ is attained in $\dseven$ and satisfies $9\leq\alphas\leq 5+3\sqrt{2}$.
\end{lemma}
\begin{proof}
It is clear that
$\inf_{\tb{x}\in\Delta^7}\{\max (f_7(\tb{x}),f_8(\tb{x}))\}\leq\inf_{\tb{x}\in\Delta^7}F^{\dagger}(\tb{x})$. We apply the substitution
 $X_1  =  x_7$, $X_2  =  x_8$, $X_3  =  x_1+x_2+x_3$, $X_4  =  x_4+x_5+x_6$.
Then we see that
$f_7(\tb{x})  =  \sigma(X_4)\sigma(X_1),\ \tnr{and}\ \ f_8(\tb{x})  =  \sigma(X_3)\sigma(X_2)$,
where $\sum_{i=1}^4 X_i=1$ and $\sigma(x)=1/x-1$ for $x\in (0,1)$. By Lemma \ref{lem2.1}, 
we obtain that $9\leq\alphas$. 

Let $\Delta_n=\{\tb{x}\in\Delta^7\co 1/n\leq x_i\leq 1-1/(2n)\ \tr{for}\ i\in I^{\dagger}\}$ of $\Delta^7$ for every
$n\geq 2$. Note that $\Delta_{n+1}\subset\Delta_n$. The function $F^{\dagger}$ has an absolute minimum $F^{\dagger}(\tb{x}_n)$ at some point $\tb{x}_n\in\Delta_n$. 
The sequence $(F^{\dagger}(\tb{x}_n))_{n=1}^{\infty}$ tends to an infimum because it is a decreasing sequence which is bounded below by $9$.

We claim that the sequence $(\tb{x}_n)_{n=1}^{\infty}$ cannot have a limit point on the boundary of $\dseven$. Assume on the contrary that
$\tb{x}_n\to \tb{b}\in\partial\dseven$ as $n\to \infty$. If $(b_1,b_2,\dots,b_8)$ denotes the coordinates of the point $\tb{b}$, then $b_i=0$ for some
$i\in I^{\dagger}$. Let us assume that $b_i=0$ for some $i\in\{1,7\}$. Then using the function $f_i$, we conclude that $b_4+b_5+b_6=1$. Because,
otherwise $(f_i(\tb{x}_n))_{n=1}^{\infty}$ would tend to infinity. But it is supposed to be tending to an infimum of $F^{\dagger}$. Therefore, we must have $b_j=0$
for every $j\in\{2,3,8\}$. Then, we get that $b_4+b_5+b_6+b_7+b_8=1$ and $b_1+b_2+b_3+b_7+b_8=1$ and $b_1+b_2+b_3=1$. Because, otherwise $(f_j(\tb{x}_n))_{n=1}^{\infty}$ would tend to
infinity when it is supposed to tend to an infimum of $F^{\dagger}$. In any case, we obtain that $b_4+b_5+b_6+b_i>1$ for some
$i\in\{1,2,3,7,8\}$. This is a contradiction. Therefore, $b_1\neq 0$ and $b_7\neq 0$. Similar arguments with suitably chosen displacement functions show that $b_i\neq 0$ for every $i\in\{2,3,4,5,6,8\}$. Hence, the sequence $(\tb{x}_n)_{n=1}^{\infty}$ cannot have a limit point on the boundary of
$\dseven$. Then there exists a positive integer $n_0$ so that $\tb{x}_n=\tb{x}_{n+1}$ for every $n\geq n_0$. Otherwise we would
have a subsequence $(\tb{x}_{n_j})_{j=1}^{\infty}$ of the sequence $(\tb{x}_n)_{n=1}^{\infty}$ so that $\tb{x}_{n_j}\in\Delta_{n_{j+1}}-\Delta_{n_j}$ for every integer $j>0$. Since we
have $\cup_{n=2}^{\infty}\Delta_n=\Delta^7$, there exists a subsequence of  $(\tb{x}_{n_j})_{j=1}^{\infty}$ which has a limit point on $\partial\Delta^7$, a contradiction. As a result, $\inf_{\tb{x}\in\dseven}F^{\dagger}(\tb{x})$ is attained at some point in $\dseven$, i.e., $\alphas=\min_{\tb{x}\in\dseven}F^{\dagger}(\tb{x})$.

Let $x_i=(\sqrt{2}-1)/2$ for $i=1,4,7,8$ and
$x_i=(3-2\sqrt{2})/4$ for $i=2,3,5,6$. Then
$\bar{\tb{x}}=(x_i)_{i\in I^{\dagger}}$
is a point in $\Delta^7$ such that $f_i(\bar{\tb{x}})=5+3\sqrt{2}$ for every $i\in I^{\dagger}$. Therefore, we get $F^{\dagger}(\bar{\tb{x}})=5+3\sqrt{2}\geq\alphas$.
\end{proof}

In the rest of this section, we will consider $\dseven$ as a submanifold of $\mathbb{R}^8$.
The tangent space $T_{\tb{x}}\Delta^7$ at any $\tb{x}\in\Delta^7$ consists of vectors whose coordinates sum to $0$. Note that each displacement function $f_i$ for $i\in I^{\dagger}$ is smooth in an open neighborhood of $\Delta^7$. Therefore, the directional derivative of $f_i$ in the direction of any $\vec v\in T_{\tb{x}}\Delta^7$ is given by $\nabla f_i(\tb{x})\cdot\vec v$ for any $i\in I^{\dagger}$.

We shall use the identity $\sum_{i=1}^{8}x_i=1$ to rewrite the formulas of the functions $f_i$ given in Proposition \ref{dispfunc} in various ways in the proofs of lemmas below to suit our purposes. Although they do not take the same values on all of $\R^8$, we will abuse notation and call the rewritten functions by $f_i$, which agree with the originals on $\dseven$.

\subsection{Relationships between the displacement functions $f_1,f_2,\dots,f_8$}\label{sec411}

By Lemma \ref{lem4.1}, we know that $\alphas$ is attained by a displacement function $f_i$ for some $i\in I^{\dagger}$.  In fact, it is possible to see that more than one function in $\{f_1,\dots,f_8\}$ attain the value $\alphas$. In other words, we have
\begin{lemma}\label{lem4-2}
Let $\alphas=\inf_{\tb{x}\in\dseven}\max\{f_1(\tb{x}),\dots,f_8(\tb{x})\}$, where $f_i$ for $i\in I^{\dagger}$ are as in Proposition \ref{dispfunc}. At any $\xs\in\dseven$ such that $F^{\dagger}(\xs)=\alphas$, there exist at least two functions $f_{i},f_{j}$
such that $f_{i}(\tb{x}^*)=f_{j}(\tb{x}^*)$ for $i\neq j$, where $i,j\in I^{\dagger}=\{1,2,3,4,5,6,7,8\}$.
\end{lemma}
\begin{proof}
Observe that for each function $f_i$ for $i\in I^{\dagger}$ there is a variable $x_j$ such that the first order partial derivative of $f_i$ with respect to $x_j$ at $\tb{x}$ is $0$ for every $\tb{x}\in\dseven$. But the first order partial derivatives of $f_i$ with respect to $x_i$ are strictly negative at every $\tb{x}\in\dseven$. These facts imply that $\nabla f_i$ is not a scalar multiple of the perpendicular $\langle 1,1,\dots,1\rangle$ to $T_{\tb{x}}\dseven$ for any $i\in I^{\dagger}$. Therefore, none of the functions $f_1,f_2,\dots,f_8$ has a local extremum on $\Delta^{7}$.

If $f_i(\tb{x}^*)\neq f_j(\tb{x}^*)$ for every $i\neq j$, then
the set $\{f_1(\tb{x}^*), f_2(\tb{x}^*),\dots, f_8(\tb{x}^*)\}$ has a unique largest element. By renumbering the functions, we may assume
that $f_1(\tb{x}^*)$ is the largest value, i.e., $f_1(\tb{x}^*)=\alphas$. By the continuity of $\mathcal{F}^{\dagger}$, there exists a neighborhood $U$ of $\tb{x}^*$
contained in $\Delta^{7}$ so that $F^{\dagger}(\tb{x})=f_1(\tb{x})$ for every $\tb{x}\in U$. Since $F^{\dagger}$ has a minimum at $\tb{x}^*$, then $f_1$ must have a local minimum at
$\tb{x}^*$, a contradiction. The lemma follows.
\end{proof}

Next, we will consider the cases in which $f_{i}$ and $f_{j}$ in Lemma \ref{lem4-2} are in the sets $\{f_1,f_2,f_3\}$, $\{f_4,f_5,f_6\}$, and $\{f_7,f_8\}$, respectively:
\begin{lemma}\label{lem4.1-1}
Let $\alphas=\inf_{\tb{x}\in\dseven}\max\{f_1(\tb{x}),\dots,f_8(\tb{x})\}$, where $f_i$ for $i\in I^{\dagger}$ are as in Proposition \ref{dispfunc}. At any $\xs\in\dseven$ such that $F^{\dagger}(\xs)=\alphas$, we have either
\begin{itemize}
\item[(1)] $f_l(\xs)=\alphas$ for all $l\in I_1=\{1,2,3\}$ or
\item[(2)] $f_l(\tb{z})<\alphas$ for all $l\in I_1$ and $f_j(\tb{z})=f_j(\xs)$ for all $j\in\mathcal{I}_1=\{4,5,6,7,8\}$ for some $\tb{z}\in\dseven$ such that $F^{\dagger}(\tb{z})=F^{\dagger}(\xs)$.
\end{itemize}
\end{lemma}
\begin{proof}
Assume that part ($1$) of the lemma does not hold at $\xs$. If $f_i(\xs)<\alphas$ for every $i\in I_1$, the point $\tb{z}=\xs$ satisfies part (2).  Then it is enough to consider the case 
$
f_i(\tb{x}^*)<f_j(\tb{x}^*)\leq f_k(\tb{x}^*)=\alphas,$ and $\ f_l(\tb{x}^*)\leq\alphas
$
for $l\in\mathcal{I}_{1}$, where $i,j,k\in I_1$ such that $i\neq j$, $j\neq k$, $i\neq k$.
Let us define the vectors $\vec u_1^2$, $\vec u_1^3$ and $\vec u_2^3$ as $\langle -1, 1, 0,0,0,0,0,0\rangle$, $\langle -1, 0, 1,0,0,0,0,0\rangle$ and $\langle 0, -1, 1,0,0,0,0,0\rangle$ in $T_{\xs}\dseven$, respectively. Also let $\vec u_2^1=-\vec u_1^2$, $\vec u_3^1=-\vec u_1^3$ and $\vec u_3^2=-\vec u_2^3$.

Using the identity $x_k=1-\sum_{n=1,n\neq k}^8x_n$,
we calculate the directional derivatives of all of the functions $f_1,f_2,\dots,f_8$ in the direction of the vector $\vec u_i^{j}$. Note that none of the functions $f_4$, $f_5$,\dots, $f_8$ contains the variables $x_1$, $x_2$ or $x_3$. For every $\tb{x}\in\Delta^{7}$ and for every $l\in\mathcal{I}_1$ we see that 
$\nabla f_i(\tb{x})\cdot\vec u_i^j>0$, $\nabla f_j(\tb{x})\cdot\vec u_i^j<0$, $\nabla f_k(\tb{x})\cdot\vec u_i^j=0$, $\nabla f_l(\tb{x})\cdot\vec u_i^j=0$,
which implies that the values of $f_j$ and $f_k$ decrease along a line segment in the direction of $\vec v=\vec u_i^j+\vec u_i^k$. The values of $f_l$ are constant along this segment, and for a short distance along $\vec v$ the values of $f_i$ is smaller than those of $f_j$ and $f_k$. Therefore there exists a point $\tb{z}$ on this line segment satisfying part (2) of the lemma.
\end{proof}

Analogous results hold for the displacement functions in $\{f_4,f_5,f_6\}$ and $\{f_7,f_8\}$. In particular, we have the followings: 
\begin{lemma}\label{lem4.1-2}
Let $\alphas=\inf_{\tb{x}\in\dseven}\max\{f_1(\tb{x}),\dots,f_8(\tb{x})\}$, where $f_i$ for $i\in I^{\dagger}$ are as in Proposition \ref{dispfunc}. At any $\xs\in\dseven$ such that $F^{\dagger}(\xs)=\alphas$, we have either
\begin{itemize}
\item[(1)] $f_l(\xs)=\alphas$ for all $l\in I_2=\{4,5,6\}$ or
\item[(2)] $f_l(\tb{z})<\alphas$ for all $l\in I_2$ and $f_j(\tb{z})=f_j(\xs)$ for all $j\in\mathcal{I}_2=\{1,2,3,7,8\}$ for some $\tb{z}\in\dseven$ such that $F^{\dagger}(\tb{z})=F^{\dagger}(\xs)$.
\end{itemize}
\end{lemma}
\begin{lemma}\label{lem4.1-3}
Let $\alphas=\inf_{\tb{x}\in\dseven}\max\{f_1(\tb{x}),\dots,f_8(\tb{x})\}$, where $f_i$ for $i\in I^{\dagger}$ are as in Proposition \ref{dispfunc}. At any $\xs\in\dseven$ such that $F^{\dagger}(\xs)=\alphas$, we have either
\begin{itemize}
\item[(1)] $f_l(\xs)=\alphas$ for all $l\in I_3=\{7,8\}$ or
\item[(2)] $f_l(\tb{z})<\alphas$ for all $l\in I_3$ and $f_j(\tb{z})=f_j(\xs)$ for all $j\in\mathcal{I}_3=\{1,2,3,4,5,6\}$ for some $\tb{z}\in\dseven$ such that $F^{\dagger}(\tb{z})=F^{\dagger}(\xs)$.
\end{itemize}
\end{lemma}
The proof of Lemma \ref{lem4.1-1} applies, mutatis mutandis, to prove Lemma \ref{lem4.1-2} and \ref{lem4.1-3}. In particular, using each identity $x_k=1-\sum_{n=1,n\neq k}^8x_i$ for $k\in I_2\cup I_3$, we perturb in the directions of the vectors $\vec u_4^5=\langle 0,0,0,-1,1,0,0,0\rangle$, $\vec u_4^6=\langle 0,0,0,-1,0,1,0,0\rangle$, and $\vec u_5^6=\langle 0,0,0,0,-1,1,0,0\rangle$ for Lemma \ref{lem4.1-2}, and perturb in the direction of the vectors $\vec u_7^8=\langle 0,0,0,0,0,0,-1,1\rangle$ and $\vec u_8^7=\langle 0,0,0,0,0,0,1,-1\rangle$ for  Lemma \ref{lem4.1-3}. Lemmas \ref{lem4.1-1}, \ref{lem4.1-2} and \ref{lem4.1-3} imply the following:
\begin{lemma}\label{lem4.1-4}
Let $\alphas=\inf_{\tb{x}\in\dseven}\max\{f_1(\tb{x}),\dots,f_8(\tb{x})\}$, where $f_i$ for $i\in I^{\dagger}$ are as in Proposition \ref{dispfunc}. There exists a point $\xs\in\dseven$ which satisfies one of the cases I, II, III, IV or V, where
\begin{itemize}
\item[I.] $f_j(\tb{x}^*)<\alphas$ for $j\in\mathcal{I}_1=\{4,5,6,7,8\}$ and $f_j(\tb{x}^*)=\alphas$ for $j\in I_1=\{1,2,3\}$,

\item[II.] $f_j(\tb{x}^*)<\alphas$ for $j\in\mathcal{I}_3=\{1,2,3,4,5,6\}$ and $f_j(\tb{x}^*)=\alphas$ for $j\in I_3=\{7,8\}$,

\item[III.] $f_j(\tb{x}^*)=\alphas$ for $j\in\mathcal{I}_3=\{1,2,3,4,5,6\}$ and $f_j(\tb{x}^*)<\alphas$ for $j\in I_3=\{7,8\}$,

\item[IV.] $f_j(\tb{x}^*)=\alphas$ for $j\in\mathcal{I}_2=\{1,2,3,7,8\}$ and $f_j(\tb{x}^*)<\alphas$ for $j\in I_2=\{4,5,6\}$,

\item[V.] $f_j(\xs)=\alphas$ for $j\in I^{\dagger}=\{1,2,3,4,5,6,7,8\}$.
\end{itemize}
\end{lemma}
\begin{proof}
Let $\tb{x}\in\dseven$ be a point such that $F^{\dagger}(\tb{x})=\alphas$. First, assume that $f_i(\tb{x})<\alphas$ for some $i\in I_1$. By Lemma \ref{lem4.1-1}, there exists a point $\tb{z}_1\in\dseven$ with $f_i(\tb{z}_1)<\alphas$ for all $i\in I_1$ and $f_j(\tb{z}_1)=\alphas$ for all $i\in \mathcal{I}_1$ with $F^{\dagger}(\tb{z}_1)=\alphas$.

If $f_j(\tb{z}_1)<\alphas$ for some $j\in I_2$, there exists a point $\tb{z}_2\in\dseven$ with $f_i(\tb{z}_1)=f_i(\tb{z}_2)<\alphas$ for all $i\in I_1$, $f_j(\tb{z}_2)<\alphas$ for all $j\in I_2$, $f_k(\tb{z}_1)=f_k(\tb{z}_2)$ for all $i\in I_3$ and $F^{\dagger}(\tb{z}_2)=\alphas$ by Lemma \ref{lem4.1-2}. We must have $f_k(\tb{z}_2)=\alphas$ for all $i\in I_3$ by Lemma \ref{lem4-2}. Thus $\xs=\tb{z}_2$ satisfies Case II. Assume that $f_j(\tb{z}_1)=\alphas$ for all $j\in I_2$. Let $T_1\co\dseven\to\dseven$ be the transformation
\begin{equation}\label{transf}
x_1\mapsto x_4,\ x_2\mapsto x_5,\ x_3\mapsto x_6,\ x_4\mapsto x_1,\ x_5\mapsto x_2,\ x_6\mapsto x_3,\ x_7\mapsto x_8\, x_8\mapsto x_7.
\end{equation}
If $f_k(\tb{z}_1)=\alphas$ for all $k\in I_3$, then $\xs=T_1(\tb{z}_1)$ satisfies Case IV. Otherwise, there exists a point $\tb{z}_2\in\dseven$ such that $f_i(\tb{z}_2)=f_i(\tb{z}_1)$ for $i\in I_1$,  $f_j(\tb{z}_2)=f_j(\tb{z}_1)$ for $j\in I_2$ and $f_k(\tb{z}_2)<\alphas$ for $k\in I_3$ by Lemma \ref{lem4.1-3}. Thus, $\xs=T_1(\tb{z}_2)$ satisfies Case I.

Consider the case $f_i(\tb{x})=\alphas$ for all $i\in I_1$. If also $f_j(\tb{x})=\alphas$ for all $j\in I_2$, then either $\xs=\tb{x}$ satisfies Case V or there exists a point $\xs=\tb{z}_1$ obtained by Lemma \ref{lem4.1-3} satisfying Case III. Therefore, assume that $f_j(\tb{x})<\alphas$ for some $j\in I_2$. By Lemma \ref{lem4.1-2}, there exists a point $\tb{z}_2\in\dseven$ with the property that $f_i(\tb{z}_1)=\alphas$ for all $i\in I_1$, $f_j(\tb{z}_1)<\alphas$ for all $j\in I_2$ and $f_k(\tb{z}_1)=f_k(\tb{x})$ for all $k\in I_3$. Then either $\xs=\tb{z}_1$ satisfies Case IV, or there exists a point $\xs=\tb{z}_2\in\dseven$ satisfying Case I by Lemma \ref{lem4.1-3}.
\end{proof}

\subsection{Calculations of the infimums}\label{sec416}

Let $\Delta_7=\left\{\tb{y}\in\Delta^7\co f_j (\tb{y})=f_k (\tb{y}),j,k\in I^{\dagger}\right\}\subset\Delta^7$. Note that $\bar{\tb{x}}\in\Delta_7$ (see Lemma \ref{lem4.1}). We aim to prove the following proposition

\begin{proposition}\label{lem4.2}
The infimum $\alphas=\min_{\tb{x}\in\dseven}F^{\dagger}(\tb{x})\in[9,5+3\sqrt{2}]$ is attained at some point $\xs\in\Delta_7$.
\end{proposition}

To this purpose, we need to show that Cases I, II, III and IV are not possible at a point $\tb{x}\in\dseven$ so that $F^{\dagger}(\tb{x})=\alphas$. We start with Case I.

\begin{lemma}\label{lem4.3}
Let $\alphas=\inf_{\tb{x}\in\dseven}\max\{f_1(\tb{x}),\dots,f_8(\tb{x})\}$, where $f_i$ for $i\in I^{\dagger}$ are as in Proposition \ref{dispfunc}. At any $\xs\in\dseven$ satisfying $f_{i} (\xs)=f_j (\xs)$ for every $i,j\in\{1,2,3\} $, and $f_2 (\xs)>f_6 (\xs)$, there exists $\vec v\in T_{\xs}\dseven$ such that each of $f_1 $, $f_2 $, and $f_3 $ decreases in the direction of $\vec v$.
\end{lemma}
\begin{proof}
Using the identity $x_8=1-\sum_{n=1}^7x_n$,  we rewrite $f_1 $, $f_2 $ and $f_3 $ as follows:
\[
\begin{array}{lll}
f_1(\tb{x})=\sigma(\Sigma_2(\tb{x}))\sigma(x_1), & f_2(\tb{x})=\dis{\frac{\sigma(x_2)}{\sigma(\Sigma_1(\tb{x}))}}, & f_3(\tb{x})=\dis{\frac{\sigma(x_3)}{\sigma(\Sigma_2(\tb{x}))}}.
\end{array}
\]  
These functions are each well-defined and smooth on an open neighborhood of $\dseven$ in $\mathbb{R}^8$. Because $\sigma(x)$ decreases in $x$ and $\Sigma_2$ is constant in all variables but $x_4$, $x_5$ and $x_6$, the following facts are clear:
\begin{itemize}
\item[(1)] $\dis{\frac{\partial f_1 }{\partial x_1}}\bigg|_{\xs}=-\frac{\sigma(\Sigma_2(\xs))}{(x_1^*)^2}<0$, $\dis{\frac{\partial f_1 }{\partial x_2}\bigg|_{\xs}=0}$, $\dis{\frac{\partial f_1 }{\partial x_3}\bigg|_{\xs}=0}$,
\item[(2)] $\dis{\frac{\partial f_3 }{\partial x_3}}\bigg|_{\xs}=-\frac{1}{(x_3^*)^2\sigma(\Sigma_2(\xs))}<0$, $\dis{\frac{\partial f_3 }{\partial x_1}\bigg|_{\xs}=0}$, $\dis{\frac{\partial f_3 }{\partial x_2}(\xs)=0}$,
\item[(3)]$\dis{\frac{\partial f_1 }{\partial x_k}}\bigg|_{\xs}=0$, $\dis{\frac{\partial f_2 }{\partial x_k}}\bigg|_{\xs}=0$ and $\dis{\frac{\partial f_3 }{\partial x_k}}\bigg|_{\xs}=0$
\end{itemize}
for every $k\in I_3 $. These facts imply that at any $\xs\in\dseven$ such that the equation below
\begin{equation}\label{eqn8}
\frac{\partial f_2 }{\partial x_2}\bigg|_{\xs}=0
\end{equation}
does not hold, there exists some $r\in\mathbb{R}$ such that each of $f_1 $, $f_2 $, and $f_3 $ decreases in the direction of the vector $\vec v_r=\langle 1,r,1,0,0,0,-r-2,0\rangle$. Note that $\vec v_r\in T_{\xs}\dseven$ since its coordinates sum to $0$. Thus it only remains to consider the case in which the equality in (\ref{eqn8}) holds.

A computation gives that
\begin{equation*}
\frac{\partial f_2 }{\partial x_2}\bigg|_{\xs}=\frac{\left(x_2^*-\Sigma_1(\xs)\right)\left(1-x_2^*-\Sigma_1(\xs)\right)}{(x_2^*)^2(\Sigma_1(\xs))^2},
\end{equation*}
which vanishes if and only if $\Sigma_2(\xs)+x_2^*=1$. Since $\Sigma_1(\xs)=x_1^*+x_2^*+x_3^*$, we conclude that equation in (\ref{eqn8}) holds if and only if $x_2^*  =  (1-x_1^*-x_3^*)/2$. By the identity $\sum_{i=1}^8x_i^*=1$, this is in turn equivalent to $x_2^*=x_4^*+x_5^*+x_6^*+x_7^*+x_8^*$. Therefore, we find that $x_6^*<x_2^*$. Then the lemma follows, because, by the definitions of $f_2 $ and $f_6 $, we obtain
$f_2 (\xs)=\left(\sigma(x_2^*)\right)^2<\sigma(x_2^*)\sigma(x_6^*)=f_6 (\xs),$
a contradiction.
\end{proof}
Before we proceed to Cases II, III and IV, we shall first prove the following statement:
\begin{lemma}\label{flipside}
For $1\leq k\leq n-1$, let $f_1$, $f_2$,\dots, $f_k$ be smooth functions on an open neighborhood $U$ of the $(n-1)-$simplex $\Delta^{n-1}$ in $\mathbb{R}^n$. If at some $\tb{x}\in\Delta^{n-1}$ the collection $\{\nabla f_1(\tb{x}), \nabla f_2(\tb{x}),\dots,\nabla f_k(\tb{x}), \langle 1,1,\dots,1\rangle\}$ of vectors in $\mathbb{R}^n$ is linearly independent, then there exists a vector $\vec u\in T_{\tb{x}}\Delta^{n-1}$ such that each $f_i$ for $i=1,\dots,k$ decreases in the direction of $\vec u$ at $\tb{x}$.
\end{lemma}
\begin{proof}
Let $\mathcal{B}=\{\vec v_1,\dots,\vec v_n\}$ be a collection of $n$
linearly independent vectors in $\mathbb{R}^n$. We claim that there exists a vector $\vec u\in\mathbb{R}^n$ such that $\vec u\cdot\vec v<0$ for every
$\vec v\in\mathcal{B}$. The assertion is clear for $n=1$. For $n>1$, assume that there exists a vector $\vec u_0\in Span\{\vec v_1,\dots,\vec v_{n-1}\}$ such that $\vec u_0\cdot\vec v<0$ for every $v\in \mathcal{S}=\{\vec v_1,\dots,\vec v_{n-1}\}$ by induction.

There is a nonzero vector $\vec v_0\in\R^n$ orthogonal to each vector in $\mathcal{S}$. If we have $\vec v_0\cdot\vec v_n=0$, then $\vec v_n$ is in the space $\vec v_0^{\perp}$ of vectors perpendicular to $\vec v_0$. Since $\tnr{dim}\ \vec v_0^{\perp}=n-1$, the set $\mathcal{S}$ spans $\vec v_0^{\perp}$. The set $\mathcal{B}$ is linearly independent therefore, we get $\vec v_n\cdot\vec v_0\neq 0$. Let $\vec u=\vec u_0-c\vec v_0$ for $c=(\vec u_0\cdot\vec v_n+1)/\vec v_0\cdot\vec v_n$. Then we see that $\vec u\cdot\vec v<0$ for every $\vec v\in\mathcal{B}$, which proves the claim.

Let $\vec w=(1,1,\dots,1)\in\R^n$. Complete the set $\{\nabla f_1(\tb{x}), \nabla
f_2(\tb{x}),\dots,\nabla f_k(\tb{x}),\vec w\}$ to a basis $\mathcal{B}=\{\nabla f_1(\tb{x}), \nabla
f_2(\tb{x}),\dots,\nabla f_k(\tb{x}),\vec u_{k+1},\dots,\vec u_{n-1},\vec w\}$ for $\mathbb{R}^n$. If we declare $\vec v_i=\tnr{proj}_{\vec w^{\perp}}\nabla f_i(\tb{x})$ for $i=1,\dots,k$ and $\vec v_j=\tnr{proj}_{\vec w^{\perp}}\vec u_j$ for $j=k+1,\dots,n-1$, then $\{\vec v_1,\dots,\vec v_{n-1},\vec w\}$ is linearly independent. This is because $\mathcal{B}$ is linearly independent. Let $\mathcal{S}=\{\vec v_1,\dots,\vec v_{n-1}\}$. Since $\vec w^{\perp}$ has dimension $n-1$ and $\mathcal{S}$ is linearly independent, we have $Span\ \mathcal{S}=\vec w^{\perp}$. By the fact above, there exists a vector $\vec u\in\ Span\ \mathcal{S}$ so that $\vec u\cdot\vec v<0$ for every $\vec v\in\mathcal{S}$. In particular, we get $\vec u\cdot\vec v_i=\vec u\cdot\nabla f_i(\tb{x})<0$ for $i=1,\dots,k$. Since $T_{\tb{x}}\Delta^{n-1}$ consists of vectors whose entries sum to $0$,
we have $\vec w^{\perp}=T_{\tb{x}}\Delta^{n-1}$, which completes the proof.
\end{proof}

The lemmas \ref{lem4.4}, \ref{lem4.5} and \ref{lem4.9} below show respectively that Cases II, III and IV are not possible at a point at which $F^{\dagger}$ takes it minimum value:
\begin{lemma}\label{lem4.4}
Let $\alphas=\inf_{\tb{x}\in\dseven}\max\{f_1(\tb{x}),\dots,f_8(\tb{x})\}$, where $f_i$ for $i\in I^{\dagger}$ are as in Proposition \ref{dispfunc}. At any $\xs\in\dseven$ such that $f_7 (\xs)=f_8 (\xs)$ and $f_4 (\xs)<f_8 (\xs)$, there exists a vector $\vec v\in T_{\xs}\dseven$ so that $f_7 $ and $f_8 $ decrease in the direction of $\vec v$.
\end{lemma}
\begin{proof}
We aim to apply Lemma \ref{flipside}. Therefore, we need to show that the set $\{\nabla f_7 (\xs),\nabla f_8 (\xs), \vec w\}$ is linearly independent, where $\vec w=\langle 1,1,1,1,1,1,1,1\rangle$. It is enough to show that the matrix below
\begin{equation}\label{case2mat}
\left[\begin{array}{c}
       \nabla f_7 \\
       \nabla f_8 \\
       \vec w
       \end{array}\right]=
       \left[\begin{array}{cccccccc}
                                                0 & 0 & 0 & \dis{\frac{\partial f_7 }{\partial x_4}} & \dis{\frac{\partial f_7 }{\partial x_5}} & \dis{\frac{\partial f_7 }{\partial x_6}} & \dis{\frac{\partial f_7 }{\partial x_7}}  & 0\\
                                                \dis{\frac{\partial f_8 }{\partial x_1}} & \dis{\frac{\partial f_8 }{\partial x_2}} & \dis{\frac{\partial f_8 }{\partial x_3}} & 0 & 0 & 0 & 0 & \dis{\frac{\partial f_8 }{\partial x_8}} \\
                                                 1& 1 & 1 & 1 & 1 & 1 & 1 & 1
\end{array}\right]
\end{equation}
has full rank at any $\xs\in\dseven$ which satisfies the hypotheses of the lemma. We have the followings
\begin{itemize}
\item[(1)] $\dis{\frac{\partial f_7 }{\partial x_4}}\bigg|_{\xs}=-\frac{\sigma(x_7^*)}{(\Sigma_2(\xs))^2}$, $\dis{\frac{\partial f_7 }{\partial x_4}\bigg|_{\xs}=\frac{\partial f_7 }{\partial x_5}\bigg|_{\xs}}$, $\dis{\frac{\partial f_7 }{\partial x_5}\bigg|_{\xs}=\frac{\partial f_7 }{\partial x_6}\bigg|_{\xs}}$,
\item[(2)] $\dis{\frac{\partial f_8 }{\partial x_1}}\bigg|_{\xs}=-\frac{\sigma(x_8^*)}{(\Sigma_1(\xs))^2}$, $\dis{\frac{\partial f_8 }{\partial x_1}\bigg|_{\xs}=\frac{\partial f_8 }{\partial x_2}\bigg|_{\xs}}$, $\dis{\frac{\partial f_8 }{\partial x_2}\bigg|_{\xs}=\frac{\partial f_8 }{\partial x_3}\bigg|_{\xs}}$,
\item[(3)] $\dis{\frac{\partial f_7 }{\partial x_7}\bigg|_{\xs}=-\frac{\sigma(\Sigma_2(\xs))}{(x_7^*)^2}}\neq 0$, and $\dis{\frac{\partial f_8 }{\partial x_8}\bigg|_{\xs}=-\frac{\sigma(\Sigma_1(\xs))}{(x_8^*)^2}}\neq 0$.
\end{itemize}
Let $A=\dis{(f_7 )_4(\xs)}$, $B=\dis{(f_7 )_7(\xs)}$, $C=\dis{(f_8 )_1(\xs)}$, and $D=\dis{(f_8 )_8(\xs)}$. We apply two row operations: first $R_2\mapsto -D\cdot R_3+R_2$, then $R_1\mapsto (B/D)\cdot R_2+R_1$. Then we obtain the row equivalent matrix
\[
       \left[\begin{array}{cccccccc}
                                                 \dis{\frac{B}{D}(C-D)}& \dis{\frac{B}{D}(C-D)} &  \dis{\frac{B}{D}(C-D)} & A-B & A-B & A-B & 0  & 0\\
                                               C-D & C-D & C-D & -D & -D & -D & -D & 0 \\
                                                 1& 1 & 1 & 1 & 1 & 1 & 1 & 1
\end{array}\right].
\]
The matrix above has full rank if and only if $A-B\neq 0$ or $C-D\neq 0$, where
\begin{equation*}
A-B=\dis{\frac{(\Sigma_2(\xs)-x_7^*)(1-\Sigma_2(\xs)-x_7^*)}{(\Sigma_2(\xs))^2(x_7^*)^2}}, \ C-D=\dis{\frac{(\Sigma_1(\xs)-x_8^*)(1-\Sigma_1(\xs)-x_8^*)}{(\Sigma_1(\xs))^2(x_8^*)^2}}.
\end{equation*}
Equivalently, it has full rank at $\xs\in\dseven$ unless $x_7^*=\Sigma_2(\xs)$ and $x_8^*=\Sigma_1(\xs)$.

At any $\xs\in\dseven$ such that $x_7^*=x_4^*+x_5^*+x_6^*$, $x_8^*=x_1^*+x_2^*+x_3^*$ and $f_7 (\xs)=f_8 (\xs)$, we derive $(\sigma(x_7^*))^2=(\sigma(x_8^*))^2$, which implies that $x_7^*=x_8^*$. Using $\sum_{n=1}^{8}x_n^*=1$, we find that $\Sigma_1(\xs)=1/4$, $\Sigma_2(\xs)=1/4$, $x_7^*=1/4$ and $x_8^*=1/4$. Since we have $\Sigma_2(\xs)=1/4$, $x_5^*>0$ and $x_6^*>0$, we get $x_4^*<1/4$, which implies that $f _4(\xs)=\sigma(\Sigma_1(\xs))\sigma(x_4^*)>f _8(\xs)=\sigma(\Sigma_1(\xs))\sigma(x_8^*)$, a contradiction. As a result, the matrix in (\ref{case2mat}) has full rank. By Lemma \ref{flipside}, the conclusion of the lemma follows.
\end{proof}
\begin{lemma}\label{lem4.5}
Let $\alphas=\inf_{\tb{x}\in\dseven}\max\{f_1(\tb{x}),\dots,f_8(\tb{x})\}$, where $f_i$ for $i\in I^{\dagger}$ are as in Proposition \ref{dispfunc}. For any $\xs\in\dseven$ such that $f_i (\xs)=f_j (\xs)$ for every $i,j\in I_1 \cup I_2$, there exists a vector $\vec v\in T_{\xs}\dseven$ such that $f_i $ decreases in the direction of $\vec v$ for each $i\in I_1 \cup I_2=\{1,2,3,4,5,6\}$.
\end{lemma}
\begin{proof}
Using the identity $x_8=1-\sum_{n=1}^7x_n$ we rewrite the formulas of $f_1 ,f_2 ,\dots, f_6 $:
\[
\begin{array}{lll}
f_1(\tb{x})=\sigma(\Sigma_2(\tb{x}))\sigma(x_1), & f_2(\tb{x})=\dis{\frac{\sigma(x_2)}{\sigma(\Sigma_1(\tb{x}))}}, & f_3(\tb{x})=\dis{\frac{\sigma(x_3)}{\sigma(\Sigma_2(\tb{x}))}},\\
 f_4(\tb{x})=\sigma(\Sigma_1(\tb{x}))\sigma(x_4), & f_5(\tb{x})=\dis{\frac{\sigma(x_5)}{\sigma(\Sigma_2(\tb{x}))}}, & f_6(\tb{x})=\dis{\frac{\sigma(x_6)}{\sigma(\Sigma_1(\tb{x}))}}.
\end{array}
\] 
All of these functions are smooth on an open neighborhood of $\dseven$ in $\mathbb{R}^8$ and do not depend on $x_7$ or $x_8$.  Since $f_2 (\xs)=f_6(\xs)$ and $f_3 (\xs)=f_5 (\xs)$, we get $x_2^*=x_6^*$ and $x_3^*=x_5^*$. Then using $f_1 (\xs)=f_4 (\xs)$, we obtain $x_1^*=x_4^*$. As a consequence, in particular, we find $\Sigma_1(\xs)=\Sigma_2(\xs)$. Using this fact together with $f_2 (\xs)=f_3 (\xs)$ yields $x_2^*=x_3^*$.  Let $\Sigma_1^*$ denote $\Sigma_1(\xs)$.

We aim to apply Lemma \ref{flipside}. To this purpose, we need to show that the matrix below
\begin{equation}\label{case3mat}
\left[\begin{array}{c}
       \nabla f_1 \\
       \nabla f_2 \\
       \nabla f_3 \\
       \nabla f_4 \\
       \nabla f_5 \\
       \nabla f_6 \\
       \vec w
       \end{array}\right]=
       \left[\begin{array}{cccccccc}
                                                \dis{\frac{\partial f_1 }{\partial x_1}} & 0 & 0 & \dis{\frac{\partial f_1 }{\partial x_4}} & \dis{\frac{\partial f_1 }{\partial x_5}} & \dis{\frac{\partial f_1 }{\partial x_6}} & 0  & 0\\
                                                 \dis{\frac{\partial f_2 }{\partial x_1}} & \dis{\frac{\partial f_2 }{\partial x_2}} & \dis{\frac{\partial f_2 }{\partial x_3}} & 0 & 0 & 0 & 0  & 0\\
                                                0 & 0 & \dis{\frac{\partial f_3 }{\partial x_3}} & \dis{\frac{\partial f_3 }{\partial x_4}} & \dis{\frac{\partial f_3 }{\partial x_5}} & \dis{\frac{\partial f_3 }{\partial x_6}} & 0  & 0\\
                                                 \dis{\frac{\partial f_4 }{\partial x_1}} & \dis{\frac{\partial f_4 }{\partial x_2}} & \dis{\frac{\partial f_4 }{\partial x_3}} & \dis{\frac{\partial f_4 }{\partial x_4}} & 0  & 0 & 0 & 0 \\
                                                0 & 0 & 0 & \dis{\frac{\partial f_5 }{\partial x_4}} & \dis{\frac{\partial f_5 }{\partial x_5}} & \dis{\frac{\partial f_5 }{\partial x_6}} & 0  & 0\\
                                                \dis{\frac{\partial f_6 }{\partial x_1}} & \dis{\frac{\partial f_6 }{\partial x_2}} & \dis{\frac{\partial f_6 }{\partial x_3}} & 0 & 0 & \dis{\frac{\partial f_6 }{\partial x_6}} & 0 & 0 \\
                                                1 & 1 & 1 & 1 & 1 & 1 & 1 & 1
\end{array}\right]
\end{equation}
has full rank at $\xs\in\dseven$. Using the equality $f_1 (\xs)=f_2 (\xs)$, the coordinates of $\xs$ and the definitions of $f_1 ,f_2 ,\dots, f_6 $ , we calculate that
\begin{itemize}
\item[(1)] $\dis{\frac{\partial f_1 }{\partial x_1}\bigg|_{\xs}=-\frac{\sigma(\Sigma_1^*)}{(x_1^*)^2}}$, $\dis{\frac{\partial f_2 }{\partial x_2}\bigg|_{\xs}=\frac{x_2^*(1-x_2^*)-\Sigma_1^*(1-\Sigma_1^*)}{(1-\Sigma_1^*)^2(x_2^*)^2}}$, $\dis{\frac{\partial f_1 }{\partial x_4}\bigg|_{\xs}=-\frac{\sigma(x_1^*)}{(\Sigma_1^*)^2}}$,
\item[(2)]  $\dis{\frac{\partial f_3 }{\partial x_3}\bigg|_{\xs}=-\frac{1}{\sigma(\Sigma_1^*)(x_2^*)^2}}$, $\dis{\frac{\partial f_2 }{\partial x_1}\bigg|_{\xs}=-\frac{\partial f_1 }{\partial x_4}\bigg|_{\xs}}$, $\dis{\frac{\partial f_4 }{\partial x_4}\bigg|_{\xs}=\frac{\partial f_1 }{\partial x_1}\bigg|_{\xs}}$,
\item[(3)]  $\dis{\frac{\partial f_2 }{\partial x_2}\bigg|_{\xs}=\frac{\partial f_5 }{\partial x_5}\bigg|_{\xs}}$, $\dis{\frac{\partial f_3 }{\partial x_3}\bigg|_{\xs}=\frac{\partial f_6 }{\partial x_6}\bigg|_{\xs}}$, $\dis{\frac{\partial f_2 }{\partial x_1}\bigg|_{\xs}=\frac{\partial f_2 }{\partial x_3}\bigg|_{\xs}}$, $\dis{\frac{\partial f_5 }{\partial x_4}\bigg|_{\xs}=\frac{\partial f_5 }{\partial x_6}\bigg|_{\xs}}$,
\item[(4)]   $\dis{\frac{\partial f_2 }{\partial x_3}\bigg|_{\xs}=\frac{\partial f_3 }{\partial x_j}\bigg|_{\xs}}$,  $\dis{\frac{\partial f_5 }{\partial x_6}\bigg|_{\xs}=\frac{\partial f_6 }{\partial x_i}\bigg|_{\xs}}$, $\dis{\frac{\partial f_4 }{\partial x_i}\bigg|_{\xs}=\frac{\partial f_1 }{\partial x_j}\bigg|_{\xs}}$, and $\dis{\frac{\partial f_5 }{\partial x_j}\bigg|_{\xs}=\frac{\partial f_6 }{\partial x_i}\bigg|_{\xs}}$
\end{itemize}
for every $i\in I_1 $, and $j\in I_2 $. Let $A=\sigma(x_1^*)$, $B=\sigma(x_2^*)$, $C=\sigma(\Sigma_1^*)$, $A'=\sigma'(x_1^*)$, $B'=\sigma'(x_2^*)$, and $C'=\sigma'(\Sigma_1^*)$. Note that $A'\neq 0$, $B'\neq 0$ and $C\neq 0$.

We perform simultaneously the following elementary row operations
$R_1\mapsto R_1+R_3,$  $R_4\mapsto R_4+R_2$, $R_5\mapsto (-1)R_3+R_5$, $R_6\mapsto (-1)R_2+R_6$, $R_3\mapsto (ACC')/B'\cdot R_6+R_3$, $R_3\mapsto (ACC')/B'\cdot R_5+R_3$, $R_3\mapsto (AC')/(A'C)\cdot R_4+R_3$,  $R_2\mapsto (AC')/(A'C)\cdot R_1+R_2$
in the matrix in (\ref{case3mat}) to obtain the matrix
\[
       \left[\begin{array}{c|cc|ccccc}
                                                A'C & 0        & \dis{\frac{B'}{C}} & 0 & 0 & 0 & 0  & 0\\ \hline
                                               0 & \dis{\frac{B'}{C}}-AC' & \dis{\frac{AC'B'}{A'C^2}}-AC' & 0 & 0 & 0 & 0  & 0\\
                                                0 & \dis{\frac{AC'B'}{A'C^2}}-AC' & \dis{\frac{B'}{C}}-AC' & 0 & 0 & 0 & 0  & 0\\\hline
                                                 0 & \dis{\frac{B'}{C}} & 0 & A'C & 0  & 0 & 0 & 0 \\
                                                0 & 0 & -\dis{\frac{B'}{C}} & 0 & \dis{\frac{B'}{C}} & 0 & 0  & 0\\
                                                0 & -\dis{\frac{B'}{C}} & 0 & 0 & 0 & \dis{\frac{B'}{C}} & 0 & 0 \\
                                                1 & 1 & 1 & 1 & 1 & 1 & 1 & 1
\end{array}\right].
\]
If we mark each partition of the matrix above by $(i,j)$ reading the first entry from top to bottom and the second entry from left to right  for $i=1,2,3$ and $j=1,2,3$, then the matrix in (\ref{case3mat}) has full rank if and only if $(2,2)$ partition has full rank. Its determinant factors as
\begin{equation}\label{determ}
\frac{B'}{C}\left(1-\frac{AC'}{A'C}\right)\left(\frac{B'}{C}\left(1+\frac{AC'}{A'C}\right)-2AC'\right),
\end{equation}
where
\[
\frac{B'}{C}\left(1-\frac{AC'}{A'C}\right)=\frac{-\Sigma_1^*}{(1-\Sigma_1^*)(x_2^*)^2}\left(1-\frac{x_1^*(1-x_1^*)}{\Sigma_1^*(1-\Sigma_1^*)}\right)\neq 0,
\]
which follows from the fact that the function $t\mapsto t(1-t)$ is increasing on $(0,1/2)$ and the inequality $0<x_1^*<\Sigma_1^*<1/2$, an implication of the equality $2\Sigma_1^*+x_7^*+x_8^*=1$.

Let us assume that the determinant in (\ref{determ}) vanishes. Then we must have the expression
\[
\frac{-\Sigma_1^*}{(1-\Sigma_1^*)(x_2^*)^2}\left(1+\frac{x_1^*(1-x_1^*)}{\Sigma_1^*(1-\Sigma_1^*)}\right)+\frac{2(1-x_1^*)}{x_1^*(\Sigma_1^*)^2} = 0.
\]
Simplifying the left-hand summand, applying the identity $f_1 (\xs)=f_2 (\xs)$ to the right, and finding a common denominator yield:

\[
\frac{-\Sigma_1^*(1-\Sigma_1^*)-x_1^*(1-x_1^*)+2x_2^*(1-x_2^*)}{(x_2^*)^2(1-\Sigma_1^*)^2}=0.
\]
We use the fact that $\Sigma_1^*=x_1^*+2x_2^*$ in the equality above. Then it simplifies to $(x_2^*)^2+2x_1^*x_2^*-x_1^*(1-x_1^*)=0$. The solutions to this quadratic are $x_2^*=-x_1^*\pm\sqrt{x_1^*}$. Since $x_2^*>0$, we get $x_2^*=-x_1^*+\sqrt{x_1^*}$. Using this formula we find that $x_2^*<x_1^*$ if and only if $x_1^*>1/4$.

Since $\sigma(x)=1/x-1$ decreases in $x$ and $\Sigma_1^*<1/2$, we get $\sigma(\Sigma_1^*)>1$. The identity $f_1 (\xs)=f_2 (\xs)$ implies that $\sigma(x_2^*)=(\sigma(\Sigma_1^*))^2\sigma(x_1^*)>\sigma(x_1^*)$, which in turn gives that $x_2^*<x_1^*$. Then we derive that $\Sigma_1^*=x_1^*+2x_2^*=2\sqrt{x_1^*}-x_1^*>3/4$, a contradiction. Hence, the matrix in (\ref{case3mat}) has full rank. By Lemma \ref{flipside}, there is a direction in which each $f_i $ decreases for  $i\in I_1 \cup I_2 $.
\end{proof}
\begin{lemma}\label{lem4.9}
Let $\alphas=\inf_{\tb{x}\in\dseven}\max\{f_1(\tb{x}),\dots,f_8(\tb{x})\}$, where $f_i$ for $i\in I^{\dagger}$ are as in Proposition \ref{dispfunc}. At any $\xs\in\dseven$ satisfying $F^{\dagger}(\xs)=\alphas$ and $f_i (\xs)=f_j (\xs)>f_k (\xs)$ for every $i,j\in I_1 \cup I_3=\{1,2,3,7,8\} $ and $k\in I_2=\{4,5,6\} $, there exists a vector $\vec v\in T_{\xs}\dseven$ such that $f_i $ decreases in the direction of $\vec v$ for each $i\in I_1 \cup I_3$.
\end{lemma}
\begin{proof}
Define the function $\Sigma_{1,3}(\tb{x})=\Sigma_1(\tb{x})+\Sigma_3(\tb{x})$ for $\tb{x}\in\dseven$. We use the identity $x_4=1-\sum_{n=1,n\neq 4}^8x_n$ to rewrite $f_i $ as follows:
\[
\begin{array}{lll}
\dis{f_1 (\tb{x})=\frac{\sigma(x_1)}{\sigma(\Sigma_{1,3}(\tb{x}))}},  &  \dis{f_2 (\tb{x})=\frac{\sigma(x_2)}{\sigma(\Sigma_1(\tb{x}))}}, & \dis{f_3 (\tb{x})=\sigma(\Sigma_{1,3}(\tb{x}))\cdot\sigma(x_3)},\\
\dis{f_4 (\tb{x})=\frac{\sigma(\Sigma_1)}{\sigma(x_5+x_6+\Sigma_{1,3}(\tb{x}))}},  & \dis{f_5 (\tb{x})=\sigma(\Sigma_{1,3}(\tb{x}))\sigma(x_5)},              &  \dis{f_6 (\tb{x})=\frac{\sigma(x_6)}{\sigma(\Sigma_1(\tb{x}))}}, \\
\dis{f_7 (\tb{x})=\frac{\sigma(x_7)}{\sigma(\Sigma_{1,3}(\tb{x}))}},  & \dis{f_8 (\tb{x})=\sigma(\Sigma_1(\tb{x}))\sigma(x_8)}, &
\end{array}
\]
which are smooth on an open neighborhood of $\dseven$ in $\mathbb{R}^8$. The functions $f_1$, $f_2$, $f_3$, $f_7$ and $f_8$ do not depend on $x_4$, $x_5$ or $x_6$. We will show that the matrix
\begin{equation}\label{case4mat}
\left[\begin{array}{c}
       \nabla f_1 \\

       \nabla f_2 \\
       \nabla f_3 \\
       \nabla f_7 \\
       \nabla f_8 \\
       \vec w
       \end{array}\right]=
       \left[\begin{array}{cccccccc}
                                                \dis{\frac{\partial f_1 }{\partial x_1}} & \dis{\frac{\partial f_1 }{\partial x_2}} & \dis{\frac{\partial f_1 }{\partial x_3}} & 0 & 0 & 0 & \dis{\frac{\partial f_1 }{\partial x_7}}  & \dis{\frac{\partial f_1 }{\partial x_8}}\\
                                                 \dis{\frac{\partial f_2 }{\partial x_1}} & \dis{\frac{\partial f_2 }{\partial x_2}} & \dis{\frac{\partial f_2 }{\partial x_3}} & 0 & 0 & 0 & 0  & 0\\
                                                 \dis{\frac{\partial f_3 }{\partial x_1}} & \dis{\frac{\partial f_3 }{\partial x_2}} & \dis{\frac{\partial f_3 }{\partial x_3}} & 0 & 0 & 0 &\dis{\frac{\partial f_3 }{\partial x_7}}   & \dis{\frac{\partial f_3 }{\partial x_8}}\\
                                                 \dis{\frac{\partial f_7 }{\partial x_1}} & \dis{\frac{\partial f_7 }{\partial x_2}} & \dis{\frac{\partial f_7 }{\partial x_3}} & 0 & 0  & 0 & \dis{\frac{\partial f_7 }{\partial x_7}} & \dis{\frac{\partial f_7 }{\partial x_8}} \\
                                                \dis{\frac{\partial f_8 }{\partial x_1}} & \dis{\frac{\partial f_8 }{\partial x_2}} & \dis{\frac{\partial f_8 }{\partial x_3}} & 0 & 0  & 0 & 0 & \dis{\frac{\partial f_8 }{\partial x_8}} \\
                                                1 & 1 & 1 & 1 & 1 & 1 & 1 & 1
\end{array}\right]
\end{equation}
has full rank at $\xs\in\dseven$. Using the coordinates of $\xs$, the definitions of $f_1 $, $f_2 $, $f_3 $, $f_7 $, $f_8 $ and the equalities  $f_1 (\xs)=f_3 (\xs)$, $f_3 (\xs)=f_7 (\xs)$ and $f_2 (\xs)=f_8 (\xs)$, we find the followings:
\begin{itemize}
\item[(1)] $\dis{\frac{\partial f_1 }{\partial x_1}\bigg|_{\xs}=\frac{x_1^*(1-x_1^*)-\Sigma_{1,3}(\xs)(1-\Sigma_{1,3}(\xs))}{(1-\Sigma_{1,3}(\xs))^2(x_1^*)^2}}$, $\dis{\frac{\partial f_1 }{\partial x_2}\bigg|_{\xs}=\frac{\sigma(x_1^*)}{(1-\Sigma_{1,3}(\xs))^2}}$,

\item[(2)] $\dis{\frac{\partial f_2 }{\partial x_1}\bigg|_{\xs}=\frac{\sigma(x_2^*)}{(1-\Sigma_1(\xs))^2}}$, $\dis{\frac{\partial f_2 }{\partial x_2}\bigg|_{\xs}=\frac{x_2^*(1-x_2^*)-\Sigma_1(\xs)(1-\Sigma_1(\xs))}{(1-\Sigma_1(\xs))^2(x_2^*)^2}}$

\item[(3)] $\dis{\frac{\partial f_3 }{\partial x_1}\bigg|_{\xs}=-\frac{\sigma(x_3^*)}{(\Sigma_{1,3}(\xs))^2}}$, $\dis{\frac{\partial f_3 }{\partial x_3}\bigg|_{\xs}=\frac{-x_3^*(1-x_3^*)-\Sigma_{1,3}(\xs)(1-\Sigma_{1,3}(\xs))}{(\Sigma_{1,3}(\xs))^2(x_3^*)^2}}$,
\item[(4)] $\dis{\frac{\partial f_8 }{\partial x_8}\bigg|_{\xs}=-\frac{\sigma(\Sigma_1(\xs))}{(x_8^*)^2}}$, $\dis{\frac{\partial f_8 }{\partial x_1}\bigg|_{\xs}=-\frac{\partial f_2 }{\partial x_1}\bigg|_{\xs}}$, $\dis{\frac{\partial f_3 }{\partial x_1}\bigg|_{\xs}=-\frac{\partial f_1 }{\partial x_3}\bigg|_{\xs}}$,
\item[(5)] $\dis{\frac{\partial f_3 }{\partial x_7}\bigg|_{\xs}=-\frac{\partial f_7 }{\partial x_3}\bigg|_{\xs}}$, $\dis{\frac{\partial f_2 }{\partial x_1}\bigg|_{\xs}=\frac{\partial f_2 }{\partial x_3}\bigg|_{\xs}}$, $\dis{\frac{\partial f_1 }{\partial x_1}\bigg|_{\xs}=\frac{\partial f_7 }{\partial x_7}\bigg|_{\xs}}$, $\dis{\frac{\partial f_8 }{\partial x_1}\bigg|_{\xs}=\frac{\partial f_8 }{\partial x_i}\bigg|_{\xs}}$,
\item[(6)] $\dis{\frac{\partial f_1 }{\partial x_2}\bigg|_{\xs}=\frac{\partial f_1 }{\partial x_j}\bigg|_{\xs}}$, $\dis{\frac{\partial f_3 }{\partial x_1}\bigg|_{\xs}=\frac{\partial f_3 }{\partial x_k}\bigg|_{\xs}}$ and $\dis{\frac{\partial f_7 }{\partial x_1}\bigg|_{\xs}=\frac{\partial f_7 }{\partial x_l}\bigg|_{\xs}}$
\end{itemize}
for every $i=2,3$, $j=3,7,8$, $k=2,7,8$, and $l=2,3,8$. Let $A=(f_1 )_1(\xs)$,  $B=(f_1 )_2(\xs)$, $C=(f_2 )_1(\xs)$, $D=(f_2 )_2(\xs)$, $E=(f_3 )_3(\xs)$, and $F=(f_8 )_8(\xs)$. Note that $B\neq 0$, $B-A\neq 0$, and $D-C\neq 0$. We also have $E+B\neq 0$, which follows from the equality $f_1 (\xs)=f_3 (\xs)$.

We simultaneously apply the row operations: $R_1\mapsto -A\cdot R_6+ R_1$, $R_2\mapsto R_5+R_2$, $R_3\mapsto R_4+R_3$, $R_4\mapsto -B\cdot R_6+R_4$, $R_5\mapsto C\cdot R_6+R_5$, $R_1\mapsto -R_4+R_1$, $R_1\mapsto (1/(B-A))R_1$, $R_2\mapsto (C-D)R_1+R_2$, $R_2\mapsto (1/(C-D))R_2$, $R_3\mapsto -(E+B)R_2+R_3$, $R_3\mapsto -(1/(E+B))R_3$, $R_5\mapsto (C/B)R_4+R_5$ and $R_4\mapsto B\cdot R_3+R_4$ to the matrix in (\ref{case4mat}) to obtain the matrix
\begin{equation}\label{case5mat}
       \left[\begin{array}{cccccc|cc}
                                                0 & 1 & 1 & 1 & 1 & 1 & 2 & 1 \\
                                                 0  & 0 & 1 & 1 & 1 & 1 & 2  & \dis{\frac{F+C-D}{C-D}}\\
                                                 0 & 0 & 0 & 1 & 1 & 1 & \dis{\frac{B-A}{E+B}+2} & \dis{\frac{F+C-D}{C-D}}\\ \hline
                                                0 & 0 & 0 & 0 & 0  & 0 & B\left(\dis{\frac{B-A}{E+B}+2}\right)+A-B & B\left(\dis{\frac{F+C-D}{C-D}}\right) \\
                                                0 & 0 & 0 & 0 & 0  & 0 & \dis{\frac{CA}{B}} & F+C \\ \hline
                                                1 & 1 & 1 & 1 & 1 & 1 & 1 & 1
\end{array}\right].
\end{equation}
The matrix above has full rank if and only if $(2,2)$ partition has non--zero determinant, where
\[
B\left(\dis{\frac{B-A}{E+B}+2}\right)+A-B=\frac{2x_1^*(1-x_1^*)-x_3^*(1-x_3^*)-\Sigma_{1,3}(\xs)(1-\Sigma_{1,3}(\xs))}{(1-\Sigma_{1,3}(\xs))^2(x_3^*)^2},
\]
\[
F+C=\frac{x_8^*(1-x_8^*)-\Sigma_{1}(\xs)(1-\Sigma_{1}(\xs))}{(\Sigma_{1}(\xs))^2(x_8^*)^2},\ \tnr{ and }\  C-D=\frac{\Sigma_1(\xs)}{(x_2^*)^2(1-\Sigma_1(\xs))^2}
\]
calculated by using the facts $f_1 (\xs)=f_3 (\xs)$ and $f_2 (\xs)=f_8 (\xs)$.  Assume that the determinant of $(2,2)$ partition of the matrix in (\ref{case5mat}) vanishes. Then the equality below 
\begin{equation}\label{determ2}
\left(B\cdot\dis{\frac{B-A}{E+B}}+A+B\right)(F+C)=CA\left(\dis{\frac{F+C-D}{C-D}}\right)
\end{equation}
must hold at any point $\xs\in\dseven$ satisfying the hypotheses of the lemma. Let $\Sigma_1^*$, $\Sigma_2^*$ and $\Sigma_{1,3}^*$ denote $\Sigma_1(\xs)$, $\Sigma_2(\xs)$ and $\Sigma_{1,3}(\xs)$, respectively.

Since $f_1 (\xs)=f_7 (\xs)$, we get $x_1^*=x_7^*$.  Using the inequalities $f_2 (\xs)>f_6 (\xs)$, $f_3 (\xs)>f_5 (\xs)$ and $f_8 (\xs)>f_4 (\xs)$, we derive that $x_2^*<x_6^*$, $x_3^*<x_5^*$ and $x_8^*<x_4^*$, which implies $\Sigma_1(\xs)<1/2$. Because otherwise we find $x_5^*+x_6^*+x_1^*>1/2$ contradicting with the fact that $\xs\in\dseven$.  

Note that $2x_8^*+x_2^*+x_3^*<x_2^*+x_3^*+x_4^*+x_8^*<1$ and $x_1^*+x_7^*+x_2^*+x_3^*=2x_1^*+x_2^*+x_3^*<1$. So we have $x_1^*,x_8^*\in (0,(1-x_2^*-x_3^*)/2)$. From the inequalities $x_2^*<x_6^*$, $x_3^*<x_5^*$, and $x_8^*<x_4^*$, we obtain  
\[f_1(\xs)=\frac{1-x_4^*-x_5^*-x_6^*}{x_4^*+x_5^*+x_6^*}\cdot\frac{1-x_1^*}{x_1^*}<\frac{1-x_8^*-x_2^*-x_3^*}{x_8^*+x_2^*+x_3^*}\cdot\frac{1-x_1^*}{x_1^*}.\]
By the equality $f_1(\xs)=f_8(\xs)$, we get $\sigma(x_8^*+x_2^*+x_3^*)\sigma(x_1^*)>\sigma(x_1^*+x_2^*+x_3^*)\sigma(x_8^*)$. Since the function $\sigma(x)/\sigma(x+x_2^*+x_3^*)$ is decreasing over the interval $(0,(1-x_2^*-x_3^*)/2)$, we find $x_1^*<x_8^*$. By the facts $\Sigma_1(\xs)<1/2$ and $(\sigma(\Sigma_1(\xs)))^2\sigma(x_8^*)=\sigma(x_2^*)$, which follows from the rearranging of the equality $f_2(\xs)=f_8(\xs)$, we also find that $x_2^*<x_8^*$.

By using the equality $f_2 (\xs)=f_8 (\xs)$, we simplify the right hand side of the equality above to 
$(x_1^*-\Sigma_{1,3}^*)\left(1-x_1^*-\Sigma_{1,3}^*\right)(x_8^*-x_2^*)\left(1-x_8^*-x_2^*\right)$,
which is nonzero because, $x_2^*<x_8^*$ and $1-x_1^*-\Sigma_{1,3}^*>0$ by the inequality $f_1 (\xs)>f_4 (\xs)$. 

Similarly, by using the equality $f_1 (\xs)=f_3 (\xs)$, we reduce the left hand side of (\ref{determ2}) to 
$\left(x_8^*-\Sigma_1^*\right)\left(1-x_8^*-\Sigma_1^*\right)((x_1^*-x_3^*)(1-x_1^*-x_3^*)+(x_1^*-\Sigma_{1,3}^*)(1-x_1^*-\Sigma_{1,3}^*))$.
We first distribute the factor $\left(x_8^*-\Sigma_1^*\right)\left(1-x_8^*-\Sigma_1^*\right)$ and move the second summand in the resulting expression to the right hand side of the equation in (\ref{determ2}).  On the right hand side of (\ref{determ2}) the term $(x_1^*-\Sigma_{1,3}^*)(1-x_1^*-\Sigma_{1,3}^*)$ is a common factor. We factor this term and, after simplifications, we obtain the following expression
\begin{equation}\label{determ5}
-(x_8^*+\Sigma_1^*)(1-x_1^*-\Sigma_{1,3}^*)(x_1^*+x_3^*)(1-x_2^*-\Sigma_{1}^*)<0.
\end{equation}
On the left hand side of (\ref{determ2}), we have 
\begin{equation}\label{determ6}
(x_8^*-\Sigma_1^*)(1-x_8^*-\Sigma_1^*)(x_1^*-x_3^*)(1-x_1^*-x_3^*).
\end{equation}
Since we assume that the expressions in (\ref{determ5}) and (\ref{determ6}) are equal, there are two cases to consider: 
\[ \textnormal{(1) $\Sigma_1^*<x_8^*$ and $x_1^*<x_3^*$},\quad\textnormal{(2)\ $\Sigma_1^*>x_8^*$ and $x_1^*>x_3^*$}. \]
Assume that (1) is the case. Note that $\Sigma_2^*>1/2$ because, $x_1^*<x_3^*$ and $f_1(\xs)=f_3(\xs)$. We claim that $\Sigma_1^*<1/4<x_8^*$. If $1/4\leq \Sigma_1^*<x_8^*$ holds, then $\sum_{n=1}^8x_n^*>1$, a contradiction. 

If $\Sigma_1^*<x_8^*\leq 1/4$, we see that $\sigma(x_8^*)\geq 3$. We get $3\sigma(\Sigma_1^*)\leq f_8 (\xs)=\alphas$, which implies $3/(\alphas+3)\leq \Sigma_1^*$. By Lemma \ref{lem4.1}, we know that $9\leq\alphas\leq 5+3\sqrt{2}$. As a consequence, we derive  $(24-9\sqrt{2})/46\leq \Sigma_1^*<x_8^*$. Since we have $\Sigma_2^*>1/2$, we find that $(71-18\sqrt{2})/46\leq\Sigma_1^*+\Sigma_2^*+x_8^*$, which in turn gives $x_7^*\leq(-25+18\sqrt{2})/46$. For a lower bound for $\Sigma_2^*$, we solve the inequality $\sigma(\Sigma_2^*)(49+36\sqrt{2})<f_7 (\xs)=\alphas$. Using this lower bound, we conclude that the sum $\Sigma_1^*+\Sigma_2^*=(873-221\sqrt{2})/483>1$, a contradiction. The claim follows.

Upon setting $f_8 (\xs)=\alphas$, substituting $1-\sum_{n=1}^7x_n^*$ for $x_8^*$ and collecting powers of $\Sigma_1^*$ in the resulting formula, we see that $\Sigma_1^*$ is a root of the quadratic
\[
q(x)=(\alphas-1)x^2+(1-\alphas)(1-\Sigma_2^*-x_7^*)x+(\Sigma_2^*+x_7^*).
\]
By the formula for the addition of the roots, we find that $x_8^*$ is the other root. Since  $\Sigma_1^*\neq x_8^*$, the discriminant of the quadratic above is strictly positive.  By the inequality $\Sigma_1^*<1/4<x_8^*$, we find $q(1/4)<0$, which implies that $\Sigma_2^*+x_7^*<(21+18\sqrt{2})/92$.

By substituting $1-\sum_{n=1,n\neq 7}^8x_n^*$ for $x_7^*$ in the equality $f_7 (\xs)=\alphas$ and collecting powers of $\Sigma_2^*$ , we find that $\Sigma_2^*$ is a root of the quadratic
\[
Q(x)=(\alphas-1)x^2+(1-\alphas)(1-\Sigma_1^*-x_8^*)x+(\Sigma_1^*+x_8^*).
\]
By the formula for the addition of the roots, we derive that $x_7^*$ is the other root. We obtain $Q(1/2)<0$ because, we have $x_1^*=x_7^*<1/2<\Sigma_2^*$. Then, we calculate that $\Sigma_1^*+x_8^*<(1+\sqrt{2})/6$, which, in turn, implies that $\sum_{n=1}^8x_n^*<1$, a contradiction. Hence the matrix in (\ref{case5mat}) has full rank in the case (1).

If (2) is the case, then consider $q_{\alpha}(x)=(\alpha-1)x^2+(1-\alpha)(1-\Sigma_2^*-x_7^*)x+(\Sigma_2^*+x_7^*)$ and $Q_{\alpha}(x)=(\alpha-1)x^2+(1-\alpha)(1-\Sigma_1^*-x_8^*)x+(\Sigma_1^*+x_8^*)$ for $\alpha\in[9,5+3\sqrt{2}]$. Note that $q_{\alphas}(x)=q(x)$ and $Q_{\alphas}(x)=Q(x)$. The solutions of $q_{\alpha}(x)+Q_{\alpha}(x)=0$ are
\[
x_{+}(\alpha)=\frac{1}{4}+\frac{1}{4}\sqrt{\frac{\alpha-9}{\alpha-1}}\ \ \tnr{or}\ \
x_{-}(\alpha)=\frac{1}{4}-\frac{1}{4}\sqrt{\frac{\alpha-9}{\alpha-1}}.
\]
Using $x_8^*>x_4^*$, $x_2^*>x_6^*$, $x_3^*>x_5^*$ and $x_1^*<x_8^*$, we derive that $x_7^*<\Sigma_1^*<\Sigma_2^*$. Then we obtain $q_{\alphas}(\Sigma_1^*)+Q_{\alphas}(\Sigma_1^*)=Q_{\alphas}(\Sigma_1^*)<0$, and $q_{\alphas}(x_8^*)+Q_{\alphas}(x_8^*)=Q_{\alphas}(x_8^*)<0$, which implies that $x_-(\alphas)<x_8^*<\Sigma_1^*<x_+(\alphas)$.  We shall use the previous inequality to produce lower and upper bounds for each of the factors in (\ref{determ5}) and (\ref{determ6}).

Since $x_-(\alpha)$ is decreasing over $[9,5+3\sqrt{2}]$, we have $a_1=(\sqrt{2}-1)/2<x_8^*$. By the assumption $x_8^*<\Sigma_1^*$ and the fact $9\leq\alphas$, we find $x_8^*<b_1=1/4$. Otherwise, we would compute that $f_8(\xs)=\alphas<9$, a contradiction. Using the inequality $\alphas=f_8 (\xs)>\sigma(\Sigma_1^*)\sigma(b_1)$, we also get $a_2<\Sigma_1^*<b_2$, where $a_2=3/(8+3\sqrt{2})$ and $b_2=x_+(5+3\sqrt{2})=(2-\sqrt{2})/2$ as $x_+(\alpha)$ is an increasing function on $[9,5+3\sqrt{2}]$. 

We find a lower and an upper bound for $\Sigma_2^*$ as follows: From the assumption $x_1^*>x_3^*$ and the equality $f_1(\xs)=f_3(\xs)$, we have $\Sigma_2^*<b_3=1/2$. By the inequalities $f_4(\xs)<\alphas$, $f_5(\xs)<\alphas$, and $f_6(\xs)<\alphas$, we obtain 
\[
\begin{array}{cc}
x_4^*+x_6^*> \dis{\frac{1-\Sigma_1^*}{(\alphas-1)\Sigma_1^*+1}}+\dis{\frac{\Sigma_1^*}{\alphas(1-\Sigma_1^*)+\Sigma_1^*}}, & x_5^*>\dis{\frac{\Sigma_2^*}{\alphas(1-\Sigma_2^*)+\Sigma_2^*}}.
\end{array}
\]
The expression on the right hand side of the first inequality above is decreasing both in $\alphas$ and $\Sigma_1^*<1/2$. So we find $x_4^*+x_6^*>1/4$ by using the bounds $5+3\sqrt{2}$ and $b_2$ for $\alphas$ and $\Sigma_1^*$, respectively. Then we have $\Sigma_2^*>x_5^*+1/4$. 

Since the expression on the right hand side of the second inequality above is decreasing in $\alphas$ and increasing in $\Sigma_2^*$, by substituting the bounds $5+3\sqrt{2}$ and $x_5^*+1/4$ and rearranging, we get $
(4+3\sqrt{2})(x_5^*)^2-(3/4)(4+3\sqrt{2})x_5^*+1/4<0$. Thus $x_5^*$ is greater than the smaller root $(3-2\sqrt{2})/4$ of the left hand side quadratic. Then it follows that 
$(2-\sqrt{2})/2<\Sigma_2^*$. Next we will consider the following two cases:
\[ \textnormal{(2a) $\frac{2-\sqrt{2}}{2}<\Sigma_2^*<\frac{1}{3}$},\quad\textnormal{(2b)\ $\frac{1}{3}\leq\Sigma_2^*<\frac{1}{2}$}. \]

Assume that (2b) is the case. By rearranging the equalities $f_1(\xs)=\alphas$, $f_2(\xs)=\alphas$ and $f_3(\xs)=\alphas$, we derive 
\begin{equation}\label{new22}
x_1^*=\dis{\frac{1-\Sigma_2^*}{(\alphas-1)\Sigma_2^*+1}},\quad  x_2^*=\dis{\frac{\Sigma_1^*}{\alphas(1-\Sigma_1^*)+\Sigma_1^*}},\quad x_3^*=\dis{\frac{\Sigma_2^*}{\alphas(1-\Sigma_2^*)+\Sigma_2^*}}.
\end{equation}
The right hand side of the expression for $x_2^*$ is increasing in $\Sigma_1^*$ and decreasing in $\alphas$. Therefore we find $x_2^*<b_4=(9\sqrt{2}-10)/62$ by substituting the relevant bounds $9$ and $b_2$ for $\alphas$ and $\Sigma_1^*$, respectively. We also find $x_1^*\leq b_5=2/11$ by plugging $9$ and $1/3=a_3\leq\Sigma_2^*$ because, the expression on the right hand side of the equality for $x_1^*$ above is decreasing in both $\Sigma_2^*$ and $\alphas$. Similarly, since the right hand side of the equality
\begin{equation}\label{new23}
x_1^*+x_3^*=\dis{\frac{1-\Sigma_2^*}{(\alphas-1)\Sigma_2^*+1}}+\dis{\frac{\Sigma_2^*}{\alphas(1-\Sigma_2^*)+\Sigma_2^*}}
\end{equation}
is decreasing in both $\alphas$ and $\Sigma_2^*$, we get $a_4=(2-\sqrt{2})/3<x_1^*+x_3^*$ by substituting $5+3\sqrt{2}$ and $b_3=1/2$ for $\alphas$ and $\Sigma_2^*$, respectively.  The right hand side of the expression 
\begin{equation}\label{new24}
x_1^*-x_3^*=\dis{\frac{1-\Sigma_2^*}{(\alphas-1)\Sigma_2^*+1}}-\dis{\frac{\Sigma_2^*}{\alphas(1-\Sigma_2^*)+\Sigma_2^*}}
\end{equation}
is again decreasing in both $\alphas$ and $\Sigma_2^*$. So we obtain $x_1^*-x_3^*<b_6=27/209$
by plugging in $9$ and $a_3=1/3$ for $\alphas$ and $\Sigma_2^*$, respectively. As a result we have the following bounds:
\[
\begin{array}{lllll}
a_1> 0.20710, & a_2> 0.24504, & a_3\geq 0.33333, & a_4>0.19526, & b_1= 0.25000, \\
b_2< 0.29290, & b_3= 0.50000, & b_4< 0.04400, & b_5< 0.18182, & b_6<0.12919. 
\end{array}
\]
Then we compute that $(a_1+a_2)a_4(a_3-b_5)(1-b_4-b_2)> 0.00886$, which implies that the expression in (\ref{determ5}) is less than $-0.00886$. Similarly, we also calculate that $(b_2-a_1)b_6(1-a_1-a_2)(1-a_4)<0.00489$, which shows that the expression in (\ref{determ6})
is greater than $-0.00489$. Hence, the determinant of the (2,2) partition of the matrix in (\ref{case5mat}) cannot be $0$. 

Assume that the inequality in (2a) holds. In this case, we have $\Sigma_2^*<b_3=1/3$. Using the equality in (\ref{new23}) we get $(1027-480\sqrt{2})/1519=a_4<x_1^*+x_3^*$. Since we have $\Sigma_1^*=x_1^*+x_2^*+x_3^*>x_2^*+a_4$, by rearranging the equality for $x_2^*$ in (\ref{new22}), we derive the inequality $(4+3\sqrt{2})(x_2^*)^2-(1-a_4)(4+3\sqrt{2})x_2^*+a_4<0$. Thus $x_2^*$ is greater than the smaller root of the left hand quadratic in the previous inequality. This implies that 
\[
\frac{2424+1698\sqrt{2}-\sqrt{9776852+6468345\sqrt{2}}}{1519 (4+3\sqrt{2})}+\frac{1027-480 \sqrt{2}}{1519}=a_2<\Sigma_1^*.
\]
Substituting the bounds $9$ for $\alphas$ and $a_3=(2-\sqrt{2})/2$ for $\Sigma_2^*$ in the expression for $x_1^*$ in (\ref{new22}) we obtain $x_1^*<b_5=1/(9\sqrt{2}-8)$. Similarly, using the previous bounds for $\alphas$ and $\Sigma_2^*$ in (\ref{new24}), we get $x_1^*-x_3^*<b_6=(369-81 \sqrt{2})/1519$. As a result we have 
\[
\begin{array}{lllll}
a_1> 0.20710, & a_2> 0.26716, & a_3>0.29289, & a_4>0.22921, & b_1= 0.25000, \\
b_2< 0.29290, & b_3= 0.33333, & b_4< 0.04400, & b_5< 0.21151, & b_6<0.16752. 
\end{array}
\]
Using these estimates we calculate $(a_1+a_2)a_4(a_3-b_5)(1-b_4-b_2)> 0.00586$ and $(b_2-a_1)b_6(1-a_1-a_2)(1-a_4)<0.00583$. Hence, the determinant of the (2,2) partition of the matrix in (\ref{case5mat}) cannot be $0$ in this case as well. Finally by Lemma \ref{flipside}, we obtain the conclusion of the lemma.
\end{proof}
\begin{proof}[Proof of Proposition \ref{lem4.2}]
It follows from Lemmas \ref{lem4.1-4}, \ref{lem4.3}, \ref{lem4.4}, \ref{lem4.5} and \ref{lem4.9}.
\end{proof}

We use Proposition \ref{lem4.2} to calculate the infimum of $G^{\dagger}$ over the simplex $\dseven$. In particular, we prove Theorem \ref{thm4.1}. First, we establish the following:
\begin{theorem}\label{thm4.14}
Let $F^{\dagger}\co\dseven\to\mathbb{R}$ be the function defined by $\tb{x}\mapsto\max\{f_i (\tb{x}): i\in I^{\dagger}\}$, where $f_i$ are defined as in Proposition \ref{dispfunc}. Then $\inf_{\tb{x}\in\Delta^7}F^{\dagger}(\tb{x})=5+3\sqrt{2}$.
\end{theorem}
\begin{proof}
By Proposition \ref{lem4.2}, we know that $\tb{x}^*\in\Delta_7$, i.e., $f_i (\xs)=\alphas$ for every $i\in I^{\dagger}$. Using the identities $f_4 (\xs)=f_8 (\xs)$, $f_1 (\xs)=f_7 (\xs)$, $f_2 (\xs)=f_6 (\xs)$, and $f_3 (\xs)=f_5 (\xs)$, we get $x_4^*=x_8^*$, $x_1^*=x_7^*$, $x_6^*=x_2^*$, and $x_3^*=x_5^*$. By the fact $f_4 (\xs)=f_1 (\xs)$,  we obtain $(x_1^*-x_4^*)(x_2^*+x_3^*)(1-x_1^*-x_2^*-x_3^*-x_4^*) = 0,$ or $x_1^*=x_4^*$. The last equality, in turn, gives that $\Sigma_1^*=\Sigma_2^*$. 

By the equality $f_2 (\xs)=f_3 (\xs)$, we see that $x_2^*=x_3^*$. 
Since $\sum_{n=1}^8x_n^*=1$, we obtain $x_2^*=1/4-x_1^*$. Using the equality $f_1 (\xs)=f_2 (\xs)$,  
we find  
\[
x_2^*-3x_2^*x_1^*-(x_1^*)^2x_2^*-4(x_2^*)^2+4(x_2^*)^2x_1^*+4(x_2^*)^3-(x_1^*)^3 =  0,
\]
which 
simplifies to $1-4x_1^*-4(x_1^*)^2=0$. 
The solutions are $x_1^*=(-1\pm\sqrt{2})/2$.
Since $x_1^*>0$, we get
$x_2^*=(3-2\sqrt{2})/4$.
In particular, we conclude that $x_i^*=x_1^*$ for every $i\in\{4,7,8\}$ and $x_j^*=x_2^*$ for every for $x_j\in\{3,5,6\}$.
Finally, we calculate that $\inf_{\tb{x}\in\Delta^7}F^{\dagger}(\tb{x})=\sigma((2-\sqrt{2})/2)\sigma((\sqrt{2}-1)/2)=5+3\sqrt{2}$. 
\end{proof}
\begin{proof}[Proof of Theorem \ref{thm4.1}]
By the definitions of $G^{\dagger}$ and $F^{\dagger}$, we have $G^{\dagger}(\tb{x})\geq F^{\dagger}(\tb{x})$ for every $\tb{x}\in\dseven$. A direct computation shows that  $G^{\dagger}(\xs)=F^{\dagger}(\xs)$. Then the conclusion of the theorem follows. 
\end{proof}

\subsection{On the uniqueness of $\xs$ in Theorem \ref{thm4.14}}

It is worth to emphasize the similarities between the statements (\ref{a}) and (\ref{b}) listed in \S\ref{S2} and the following two statements
\begin{enumerate}[label=\arabic*]
\renewcommand{\labelenumi}{\alph{enumi}}
\item[c.] $\inf_{\tb{x}\in\Delta^7}F^{\dagger}(\tb{x})=\min_{\tb{x}\in\Delta^7}F^{\dagger}(\tb{x})$,
\item[d.] There exists $\tb{x}^*\in\Delta_7\subset\Delta^7$ such that $\min_{\tb{x}\in\Delta^7}F^{\dagger}(\tb{x})=f_i (\tb{x}^*)$ for
$i\in I^{\dagger}$,
\end{enumerate}
used in the proof of Theorem \ref{thm4.14} to calculate the number $(1/2)\log(5+3\sqrt{2})$.
Although it is straight forward to observe the fact in (\ref{b}), it takes considerable amount of calculations to prove the statement given in (d). Analogous
to Lemma \ref{lem2.1}, Theorem \ref{thm4.14} shows that the point $\tb{x}^*$ is unique. Assuming the uniqueness of the point $\tb{x}^*$ a priori together
with (c) suggests an alternative way of finding the coordinates of the point $\tb{x}^*\in\dseven$.

Let $T_1$ be the transformation defined in (\ref{transf}). Since we have $f_i (T_1(\tb{x}))=f_j (\tb{x})$ for every
$\tb{x}\in\dseven$ for every pair $(i,j)\in\{(1,4),(4,1),(3,5),(5,3),(7,8),(8,7),(2,6),(6,2)\}$, we conclude that $\{f_i \co i\in
I^{\dagger}\}=\{f_i \circ T_1\co i\in I^{\dagger}\}$. Let us define $H_1\co \dseven\to \mathbb{R}$, where $H_1(\tb{x})=\max\{(f_i \circ T_1)(\tb{x})\co i\in I^{\dagger}\}$.
We see that $F^{\dagger}(\tb{x})=H_1(\tb{x})$ for every $\tb{x}\in\dseven$ and $\min_{\tb{x}\in\Delta^7}F^{\dagger}(\tb{x}) = \min_{\tb{x}\in\Delta^7}H_1(\tb{x})$. Since $F^{\dagger}$
takes its minimum value at the point $\tb{x}^*$ and $\{f_i (\tb{x}^*)\co i\in I^{\dagger}\}=\{(f_i \circ T_1)(\tb{x}^*)\co i\in I^{\dagger}\}$, the function $H_1$ takes its minimum value at
the point $T_1^{-1}(\tb{x}^*)$. Then we obtain $T_1^{-1}(\tb{x}^*)=\tb{x}^*$ by the uniqueness of $\xs$. This means that $x_1^*=x_4^*$,
$x_2^*=x_5^*$, $x_3^*=x_6^*$ and $x_7^*=x_8^*$. 

Let $\Delta^3=\{\tb{x}\in\dseven\co x_1=x_4,x_2=x_5,x_3=x_6,x_7=x_8\}$. Note that
 $f_1 (\tb{x})=f_4 (\tb{x})$, $f_2 (\tb{x})=f_5 (\tb{x})$, $f_3 (\tb{x})=f_6 (\tb{x})$ and $f_7 (\tb{x})=f_8 (\tb{x})$ for every $\tb{x}\in\Delta^3$. Define the continuous function $F_1\co\Delta^3\to\mathbb{R}$ such that $\tb{x}\mapsto\max(g_1(\tb{x}),g_2(\tb{x}),g_3(\tb{x}),g_7(\tb{x}))$, where $g_i=f_i |_{\Delta^3}$ for $i=1,2,3,7$. Then we have $\min_{\tb{x}\in\dseven}F^{\dagger}(\tb{x})=\min_{\tb{x}\in\Delta^3}F_1(\tb{x})$. 

Consider $T_2\co\dseven\to\dseven$ defined by $x_2\mapsto x_3$, $x_3\mapsto x_2$ and $x_i\mapsto x_i$ for every $i\in I^{\dagger}-\{2,3\}$.
The map $T_2$ preserves $\dseven$ and $\Delta^3$. Then we have $g_i(T_2(\tb{x}))=f_j(\tb{x})$ for every
$\tb{x}\in\Delta^3$ for every pair $(i,j)\in\{(1,1),(2,3),(3,2),(7,7)\}$. An argument similar to the one above for $H_2\co\dseven\to \mathbb{R}$ defined by
$H_2(\tb{x})=\max\{(f_i \circ T_2)(\tb{x})\co i=1,2,3,7\}$ shows that $T_2^{-1}(\tb{x}^*)=\tb{x}^*$. This means that $x_2^*=x_3^*$.

Let $\Delta^2=\{\tb{x}\in\Delta^3\co x_2=x_3\}$. Note that $g_2(\tb{x})=g_3(\tb{x})$ for every $\tb{x}\in\Delta^2$. Define
the functions $h_i\co \Delta^2\to \mathbb{R}$ such that $h_i=g_i|_{\Delta^2}$ for $i=1,2,7$.
Introduce the continuous function $F_2\co \Delta^2 \to \mathbb{R}$, where $\tb{x}\mapsto\max(h_1(\tb{x}),h_2(\tb{x}),h_7(\tb{x}))$. Then we have $\min_{\tb{x}\in\Delta^3}F_1(\tb{x})=\min_{\tb{x}\in\Delta^2}F_2(\tb{x})$. Note that  $\tb{x}^*=(x_1^*,x_2^*,x_2^*,x_1^*,x_2^*,x_2^*,x_7^*,x_7^*)$ with $x_1^*+2x_2^*+x_7^*=1/2$.

In the rest of the discussion, we will consider $\Delta^2$ as a submanifold of $\mathbb{R}^8$. 
Then the tangent space $T_{\tb{x}}\Delta^2$ at any $\tb{x}\in\Delta^2$ is a subspace of $T_{\tb{x}}\Delta^7$ generated by the vectors $\vec u_1=\langle 1,0,0,1,0,0,-1,-1\rangle$ and $\vec u_2=\langle 0,1,1,0,1,1,-2,-2\rangle$. Note that $h_1(\tb{x})$, $h_2(\tb{x})$, and $h_7(\tb{x})$ are smooth in an open neighborhood of $\Delta^2$. Therefore, $\nabla h_i(\tb{x})\cdot\vec v$ is the derivative of $h_i$ in the direction of $\vec v\in T_{\tb{x}}\Delta^2$ for each $i\in\{1,2,7\}$.

Using the identity $x_7=1/2-x_1-2x_2$, we rewrite the formulas of $h_1(\tb{x})$, $h_2(\tb{x})$, and $h_7(\tb{x})$ as follows: $h_1(\tb{x})=\sigma(x_1+2x_2)\sigma(x_1)$, $h_2(\tb{x})=\sigma(x_2)/\sigma(x_1+2x_2)$, and $h_7(\tb{x})=\sigma(x_1+2x_2)\sigma(1/2-x_1-2x_2)$. Then we find the following partial derivatives:
\begin{equation*}
\frac{\partial h_1}{\partial x_1}=  \frac{2(x_1^2+2x_1x_2-x_1-x_2+2x_2^2)}{x_1^2(x_1+2x_2)^2},\ \ \frac{\partial h_1}{\partial x_2} =  -\frac{2(1-x_1)}{x_1(x_1+2x_2)^2},
\end{equation*}
\begin{equation*}
\frac{\partial h_2}{\partial x_1}=  \frac{1-x_2}{x_2(1-x_1-2x_2)^2},\ \ \ \tnr{and}\ \ \ \frac{\partial h_2}{\partial x_2} = \frac{4x_1x_2+2x_2^2-x_1+x_1^2}{x_2^2(1-x_1-2x_2)^2}.
\end{equation*}
It is clear that $h_1$ and $h_2$
have no critical points in $\Delta^2$. On the other hand, every point on the set
$L=\{\tb{x}\in \Delta^2\co -1+4x_1+8x_2=0\}$ is a critical point for the function $h_7$. Because, we have
\begin{equation*}
\frac{\partial h_7}{\partial x_1}=  \frac{-1+4x_1+8x_2}{(x_1+2x_2)^2(-1+2x_1+4x_2)^2},\ \ \ \tnr{and}\ \ \ \frac{\partial h_7}{\partial x_2} = \frac{2(-1+4x_1+8x_2)}{(x_1+2x_2)^2(-1+2x_1+4x_2)^2}.
\end{equation*}
Let $\tilde{h}_i=h_i|_{L}$ for $i=1,2,7$. Then we get
\begin{displaymath}
\tilde{h}_1(x_2)  =  3\cdot\frac{3+8x_2}{1-8x_2},\ \ \tilde{h}_2(x_2)  = \frac{1}{3}\cdot\frac{1-x_2}{x_2},\ \ \tnr{and}\ \ \tilde{h}_7(x_2)  =  9
\end{displaymath}
for $0<x_2<1/8$. Let $I=(0,1/8)$. Then $\tilde{h}_1$ and $\tilde{h}_2$ have no critical points in $I$. We see that $\tilde{h}_1(x_2)>\tilde{h}_7(x_2)$ for every $x_2\in I$, because we have
\begin{displaymath}
\lim_{x_2\to  0^+}3\cdot\frac{3+8x_2}{1-8x_2}=9^+\ \tr{and}\ \lim_{x_2\to 1/8^-}3\cdot\frac{3+8x_2}{1-8x_2}=\infty.
\end{displaymath}
Therefore, it is enough to calculate the infimum of the maximum of $\tilde{h}_1$ and $\tilde{h}_2$ over $I$ to calculate the infimum of the
maximum of $h_1$, $h_2$ and $h_7$ over $L$. 

Since $\tilde{h}_1$ and $\tilde{h}_2$ have no critical points in $I$, the
infimum of the maximum of $\tilde{h}_1$ and $\tilde{h}_2$ is attained at a point $x_2^*$ such that $\tilde{h}_1(x_2^*)=\tilde{h}_2(x_2^*)$. In other words, we need to solve the equation $64(x_2^*)^2+36x_2^*-1=0$. We get $x_2^*=(-9\pm\sqrt{97})/32.$ Since $x_2^*$ is positive, we calculate that $\inf_{x_2\in I}\left\{\max\left(\tilde{h}_1(x_2),\tilde{h}_2(x_2)\right)\right\}=(17+2\sqrt{97})/3$. Note that the point $x_2^*$ at which the infimum of the maximum of $\tilde{h}_1$ and $\tilde{h}_2$ over $L$ is unique.

We claim that there exist $i,j\in\{1,2,7\}$ with $i\neq j$ such that $h_i(\xs)=h_j(\xs)$. Assume otherwise that
$h_i(\xs)\neq h_j(\xs)$ for every $i,j\in\{1,2,7\}$ for $i\neq j$. Then we have either $h_1(\xs)> h_i(\xs)$ for $i=2,7$
or $h_2(\xs)>h_i(\xs)$ for $i=1,7$ or $h_7(\xs)>h_i(\xs)$ for $i=1,2$. Since $h_1$ and $h_2$ have no
critical points in $\Delta^2$, we cannot have $h_1(\xs)> h_i(\xs)$ for $i=2,7$ or $h_2(\xs)>h_i(\xs)$ for $i=1,7$.

Assume that $h_7(\xs)>h_i(\xs)$ for $i=1,2$. It follows from the fact that $h_7$ is continuous on $\Delta^2$, there exists a neighborhood $V$ of $\xs$ in $\Delta^2$ so that $h_7(\tb{x})>
h_i(\tb{x})$ for every $\tb{x}\in V$ for $i=1,2$. Hence, we get $F_2(\tb{x})=h_7(\tb{x})$
for every $\tb{x}\in V$. Since we have $h_7(\xs)=F_2(\xs)$, the function $h_7$ has a local minimum at $\xs\in V$. This means that
$\xs\in L$ so that $x_1^*=(13-\sqrt{97})/16$, $x_2^*=(-9+\sqrt{97})/32$ and $x_7^*=1/4$.
This is a contradiction. Because, we know by Lemma \ref{lem4.1} that $F^{\dagger}(\bar{\tb{x}})=5+3\sqrt{2}<h_2(\xs)$.
Hence, there exist $i,j\in\{1,2,7\}$ with $i\neq j$ such that $h_i(\xs)=h_j(\xs)$. We need to consider the cases
\begin{itemize}
\item[(I)] $h_1(\xs)=h_7(\xs)>h_2(\xs)$,\quad (III) $h_1(\xs)=h_2(\xs)>h_7(\xs)$,
\item[(II)] $h_7(\xs)=h_2(\xs)>h_1(\xs)$,\quad (IV) $h_1(\xs)=h_2(\xs)=h_7(\xs)$.
\end{itemize}

Assume that (I) is the case. 
Let $\vec v_1=\vec u_1+(-1/2)\vec u_2$. The first order partial derivatives of
$h_1$, $h_2$ and $h_7$ show that
$\nabla h_1(\tb{x})\cdot\vec v_1<0$, $\nabla h_7(\tb{x})\cdot\vec v_1=0$,   $\nabla h_2(\tb{x})\cdot\vec v_1>0$
for every $\tb{x}\in\Delta^2$. Therefore, $h_1$ is decreasing, $h_2$ is increasing and $h_7$ is constant along a line segment in the direction of $\vec v_1$.
This means that if we move along the line segment starting at $\xs$ in
the direction of the vector $\vec v_{1}$ for a sufficiently small amount, we obtain a point $\tb{y}\in\Delta^2$ such that $h_i(\tb{y})<\alphas$
for $i=1,2$ and $h_7(\tb{y})=\alphas$. This is a contradiction. Therefore, we cannot have $h_1(\xs)=h_7(\xs)>h_2(\xs)$.

Assume that (II) is the case. 
Let $\vec v_2=-\vec v_1$. The first order partial derivatives of $h_1$,
$h_2$ and $h_7$ show that
$
\nabla h_1(\tb{x})\cdot\vec v>0$, $\nabla h_7(\tb{x})\cdot\vec v=0$,  $\nabla h_2(\tb{x})\cdot\vec v<0$
for every $\tb{x}\in\Delta^2$. An argument similar to the above applies, mutatis mutandis, in this case to show that we cannot have $h_2(\xs)=h_7(\xs)>h_1(\xs)$.
We already know that (III) is not possible at the point $\xs$. Because this case corresponds to Case III in Lemma \ref{lem4.5}. As a result, we have $h_2(\xs)=h_7(\xs)$ and $h_7(\xs)=h_1(\xs)$.

Using the equality $h_1(\xs)=h_2(\xs)$, we see that $\sigma(x_1^*)=\sigma(1/2-x_1^*-2x_2^*)$, which implies $x_1^*=1/4-x_2^*$. We use $h_2(\xs)=h_7(\xs)$ to obtain $16(x_2^*)^2-24x_2^*+1=0$ or $x_2^*=(3\pm 2\sqrt{2})/4$.
Since $x_2^*$ is positive, we find $x_2^*=(3-2\sqrt{2})/4$, $x_1^*=(\sqrt{2}-1)/2$ and $x_1^*=x_7^*$.
Finally, we calculate that 
$\min_{\tb{x}\in\Delta^2}F_2(\tb{x})=5+3\sqrt{2}$.

In the discussion above, we don't refer to the statement $\tb{x}^*\in\Delta_7$ given in (d). The assumption that the point $\xs$ is unique reduces the necessary calculations to obtain $\alphas$ considerably. Notice the fact that $\Gamma_{\mathcal{D}^{\dagger}}$ is not a symmetric decomposition of $\Gamma$. The investigation of the likely conditions such as the convexity properties of the
displacement functions for the decompositions $\Gamma_{\mathcal{D}^1}$ and $\Gamma_{\mathcal{D}^{\dagger}}$ that might lead to a proof of the uniqueness of the points $\tb{x}^*$ will be left to future studies.


\section{Proof of The Main Theorem}\label{sec5}

In this final section, we present a detailed proof of the main result of this paper, stated below:
\begin{theorem}\label{thm5.1}
Let $\xi$ and $\eta$ be non--commuting isometries of $\hyp$. Suppose that $\xi$ and $\eta$ generate a torsion--free discrete group which is not
co--compact and contains no parabolic. Let $\Gamma_{\dagger}$ and $\alpha_{\dagger}$ denote the set of isometries $\{\xi,\eta,\xi\eta\}$ and the real number $5+3\sqrt{2}$, respectively. Then for any $z_0\in\hyp$ we have
$$
e^{\left(\dis{2\max\nolimits_{\gamma\in\Gamma_{\dagger}}\left\{\dgamma\right\}}\right)}\geq\alpha_{\dagger}.
$$
\end{theorem}
\begin{proof}
By Proposition $9.2$ in \cite{CSParadox}, the group $\langle\xi,\eta\rangle$ is a free
group on the generators $\xi$ and $\eta$. Let $z_0$ be a point in $\hyp$. If $\Gamma=\langle\xi,\eta\rangle$ is geometrically infinite, then Theorem \ref{thm3.4} and Lemma \ref{lem1.2} imply that
$$
\max_{\gamma\in\Gamma_{\dagger}}\left\{\dgamma\right\} \geq  \frac{1}{2}\log G^{\dagger}(\tb{m})
                         \geq  \frac{1}{2}\log\left(\inf_{\tb{x}\in\Delta^7} G^{\dagger}(\tb{x})\right)
$$
for $\tb{m}=\left(m_{p(\psi)}\right)_{\psi\in\Psi^{\dagger}}\in\Delta^7$, where $p$ and $m_{p(\psi)}$ are the bijection and the total measures defined in (\ref{sigma}) and \S \ref{Sec3}, respectively. The function $G^{\dagger}$ is defined in Theorem \ref{thm4.1}, which implies the conclusion of the theorem.

Assume that $\Gamma=\langle\xi,\eta\rangle$ is geometrically finite. Then $(\xi,\eta)$ is in $\mathfrak{GF}$, an open subset of the character variety $\mathfrak{X}=\tnr{Isom}^+(\hyp)\times\tnr{Isom}^+(\hyp)$, consisting of $(\xi,\eta)$ such that $\langle\xi,\eta\rangle$ is free, geometrically finite and without any parabolic.  Let $f_{z_0}^{\dagger}\co\mathfrak{X}\to\R$ be the function defined as
\[(\xi,\eta)\mapsto\max\{\dxi,\deta,\dxieta\}.\]
It is straightforward to see that $f^{\dagger}z_0$ is a proper, continuous and non--negative valued function on $\mathfrak{X}$. Therefore, it takes a minimum value at some point $(\xi_0,\eta_0)\in\overline{\mathfrak{GF}}$. We claim that $(\xi_0,\eta_0)$ is in $\overline{\mathfrak{GF}}-\mathfrak{GF}$.

Assume on the contrary that $(\xi_0,\eta_0)$ is in $\mathfrak{GF}$. Since $\xi_0$, $\eta_0$ and
$\xi_0\eta_0$ have infinite orders in $\langle\xi_0,\eta_0\rangle$, we have $\xi_0\cdot z\neq z$, $\eta_0\cdot z\neq z$ and $\xi_0\eta_0\cdot z\neq z$ for every $z\in\hyp$. In particular, we get that $\xi_0\cdot z_0\neq z_0$, $\eta_0\cdot z_0\neq z_0$ and $\xi_0\eta_0\cdot z_0\neq z_0$. Therefore, there
exists hyperbolic geodesic segments joining $z_0$ to $\xi_0\cdot z_0$, $z_0$ to $\eta_0\cdot z_0$ and $z_0$ to $\xi_0\eta_0\cdot z_0$. Note that we have the equalities
$\tr{dist}(z_0,\xi_0\eta_0\cdot z_0)=\tr{dist}(\xi_0^{-1}\cdot z_0,\eta_0\cdot z_0)$ and $\dxi=\dxiinv$. We consider the geodesic triangle
$\Delta=\Delta_{P_2P_0P_1}$, where $P_1=\xi_0^{-1}\cdot z_0$, $P_0=z_0$ and $P_2=\eta_0\cdot z_0$. The value $f_{z_0}^{\dagger}(\xi_0,\eta_0)$ is the longest side length of $\Delta$. There are two cases to consider: $\Delta$ is acute or $\Delta$ is not acute.

Assume that the latter is the case. Let $\gamma$ be the unique longest edge of $\Delta$. By the hyperbolic law of sines, $\gamma$ is opposite to the non--acute angle. If $P_1$ lies in $\gamma$, we let $P_1^{(i)}$ be a sequence of points in the interior of $\gamma$ so that $P_1^{(i)}\to P_1$.
Let $P_j^{(i)}=P_j$ for $j\in\{0,2\}$ and $i\in\mathbb{N}$. Otherwise, we let $P_2^{(i)}$ be a sequence of points in the interior of $\gamma$ so that $P_2^{(i)}\to P_2$ and define $P_j^{(i)}=P_j$ for $j\in\{0,1\}$ and $i\in\mathbb{N}$. Let $\Delta_i$ be the geodesic triangle contained in $\Delta$ with vertices $P_0^{(i)}$, $P_1^{(i)}$ and $P_2^{(i)}$. By the construction, the unique longest side $\gamma_i$ of $\Delta_i$ is contained in $\gamma$ for all but finitely many $i$. Let $\{\xi_i\}$ be a sequence of isometries such that $\xi_i\to\xi$ and $\xi_i^{-1}\cdot z_0=P_1^{(i)}$. Similarly, Let $\{\eta_i\}$ be a sequence of isometries such that $\eta_i\to\eta$ and $\eta_i\cdot z_0=P_2^{(i)}$. Then we have $(\xi_i,\eta_i)\in\mathfrak{GF}$ for all but finitely many $i$ and $f_{z_0}^{\dagger}(\xi_i,\eta_i)=l(\gamma_i)<f_{z_0}^{\dagger}(\xi_0,\eta_0)$, a contradiction.

Assume that $\Delta$ is acute. Then the perpendicular arc $\gamma_1$ from $P_1$ to the geodesic containing $P_0$ and $P_2$ meets it in the interior of the edge of $\Delta$ opposite to $P_1$. Let $P_1^{(i)}$ be a sequence of points in the interior of $\gamma_1$ so that $P_1^{(i)}\to P_1$. For each $i$, we see that
\begin{center}
\includegraphics[scale=.4]{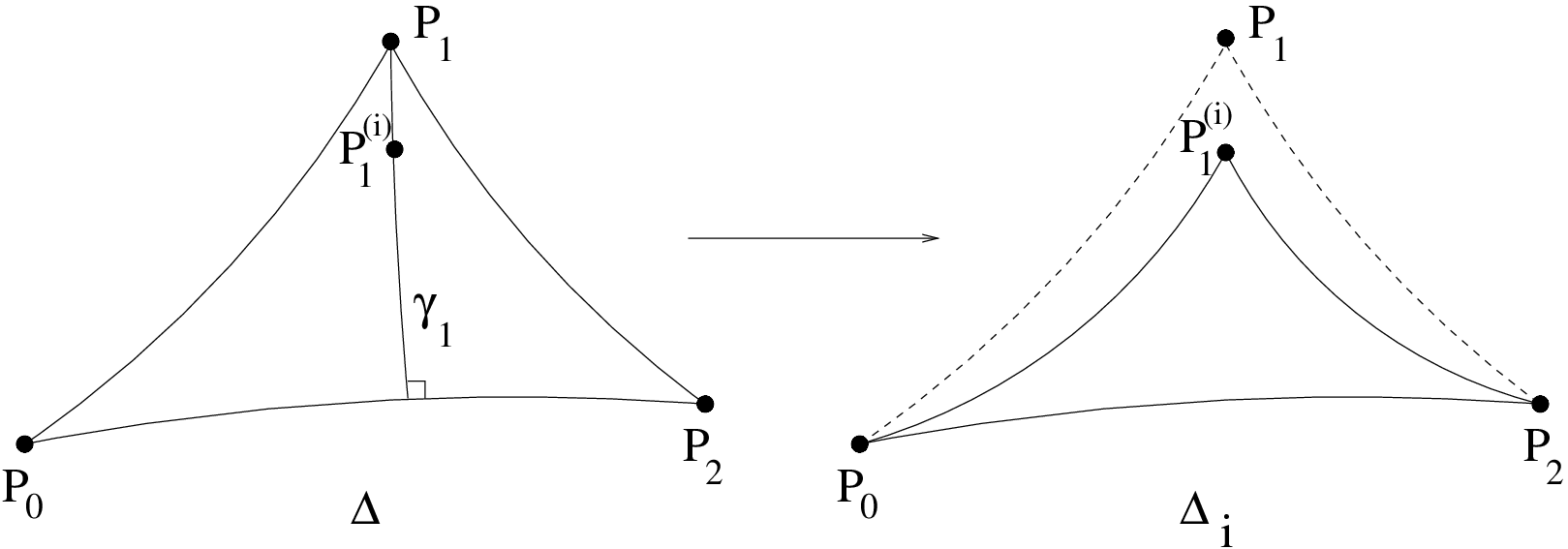}
\end{center}
$d(P_1^{(i)},P_0)<d(P_1,P_0)$ by applying the hyperbolic law of cosines to the right triangle containing $P_1^{(i)}$, $P_0$ and a sub--arc of $\gamma_1$. Similarly, we have $d(P_1^{(i)},P_2)<d(P_1,P_2)$.

The triangle $\Delta_i$ with vertices $P_0$, $P_1^{(i)}$ and $P_2$ is itself acute because, its angles at $P_0$ and $P_2$ are less than those of $\Delta$, and its angle at $P_1^{(i)}$ limits to the angle of $\Delta$ at $P_1$. Thus, the perpendicular arc $\gamma_2^{(i)}$ from $P_2$ to the geodesic containing $P_0$ and $P_1^{(i)}$ meets
\begin{center}
\includegraphics[scale=.4]{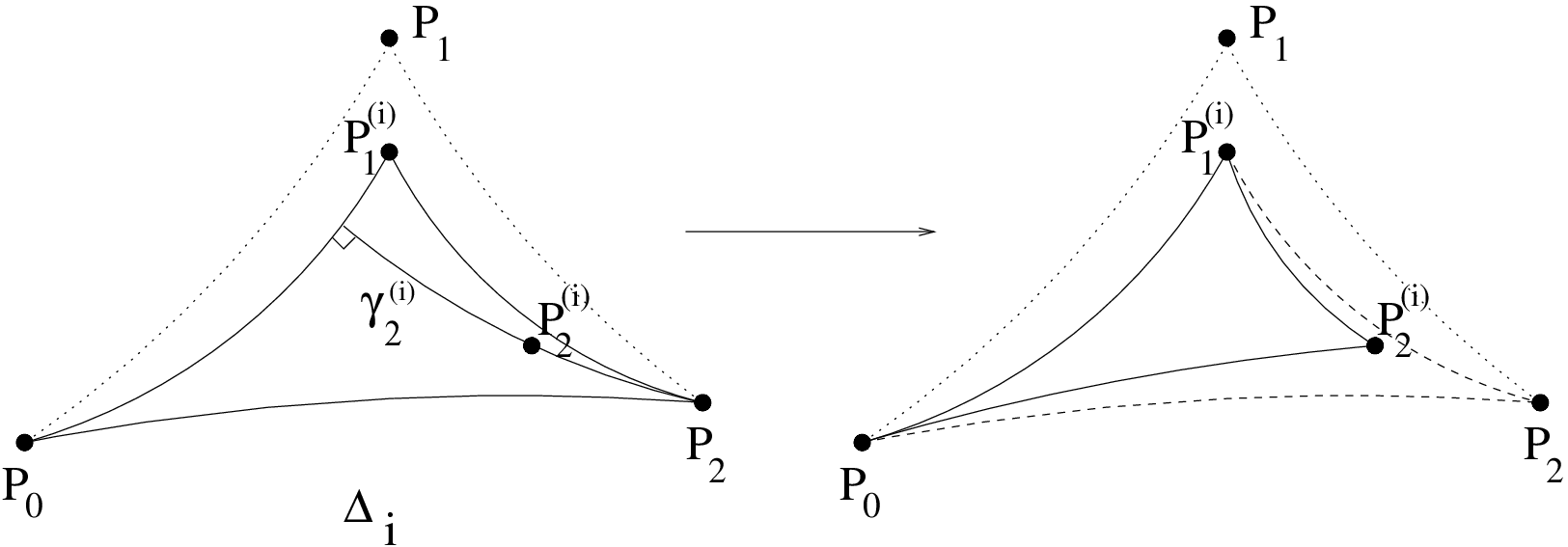}
\end{center}
this geodesic inside of $\Delta_i$. Let $P_2^{(i)}$ be the point on $\gamma_2^{(i)}$ at distance $1/i$ from $P_2$. We find that $d(P_2^{(i)},P_0)<d(P_2,P_0)$ and $d(P_2^{(i)},P_1^{(i)})<d(P_2,P_1^{(i)})<d(P_2,P_1)$ by the hyperbolic law of cosines. In other words, by the two--step process described above, we obtain a triangle with vertices at $P_0$, $P_1^{(i)}$ and $P_2^{(i)}$ so that all edge lengths are less than those of $\Delta$. Let $\{\xi_i\}$ and $\{\eta_i\}$ be the sequences such that $\xi_i^{-1}\cdot z_0=P_1^{(i)}$ and $\eta_i\cdot z_0=P_2^{(i)}$. Then we have $f_{z_0}^{\dagger}(\xi_i,\eta_i)<f_{z_0}^{\dagger}(\xi_0,\eta_0)$ for all but finitely many $i$, a contradiction. Hence, we conclude that $(\xi_0,\eta_0)\in\mathfrak{GF}$.

Finally, the facts that the set of $(\xi,\eta)$ such that $\langle\xi,\eta\rangle$ is free, geometrically infinite and without any parabolic is dense in $\overline{\mathfrak{GF}}-\mathfrak{GF}$ (\cite{CSParadox}, Proposition 8.2) and every $(\xi,\eta)\in\mathfrak{X}$ with $\langle\xi,\eta\rangle$ is free and without any parabolic is in $\overline{\mathfrak{GF}}$ (\cite{CSParadox}, Proposition 9.3) reduce geometrically finite case to geometrically infinite case proving the theorem.
\end{proof} 

%
%
%
%

\end{document}